\documentclass[
a4paper, 
11pt, 
parskip=half,
]{scrartcl}

	\usepackage[colorlinks,plainpages,unicode]{hyperref}
	\hypersetup{urlcolor=gray, citecolor=gray, linkcolor=gray,hypertexnames=false}
	\usepackage {amsmath}
	\usepackage {amssymb}
	\usepackage {amsthm}
	\usepackage {dsfont}
	\usepackage {bbm}
	\usepackage {enumitem}
	\setlist[itemize]{leftmargin=*, partopsep=0pt, itemsep=0pt,parsep=0pt,topsep=0pt}
	\setlist[enumerate]{leftmargin=*, partopsep=0pt, itemsep=0pt,parsep=0pt,topsep=0pt}
	
	\usepackage{amscd}
	\usepackage {soul}
	\usepackage{pgf,tikz}
	\usetikzlibrary{arrows, matrix, decorations.markings, calc}
	\usepackage{standalone}
	\usepackage{wasysym}
	\usepackage{graphicx}
	\usepackage{url}
	\usepackage{float}
	\usepackage[labelfont=bf]{caption}
	\usepackage{geometry}
\usepackage{tocstyle} 

\usepackage[toc]{multitoc}

\usepackage[british]{babel}              	 
\usepackage[T1]{fontenc} 
\usepackage[utf8]{inputenc}        
\usepackage[british]{isodate}
\origdate

	\usepackage{fancyhdr}
 	
\expandafter\let\expandafter\oldproof\csname\string\proof\endcsname
\let\oldendproof\endproof
\renewenvironment{proof}[1][\proofname]{%
  \oldproof[\textbf{#1}]%
}{\oldendproof}

	\newcommand{\R}{\mathbb{R}}
	\renewcommand{\C}{\mathbb{C}}
	
	\newcommand{\N}{\mathbb{N}}
	\newcommand{\norm}[1]{\left \lVert #1 \right\rVert}
	\newcommand{\normw}[1]{\left \lVert #1 \right\rVert_{W^{1,2}_\epsilon}}
	\newcommand{\norml}[1]{\left \lVert #1 \right\rVert_{L^2_\epsilon}}
	\newcommand{\normlinfty}[1]{\left \lVert #1 \right\rVert_{L^\infty_\epsilon}}
	\newcommand{\normlinftyy}[1]{\left \lVert #1 \right\rVert_{L^\infty}}
	\newcommand{\normltwo}[1]{\left \lVert #1 \right\rVert_{L^2}}
	\newcommand{\scal}[2]{\left \langle #1, #2 \right \rangle}
	\newcommand{\geps}[2]{g^0_\epsilon \left( #1, #2 \right)}
	\newcommand{\res}[2]{\left. #1 \right|_{#2}}
	\newcommand{\abs}[1]{\left \lvert #1 \right\rvert}
	\newcommand{\inv}{^{-1}}
	\newcommand{\Y}{L^2(\R,\R^n)}
	\newcommand{\X}{W^{1,2}(\R,\R^n)}
	\newcommand{\W}{W^{2,2}(\R,\R^n)}
	\newcommand{\vect}[2]{\begin{pmatrix} #1 \\ #2 \end{pmatrix}}
	\newcommand{\fracb}{\frac{1}{1- \norm b ^2}}
	\newcommand{\piplus}{\pi_\epsilon^+(z)}
	\newcommand{\piminus}{\pi_\epsilon^-(z)}

	\DeclareMathOperator{\diag}{diag}

	\DeclareMathOperator{\image}{im}
	
	\DeclareMathOperator{\ind}{ind}
	\DeclareMathOperator{\spann}{span}
	
	\DeclareMathOperator{\sgn}{sgn}
	
	\DeclareMathOperator{\energy}{E}

	\DeclareMathOperator{\Crit}{Crit}

\newtheoremstyle{main}
  {}
  {}
  {\itshape}
  {}
  {\bfseries}
  {.}
  { }
  {\thmname{#1}}				%

	\theoremstyle{main}
	\newtheorem {TheoremL}{Main Theorem}
	
	\theoremstyle{plain}
	\newtheorem {Corollary}{Corollary}[section]
	\newtheorem {Theorem}[Corollary]{Theorem}
	\newtheorem {Lemma}[Corollary]{Lemma}
	\newtheorem {Proposition}[Corollary]{Proposition}
	\newtheorem {Definition}[Corollary]{Definition}

	\theoremstyle{definition}
	\newtheorem{Standing}[Corollary]{Assumption}
	\newtheorem{Notation}[Corollary]{Notation}
	\newtheorem {Remark}[Corollary]{Remark}

\title{Gradient Flow Line Near Birth-Death Critical Points}
\date{ \today}
\author{Charel Antony\thanks{The author was partially supported by the Swiss National
Science Foundation (grant number  200021\_156000).}} 

\begin {document}
\maketitle
%
%
%
%
\begin{abstract}
Near a birth-death critical point in a one-parameter family of gradient flows, there are precisely two Morse critical points of index difference one on the birth side. This paper gives a self-contained proof of the folklore theorem that these two critical points are joined by a unique gradient trajectory up to time-shift. The proof is based on the Whitney normal form, a Conley index construction, and an adiabatic limit analysis for an associated fast-slow differential equation. 
\end{abstract}
%
\thispagestyle{empty}
%
%
%
%
\section{Introduction}
Let $M$ be a smooth $n$-manifold and let $(F_\lambda)_{\lambda\in \R}$ be a smooth family of functions having $p_0\in M$ as a birth-death critical point for $\lambda=0$. This means that at the point $p_0$, the differential of $F_0$ vanishes, the Hessian of $F_0$ has a one dimensional kernel, and two numbers,
the third derivative $ d^3 F_{0}$ and the mixed derivative $d(\partial_\lambda F|_{\lambda=0})$ in the direction of the kernel of the Hessian, do not vanish. (See Definition \ref{def:normal1}.) Assume these two latter numbers have opposite signs. Then the birth side is $\lambda>0$, i.e. for $\lambda>0$ sufficiently small, there exist precisely two critical points $p_{\pm}(\lambda)$ of $F_\lambda$ near $p_0$, which are Morse and have index difference one. 

\begin{TheoremL}[Main Theorem]\label{thm:goal1}
Assume $F_\lambda:M\to\R$, $p_0$, $p_{\pm}(\lambda)$ are as above and let $\left\{G_{\lambda}\right\}_{\lambda\in \R}$ be a smooth family of Riemannian metrics on $M$.

\vspace{-\parskip}
Then, for $\lambda>0$ sufficiently small, there exists a unique (up to time-shift) gradient trajectory $\Gamma_\lambda:\R\to M$ solving $\dot \Gamma_\lambda~=~-\nabla_{G_\lambda} F_\lambda~\circ~ \Gamma_\lambda$ and connecting the critical points $p_-(\lambda)$ to $p_+(\lambda)$. Moreover, the unstable manifold of $p_-(\lambda)$ and the stable manifold of $p_+(\lambda)$ intersect transversally along $\Gamma_\lambda$.
\end{TheoremL}
\vspace{-1em}
\begin{proof}
See page  \pageref{proof:goal1}.
\end{proof}
\vspace{-1em}
This is a folklore theorem and a version of it was stated, but not proven in \cite[(3.10)]{Flo1}. The main difficulty lies in working with \textbf{general metrics}, as the metric can couple the system of non-linear ODE's. In \cite[Section 3.6]{Hut1} to circumvent this difficulty, the author chooses a standard Euclidean metric in the chart, which represents the trivial case of our result. Our approach to the proof relies heavily on the Whitney normal form in a local coordinate chart (Appendix~\ref{sec:app}). Another approach was outlined independently in an as yet unpublished paper by Suguru Ishikawa.

The Main Theorem can be seen as a converse to Smale's Cancellation Lemma \cite[5.4]{Mil1}, which asserts that a pair of critical points of index difference one with precisely one (transverse) gradient trajectory joining them can be eliminated. The Main Theorem says that conversely if such a pair of critical points appears in a generic family, they are necessarily joined by exactly one (transverse) gradient trajectory.

Generic families~$F_\lambda$ joining two Morse functions \cite{Cer1} are used to construct continuation morphisms between the Morse homology groups. There are two such constructions, bifurcation and cobordism, both used by A. Floer \cite{Flo1,Flo2}, and a longstanding conjecture asserts that both morphisms agree \cite[Remark 2.1]{Hut2}. In the absence of birth-death critical points a proof of them agreeing on the homology level was given by D. Komani \cite{Kom1}. The Main Theorem will be needed for a proof of this conjecture on the homology level in full generality.

\begin{figure}[h!tb]
\centering
\includegraphics[width=0.9\textwidth]{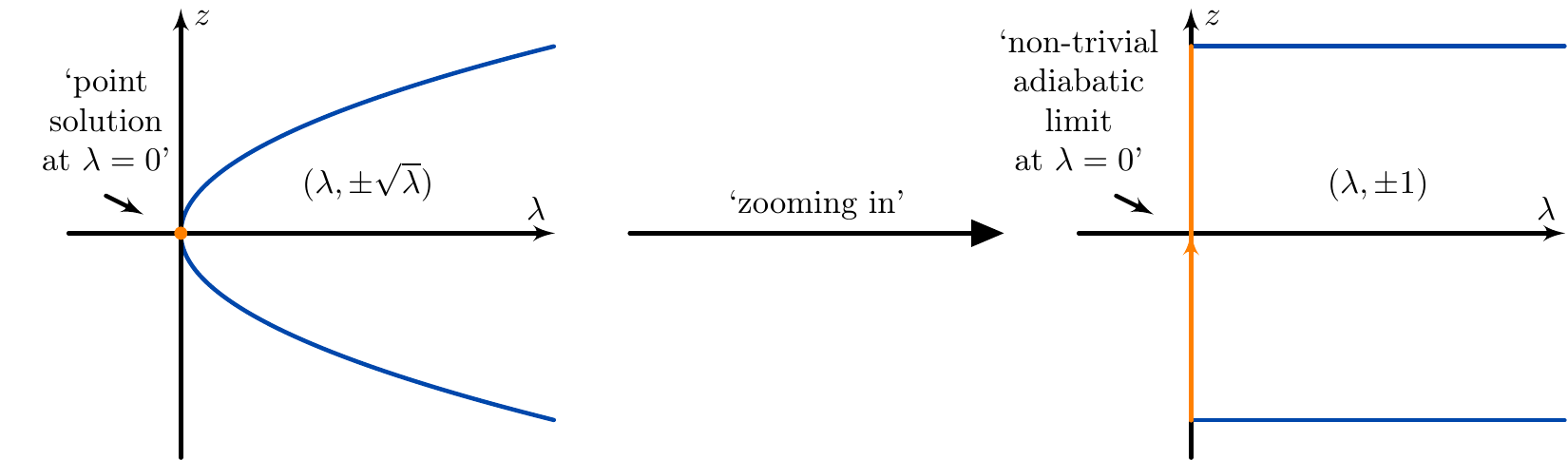}
\caption{Critical points in the $z$ direction before and after the `zooming in'. The point solution at $\lambda=0$ becomes a non-trivial limit solution, as discussed below. }
\label{fig:intro1}
\end{figure}

\textbf{Constant Metric Case. }%
The starting point for the proof of the Main Theorem is Whitney's normal form (Theorem \ref{thm:normal2}). This theorem says that for $\lambda$ small in a smooth $\lambda$-dependent family of local coordinates, $F_\lambda$ is equal to a Morse part in $x$ plus a one dimensional birth-death family in $z$, i.e. in these local coordinates
\[
F_\lambda(x,z) = \frac{1}{2} x^\top A x + \frac{1}{3} z^3 - \lambda z 
\]
for $(x,z) \in \R^{n-1}\times \R$ and where $A\in \R^{(n-1)\times (n-1)}$ is an invertible symmetric matrix. This normal form also requires a $\lambda$-dependent change of coordinates on the target $\R$ which induces a rescaling in time, however this does not affect the statements in the Main Theorem. Now assume the metric is constant in these coordinates and $\lambda>0$. After `zooming in' (Lemma \ref{lem:trans1}), the gradient flow equation in local coordinates with $\epsilon = \sqrt \lambda >0$ becomes a fast-slow system%
\footnote{We denote by $\scal{u}{v} = {u}^\top v$ the standard Euclidean inner product.} 
\begin{equation}\label{eq:intro1}\tag{$\ast_\epsilon$}
\begin{cases}
\epsilon \, \dot x 		&=- A x + b(1-z^2), \\
\dot z 			&=- \scal{Ax}{b}+(1-z^2),
\end{cases}
\end{equation}
\nopagebreak
where $A \in \R^{(n-1)\times (n-1)}$ is a symmetric invertible matrix and $b \in \R^{n-1}$ is a column vector of norm $\norm b <1$. (For general Riemannian metrics see equation \eqref{eq:eq4}.) The right hand side has exactly two zeroes at $x=0$ and $z=\pm 1$. See Figure \ref{fig:intro1}.

Our first step is to see what happens if we set $\epsilon=\sqrt \lambda$ to zero in \eqref{eq:intro1}. The first line becomes a constraint which can be solved for $x$ since $A$ is invertible. So the equation reads
\[\tag{$\ast_0$}
\begin{cases}
x				&= A^{-1} b(1-z^2), \\
\dot z 			&=- \scal{Ax}{b}+(1-z^2).
\end{cases}
\]
Substituting the first equation into the second line yields the differential equation 
\[
\dot z 			= (1-\norm b^2)(1-z^2)
\]
which has the unique solution $\tanh((1-\norm b^2)t)$ (up to time-shift) under the boundary conditions $\lim_{t\to \pm \infty} z(t) = \pm 1$. With the corresponding $x$ from the first line, we get a limit solution $\gamma_0: \R \to \R^n$ at $\epsilon=\sqrt \lambda=0$. 

The name adiabatic limit was introduced to mathematics by E. Witten in \cite{Wit1}. This method can be used to construct solutions to certain geometric differential equations; typical examples are \cite{Don3}, \cite{Sal1}, \cite{Fin1}, \cite{Sal4} and \cite{Sal5}. These equations always present a certain rescaling parameter. The leitmotif is that, as this parameter goes to zero, the underlying metric degenerates and the equations change into those of a different geometric type.  A solution of the limit equations will then give an approximate solution for the original equations for small parameters. In this paper, we showcase this method in the simpler setting of ODE's as opposed to the PDE setting of most other authors.

From the existence part of the argument (described in more details below), we will construct a solution $\gamma_\epsilon$ to \eqref{eq:intro1} close to the limit solution $\gamma_0$ for $\epsilon>0$ small. Here close depends on the right choice of $\epsilon$-weighted norms, which is dictated by the geometry of the problem. These norms are crucial to making the analysis work as we cannot expect that $\gamma_\epsilon$ converges in $C^\infty$ to $\gamma_0$, as $\epsilon$ tends to zero. It is worth noting that for $b=0$ (a~metric adapted to the $(x,z)$ splitting) $\gamma_0 = (0, \tanh)$ solves \eqref{eq:intro1} for all $\epsilon>0$. So this case is trivial. On the other hand, for $b \neq 0$, even this constant metric case is surprisingly hard to solve. The general metric case only adds error terms to \eqref{eq:intro1} which do not add any new hurdles in the proof of the Main Theorem.

\begin{figure}[h!tb]
\centering
\includegraphics[width=0.9\textwidth]{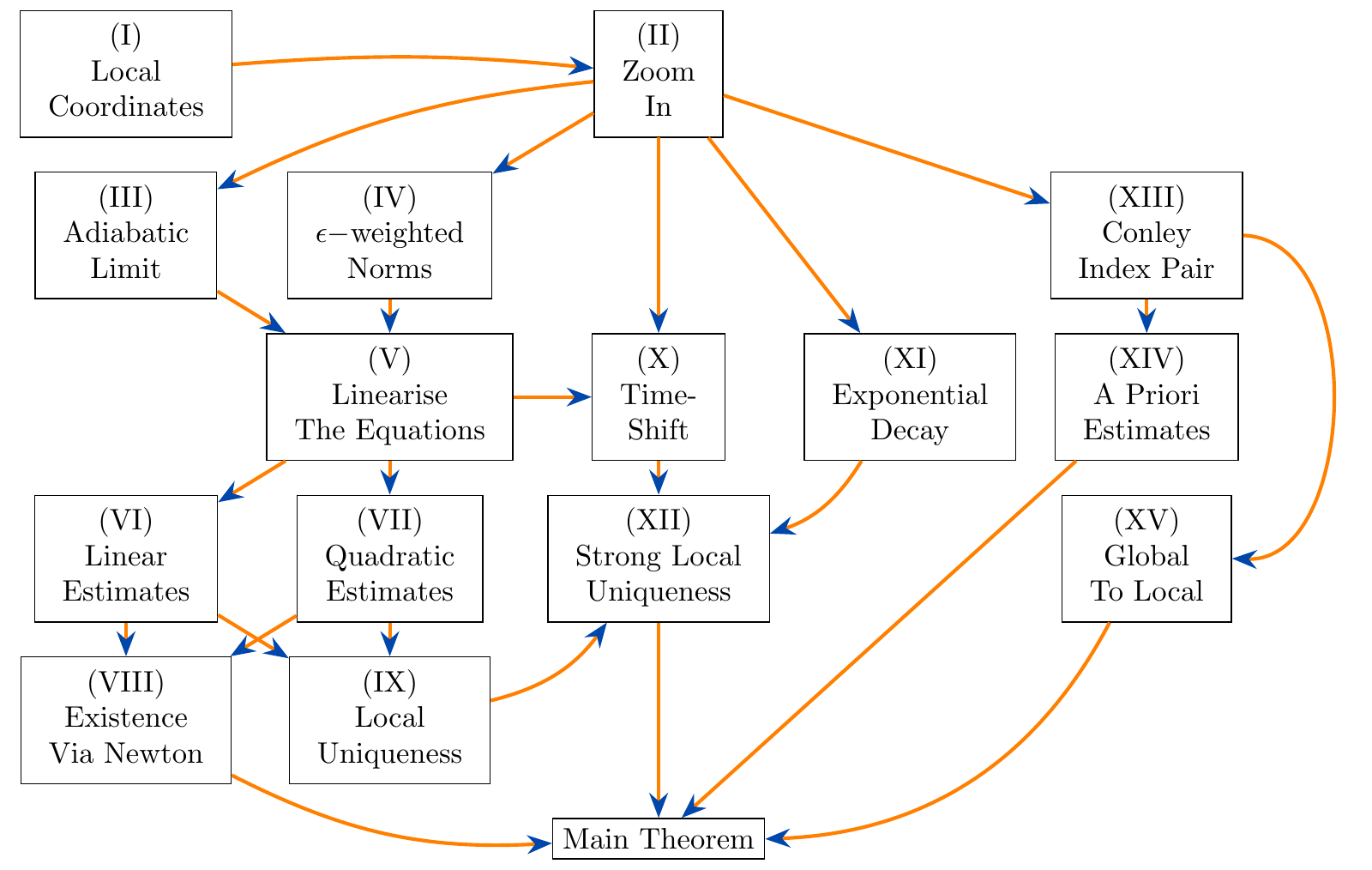}
\caption{Structure of the argument. (I)-(II) in Section \ref{sec:two}, (III)-(IX) in Section \ref{sec:three}, (X)-(XII) in Section \ref{sec:four} and (XIII)-(XV) in Section \ref{sec:five}.}
\label{fig:intro2}
\end{figure}

\textbf{Outline Of Proof. }%
Let us now outline the structure of the proof of the Main Theorem. We provide an overview in Figure \ref{fig:intro2}. The argument can be broadly separated into three parts: local existence/uniqueness, strong local uniqueness and global uniqueness. As stated before, the starting point are local coordinates \textbf{(I)} adapted to both function and metric. After zooming in \textbf{(II)}, we can state a local version of our result (Theorem \ref{thm:goal2}) which will imply the Main Theorem, as described below in (XV).
\begin{itemize}
\item \textbf{Local Existence And Uniqueness. } This is an adiabatic limit argument. We start from a solution $\gamma_0$ of the equations for $\epsilon=0$ \textbf{(III)}, as illustrated in Figure \ref{fig:intro1}, which is interpreted as an approximate solution to the $\epsilon$-equations. The goal is to find a true solution nearby. The key point is the correct choice of $\epsilon$-weighted norms \textbf{(IV)}. This also leads to $\epsilon$-weighted norms on $L^2(\R, \R^n), W^{1,2}(\R, \R^n)$ and $L^\infty(\R, \R^n)$. By linearising the equations \eqref{eq:intro1} at $\gamma_0$ \textbf{(V)}, we obtain a linear Fredholm operator $D_\epsilon$ of index one. This operator is surjective and has a one dimensional kernel, since our equations \eqref{eq:intro1} are invariant under time-shift. The Newton iteration for \eqref{eq:intro1} requires a bijective linearised operator which is obtained by restricting the domain to the range of the $\epsilon$-~adjoint operator~$D_\epsilon^\ast$ \textbf{(VIII)}. To establish convergence of the Newton scheme, we will require linear estimates \textbf{(VI)} for $D_\epsilon$ and $D_\epsilon^\ast$ with norms depending on $\epsilon>0$ and constants independent of $\epsilon$. Quadratic estimates \textbf{(VII)} are also needed. This gives existence of a solution $\gamma_\epsilon$ to \eqref{eq:intro1} near the limit solution $\gamma_0$ for $\epsilon>0$ small. Furthermore, the same estimates imply local uniqueness \textbf{(IX)} under the following constraints.
\begin{enumerate}[label = \textbf{(C\arabic*)}]
\item The solution $\gamma$ is close to $\gamma_0$ in the $\epsilon$-weighted $L^\infty$-norm, and the difference $\gamma-\gamma_0$ satisfies a weak $L^2$ estimate. 
\item The difference $\gamma-\gamma_0$ belongs to the codimension one subspace $\image D_\epsilon^\ast$.
\end{enumerate}

\phantom{'}
\item \textbf{Strong Local Uniqueness.} We strengthen the local uniqueness by eliminating constraints (C1) and (C2). Exponential decay at the ends \textbf{(XI)} will eliminate $(C1)$, and a time-shift will deal with the constraint (C2). These results are combined in \textbf{(XII)} to get the desired strong local uniqueness.  

\phantom{'}
\item \textbf{Global Uniqueness. }  We are left with proving global uniqueness. To this end we construct a Conley index pair \textbf{(XIII)} to go from the global problem on the manifold to a local problem in the charts \textbf{(XV)}.
More precisely, a gradient trajectory in the unstable manifold of $p_-(\lambda)$ that leaves the local coordinate chart where equation~\eqref{eq:intro1} is valid cannot return to $p_+(\lambda)$, because the value of $F_\lambda$ has become too small.

Furthermore, the Conley index pair construction will also provide us with a priori estimates \textbf{(XIV)} on all solutions of our equations \eqref{eq:intro1}. These let us apply the strong uniqueness to finish the proof of the Main Theorem.      
\end{itemize}
\textbf{Acknowledgements:} The author thanks his adviser D.A. Salamon for all the crucial support and is also grateful to S. Trautwein for some helpful comments. 
%
%
%
%
%
\section{The Gradient Equation In Local Coordinates}\label{sec:two}
In this section, we execute steps (I) and (II) as outlines in the introduction. Namely, we write out the gradient equation in adapted local coordinates. Then, we will `zoom in' and derive the local set of equations \eqref{eq:eq4} that we will be used throughout the rest of this paper. Thereafter, we will state a local version of the Main Theorem, Theorem \ref{thm:goal2}. 

In the Main Theorem, we are interested in $F:\R\times M \to \R$ a smooth function having $p_0\in M$ as a birth-death critical point of $(F_\lambda)_{\lambda \in \R}$ for $\lambda=0$ as in Definition \ref{def:normal1} together with a smooth family of Riemannian metrics $(G_\lambda)_{\lambda\in \R}$ on $M$. By Theorem \ref{thm:normal2} and up to flipping $\lambda$ to $-\lambda$, we get
 $\mathfrak c, \epsilon_0>0$, a family of charts $\varphi_\lambda: U\subset M \to \R^n$ with $\varphi_0(p_0)=0$ and a family of affine invertible maps $\chi_{\lambda}:\R \to \R$ such that\footnote{Here $G^{-1}$ stands for the inverse of the matrix associated with the metric.}
\begin{equation}\label{eq:trans1}
\begin{aligned}
\chi_{\lambda} \circ F_{\lambda} \circ \varphi_{\lambda}^{-1}(x,z) 	&=\frac 1 {\mathfrak c}\bigg(\frac{1}{2} x^\top A x + \frac{1}{3} z^3 - \lambda z\bigg),\\
((\varphi_{\lambda})_\ast G_\lambda)_{(x,z)}^{-1} 					&= \mathfrak c\bigg(\begin{pmatrix} \mathbbm 1 & b \\   b^\top & 1 \end{pmatrix} + h_{\lambda,(x,z)}\bigg), 
\end{aligned}
\end{equation}
for all $(x,z)\in (\R^{n-1}\times \R)\cap \varphi_\lambda(U)$ and all $\lambda \in \R$ with $\abs{\lambda} \leq \epsilon_0^2$. Here $A \in \R^{(n-1)\times (n-1)}$ is a symmetric invertible matrix, $\mathbbm 1 \in \R^{(n-1)\times (n-1)}$ is the identity matrix, $b \in \R^{n-1}$ is a column vector with $\norm b <1$ and $h_\lambda:  \varphi_\lambda(U) \to \R^{n\times n}$ with $h_{0,(0,0)} = 0$ and $h_\lambda$ symmetric. For $\lambda>0$, the critical points%
\footnote{All vectors in this paper will be column vectors. However, we will often abuse notation to write $(x,z)$ instead of $\vect{x}{z}$ for readability.}
 $(x,z) = (0,\pm\sqrt \lambda)$ are Morse of index difference one and for $\lambda<0$, there are no critical points. 

For $0<\lambda<\epsilon_0^2$, we study the negative gradient trajectories $\Gamma$ with 
$
\dot \Gamma = - \nabla_{G_\lambda} F_\lambda \circ \Gamma
$ 
and  
$
\lim_{t\to \pm \infty} \Gamma(t) =  p_\pm(\lambda) := \varphi_\lambda^{-1} (0, \pm \sqrt \lambda). 
$
Assume $\Gamma(\R) \subset U$, then $\gamma: \R \to \R^n$ given by $\gamma(t) = (\gamma_x(t), \gamma_z(t)) := \varphi_\lambda \circ \Gamma ( a_\lambda t)$, for $t \in \R$ and $\chi_\lambda(w) = a_\lambda w + b_\lambda$, is a solution of 
\begin{equation}\label{eq:trans2}
\begin{cases}
\dot \gamma_x 		&=- A \gamma_x + b(\lambda-\gamma_z^2) + \tilde R^x_\lambda(\gamma),\\
\dot \gamma_z 		&=- \scal{A\gamma_x}{b}+(\lambda-\gamma_z^2) + \tilde R^z_\lambda(\gamma), \\
&\lim_{t\to \pm \infty} \gamma(t) =(0,\pm \sqrt \lambda),
\end{cases}
\,\text{where} \vect{\tilde R_\lambda^x(x,z)}{\tilde R_\lambda^z(x,z)} := - h_{\lambda,(x, z)} \vect{Ax}{\lambda-z^2}.
\end{equation}
Next, we `zoom in', in order to get fast-slow system with fixed limits. Also, the metric will degenerate for $\epsilon = \sqrt \lambda$ goes to zero, which is the standard setting for adiabatic~limits.
%
%
\begin{Lemma}\label{lem:trans1}
 For (A,b,h) as in \eqref{eq:trans1} and $\epsilon_0>0$ as above, we get that $\tilde \gamma:\R \to \R^{n-1} \times \R$  is solution to \eqref{eq:trans2} for $0< \lambda< \epsilon_0^2$ exactly if
$\gamma:\R \to \R^n$ given by 
\begin{equation}\label{eq:trans3}
\gamma_x(t) = \frac{\tilde \gamma_x(t/\epsilon)}{\epsilon^2}, \quad \gamma_z(t) = \frac{\tilde \gamma_z(t/\epsilon)}{\epsilon}, \quad \text{ with } 0< \epsilon:= \sqrt \lambda <\epsilon_0
\end{equation}
is solution to
\begin{equation}\label{eq:eq4}
\begin{cases}
\epsilon \, \dot \gamma_x 		\!\!\!\!\!&=\!- A \gamma_x + b(1-\gamma_z^2) + \epsilon R^x_\epsilon(\gamma),\\
\dot \gamma_z 			\!\!\!\!\!&=\!- \! \scal{A\gamma_x}{b}+(1-\gamma_z^2) + R^z_\epsilon(\gamma), \\
&\lim_{t\to \pm \infty} \gamma(t)=(0,\pm 1),
\end{cases}
\,\text{where} \vect{\epsilon R_\epsilon^x(x,z)}{R_\epsilon^z(x,z)} \! := \! - h_{\epsilon^2,(\epsilon^2 x , \epsilon z)} \! \vect{Ax}{z^2-1}.
\end{equation}
Furthermore, \eqref{eq:eq4} is the negative gradient flow equation of the pair $(f_\epsilon, g_\epsilon)$ given by
\begin{equation}\label{eq:eq5}
\begin{cases}
&f_\epsilon(x,z)= \frac{ \epsilon}{2} x^\top A x + \frac{1}{3}z^3-z,\\[0.5 em]
&(g_{\epsilon,(x,z)})^{-1} = \begin{pmatrix} \mathbbm 1/\epsilon^2 & b/\epsilon \\ b^\top/\epsilon & 1 \end{pmatrix} + \begin{pmatrix} \mathbbm 1/\epsilon & 0 \\ 0 & 1 \end{pmatrix} h_{\epsilon^2, (\epsilon^2 x,\epsilon z)} \begin{pmatrix} \mathbbm 1/\epsilon & 0 \\ 0 & 1 \end{pmatrix},
\end{cases}
\end{equation}
for $(x,z) \in \R^{n-1}\times \R$ and $0< \epsilon < \epsilon_0$.
\end{Lemma}
%
%
\begin{proof}
We plug in the expressions in \eqref{eq:trans3} into \eqref{eq:trans2} and obtain \eqref{eq:eq4}. Alternatively, with charts.
By taking $\varphi_\epsilon^{(3)}(x,z) = ( x/ \epsilon^2, z/ \epsilon)$ and $\chi_\epsilon^{(3)}(w) = z/\epsilon$, we get for $0< \epsilon < \epsilon_0$
\[
(\chi_\epsilon^{(3)} \circ \chi_{\epsilon^2}) \circ F_{\epsilon^2} \circ (\varphi_\epsilon^{(3)} \circ \varphi_{\epsilon^2})^{-1}=\frac{ \epsilon^2}{\mathfrak c} f_\epsilon \text{ and }((\varphi^{(3)}_\epsilon \circ\varphi_{\epsilon^2})_\ast G_{\epsilon^2})^{-1}= \frac {\mathfrak c} {\epsilon^2} g_\epsilon^{-1}.
\]
So $(f_\epsilon, g_\epsilon)$ in \eqref{eq:eq5} is also a gradient pair for \eqref{eq:eq4}.
\end{proof}
We can now give standing assumptions on triples $(A,b,h)$ and formulate the Main Theorem in local coordinates. We will show in Section \ref{sec:globloc} how we go from the global Main Theorem to its local version Theorem \ref{thm:goal2}.
\begin{Standing}\label{sta:trans1}
Here $A \in \R^{(n-1)\times (n-1)}$ is a symmetric invertible matrix, $b \in \R^{n-1}$ is a column vector of norm $\norm b <1$, and $h: \R \times \R^{n-1}\times \R \to \R^{n\times n}:\, (\epsilon, x, z) \mapsto h_{\epsilon^2,(x,z)}$ is a smooth function with values in the space of symmetric matrices satisfying $h_{0,(0,0)}=0$. 
\end{Standing}
%
%
\begin{Theorem}\label{thm:goal2}
For every triple $(A,b,h)$ as in Standing Assumption \ref{sta:trans1}, there exists a constant $\epsilon_0>0$ such that for all $0<\epsilon<\epsilon_0$, there exists a unique (up to time-shift) solution to \eqref{eq:eq4}. This solution is transverse, i.e. the linearised equation over the solution is surjective.
\end{Theorem}
\begin{proof}
The existence part is proven in Theorem \ref{thm:existence} and uniqueness is proven in Theorem \ref{thm:uniq2}. Transversality is proven in Corollary \ref{cor:existence}.
\end{proof}
%
%
%
\section{Local Existence And Uniqueness} \label{sec:three}
This section contains the steps (III) - (IX) from the introduction, each in a separate subsection. Namely, we use local adiabatic analysis to derive existence in Theorem \ref{thm:existence} and local uniqueness in Proposition \ref{prop:loc1}. $(A,b,h)$ will always satisfy Assumptions \ref{sta:trans1}. 
%
%
%
\subsection{Adiabatic Limit}\label{sec:adia}
As illustrated in Figure \ref{fig:intro1}, after `zooming in' Lemma \ref{lem:trans1}, we blew up the point solution at $\epsilon= \sqrt \lambda =0$ to a non-trivial solution $\gamma_0:\R \to \R^n$, the limit solution. This solution will turn out to be the adiabatic limit of solutions $\gamma$ of \eqref{eq:eq4} for $\epsilon>0$ small. We will now derive a formula for $\gamma_0$.

We can look at the equation \eqref{eq:eq4} for $\epsilon=0$ to formally get
\begin{equation}\label{eq:adia1}
\begin{cases}
0= - A(\gamma_0)_x + b \left(1-(\gamma_0)_z^2 \right), \\
(\dot \gamma_0)_z = - \scal b {A(\gamma_0)_x} + 1- (\gamma_0)_z^2,\\
\lim_{t\to \pm \infty}\gamma_0(t)= (0, \pm 1).
\end{cases}
\end{equation}
Here, $R_\epsilon$ disappears since $h_{0,(0,0)}=0$.
The first line $(\gamma_0)_x =A^{-1} b  (1-(\gamma_0)_z^2)$ is a constraint. Plugging this constraint into the second line, we get
\begin{equation}\label{eq:adia2}
\begin{cases}
(\dot \gamma_0)_z= (1-\norm b^2) (1-(\gamma_0)_z^2),\\
\lim_{t\to \pm \infty} (\gamma_0)_z (t)= \pm 1.
\end{cases}
\end{equation}
At this point, it is worth noting that $(1-\norm b^2)>0$ as $g_\epsilon$ in \eqref{eq:eq5} is a metric. Hence, this differential equation has a unique solution up to time-shift. The limit solution equals
\begin{equation}\label{eq:adia3}
\gamma_0:\R \to \R^n \text{ given by }\gamma_0(t) = \left(\frac{A^{-1} b }{\cosh^2((1-\|b\|^2) t)}, \,  \tanh( (1-\|b\|^2) t)\right).
\end{equation}
%
\subsection{$\epsilon$-weighted Norms}\label{sec:met}

We have an $\epsilon$-dependent Riemannian metric $g_\epsilon$ whose inverse 
\[
g_\epsilon^{-1} = \begin{pmatrix} \mathbbm 1/\epsilon^2 & b/\epsilon \\ b^\top/\epsilon & 1 \end{pmatrix} + \text{term in } h
\] 
was defined in \eqref{eq:eq5}. This will let us define the relevant $\epsilon$-weighted norms in this section.

For $h=0$, we can readily invert the above matrix
\begin{equation}\label{eq:norm1}
\begin{aligned}
g_\epsilon^0= \begin{pmatrix} \epsilon^2 \left( \mathbbm 1+ \frac{b  b ^\top}{1-\norm b^2} \right) & -\epsilon \frac{b }{1-\norm b^2} \\ - \epsilon \frac{ b ^\top}{1-\norm b^2} & \frac{1}{1-\norm b^2} \end{pmatrix}.
\end{aligned}
\end{equation}
The corresponding inner product equals for $\eta=(\xi, \zeta),\, \tilde\eta=(\tilde \xi, \tilde \zeta) \in \R^{(n-1)}\times \R$ to
\begin{equation}\label{eq:norm2}
\begin{aligned}
g_\epsilon^0(\eta,\tilde \eta)=\epsilon^2 \scal \xi {\tilde\xi}+ \fracb \left( \epsilon^2 \scal{\xi}{b} \langle \tilde\xi, b \rangle - \epsilon \scal \xi b \tilde\zeta -\epsilon\langle \tilde \xi, b\rangle \zeta+ \zeta \tilde\zeta \right)
\end{aligned}
\end{equation}
and the associated norm is given by
\begin{equation}\label{eq:norm3}
\begin{aligned}
\norm \eta ^2_\epsilon = \epsilon^2 \norm{\xi}^2 + \fracb (\zeta- \epsilon \scal {b } \xi )^2.
\end{aligned}
\end{equation}
%
This norm satisfies the following equivalence of norms for all $\eta=(\xi, \zeta) \in \R^{n-1} \times \R$
\begin{equation}\label{eq:equiv}
\begin{aligned}
\frac{1}{2} (\epsilon^2 \norm \xi^2 + \abs \zeta^2) \leq \norm{\eta}^2_\epsilon \leq \frac{2}{1-\norm b^2}  (\epsilon^2 \norm \xi^2 + \abs \zeta^2).
\end{aligned}
\end{equation}
\begin{proof}[Proof of \eqref{eq:equiv}]
We use $\norm b <1$. We obtain the second inequality by
\[
\norm{\eta}_\epsilon^2 \leq \epsilon^2 \norm \xi^2 + 2\epsilon^2 \frac{\norm{b}^2}{1-\norm b^2} \norm \xi^2 + \frac{2}{1-\norm b^2} \zeta^2 \leq \frac{2}{1-\norm{b}^2} (\epsilon^2 \norm \xi^2+\zeta^2).
\] 
Whereas the first inequality comes from
\begin{align*}
\norm{\eta}_\epsilon^2 	&= \epsilon^2 \norm \xi^2 + \frac{1}{1-\norm b^2} (\zeta^2-2\epsilon \scal \xi b \zeta + \epsilon^2 \scal{\xi}{b})\\
							&\hspace{-2em}\geq \epsilon^2\norm{\xi}^2+\fracb \left(\zeta^2 -\frac{1+\norm b^2}{2} \zeta^2 -\frac{2}{1+\norm b^2}\epsilon^2\scal{\xi}{b}^2 +\epsilon^2 \scal{\xi}{b}^2\right) \\
							&\hspace{-2em}\geq \epsilon^2 \norm{\xi}^2 \frac{1}{1+\norm{b}^2} + \frac{1}{2} \zeta^2
							\geq \frac{1}{2}(\epsilon^2 \norm \xi^2+\zeta^2). \qedhere
\end{align*}
\end{proof}
%
Next we define Hilbert norms for $\eta \in L^2(\R, \R^n)$ respectively $\eta \in \X$ by
\begin{equation}\label{eq:funcnorm1}
\begin{aligned}
\norm \eta_{L^2_\epsilon}^2 := \int_\R \norm{\eta(t)}_\epsilon^2 \, dt \text{ respectively } \norm{\eta}_{W^{1,2}_\epsilon}^2:= \norm{\dot \eta}^2_{L^2_\epsilon} + \norm{\eta}^2_{L^2_\epsilon}.
\end{aligned}
\end{equation}
Furthermore, we define for $\eta \in L^\infty(\R, \R^n)$ the norm 
\begin{equation}\label{eq:funcnorm3}
\begin{aligned}
\norm{\eta}_{L^\infty_\epsilon} = \epsilon \norm \xi_{L^\infty}  + \norm \zeta_{L^\infty}.
\end{aligned}
\end{equation}
These are related. Namely, there is $C>0$ such that for all $0<\epsilon<1$ and $\eta \in \X$
\begin{equation}\label{eq:funcnorm2}
\begin{aligned}
\epsilon^{1/2} \norm{\eta}_{L_\epsilon^\infty} \leq C \norm{\eta}_{W^{1,2}_\epsilon}.
\end{aligned}
\end{equation}

\begin{proof}[Proof of \eqref{eq:funcnorm2}]
Put $\tilde \xi(t):=\epsilon \xi(\epsilon t)$ and $\tilde \zeta(t):=\zeta(\epsilon t)$. Then we have from \eqref{eq:equiv} that
$
\norm {\tilde \eta}^2_{L^2}  \leq 2 \epsilon \inv \norm {\eta}^2_{L^2_\epsilon},
$
and
$
\norm {\dot{\tilde \eta}}^2_{L^2} \leq 2 \epsilon \norm {\dot \eta}^2_{L^2_\epsilon}.
$
Therefore, using the Sobolev embedding Theorem \cite[Theorem VIII.7]{Bre1}, we get
\[
\norm {\eta}_{L^\infty_\epsilon}^2 = \norm{\tilde\eta}_{L^\infty}^2 \leq C_1 \norm {\tilde \eta}_{W^{1,2}}^2 \leq 2 C_1( \epsilon\inv \norm {\eta}_{L^2_\epsilon}^2 + \epsilon \norm {\dot \eta}_{L^2_\epsilon}^2).\qedhere
\]
\end{proof}
For later use, we also compare the Riemannian metric $g_\epsilon$ in \eqref{eq:eq5} to the norm in \eqref{eq:equiv}.
\begin{Lemma}\label{lem:equiv1}
Fix $M>0$. Then there is $\epsilon_0>0$ such that for all $0<\epsilon< \epsilon_0$, for all $(x,z) \in \R^{n-1}\times \R$ with
$
\norm{h_{\epsilon^2, (\epsilon^2 x, \epsilon z)}} \leq M \epsilon,
$
and for all $  \eta = (  \xi,   \zeta) \in \R^{n-1}\times \R$, we have
\[
\frac{1}{4} (\epsilon^2 \|   \xi \|^2 + |  \zeta|^2) \leq g_{\epsilon,(x,z)} ( \eta,   \eta)  \leq \frac{4}{1-\norm b^2}  (\epsilon^2 \|   \xi \|^2 + |  \zeta|^2).
\]

\end{Lemma}

\begin{proof}
Let us introduce the following notations for $(x,z) \in \R^{n-1}\times \R$
\begin{equation}\label{eq:isol11}
\begin{aligned}
(g_\epsilon^0)^{-1} = \begin{pmatrix} \mathbbm 1/ \epsilon^2 & b/\epsilon \\   b ^\top/\epsilon & 1 \end{pmatrix}, %
& \; (g_{\epsilon,(x,z)})^{-1} =(g_\epsilon^0)^{-1} + \begin{pmatrix} \mathbbm 1/\epsilon & 0\\ 0 & 1 \end{pmatrix} h_{\epsilon^2,(\epsilon^2 x, \epsilon z)} \begin{pmatrix}\mathbbm 1/\epsilon& 0\\ 0 & 1 \end{pmatrix},  \\
\big( \hat g^0\big)^{-1} = \begin{pmatrix} \mathbbm 1 & b \\   b ^\top & 1 \end{pmatrix}, %
&\; \big(\hat g_{\epsilon, (x,z)}\big)^{-1} =  \big( \hat g^0 \big)^{-1} + h_{\epsilon^2,(\epsilon^2 x, \epsilon z)}.
\end{aligned}
\end{equation}
Then we have by assumption that
$
\norm {h_{\epsilon^2,(\epsilon^2 x, \epsilon z)} } \norm{\hat g^0} \leq M \epsilon  \norm{\hat g^0} <\frac{1}{2},
$
for all $0<\epsilon<\epsilon_0$ for $\epsilon_0$ small. So we may apply Theorem 1.5.5 from \cite{DietmarFA}, to get
\begin{align*}
\norm{\hat g_{\epsilon, (x,z)}-\hat g^0} \leq \frac{\norm {h_{\epsilon^2,(\epsilon^2 x, \epsilon z)}} \norm{\hat g^0}^2 }{1- \norm {h_{\epsilon^2,(\epsilon^2 x, \epsilon z)}} \norm{\hat g^0}} \leq 2 M \epsilon \norm{\hat g^0}^2.
\end{align*}
By \eqref{eq:equiv}, we have
$
1/2 \abs{  \eta }^2 \leq \hat g^0\left(   \eta ,   \eta  \right) \leq 2/(1-\norm b^2) \abs{  \eta }^2
$
 for all $  \eta \in \R^n$.
So in the end, we get for all $  \eta = (  \xi,   \zeta) \in \R^{n-1}\times \R$,
\begin{equation*}
\begin{aligned}
(g_{\epsilon,(x,z)}) \left(  \eta,   \eta \right) 	&= (\hat g_{\epsilon,(x,z)}) \left(\vect {\epsilon   \xi}{  \zeta}, \vect {\epsilon   \xi}{  \zeta} \right)
															=\left( \hat g^0 + (\hat g_{\epsilon, (x,z)} - \hat g^0) \right) \left(\vect {\epsilon   \xi}{  \zeta}, \vect {\epsilon   \xi}{  \zeta} \right)\\
															&\geq \big(\frac{1}{2}  - 2 M \epsilon \norm{\hat g^0}^2 \big) \abs{\vect {\epsilon   \xi}{  \zeta} }^2 
															\geq \frac{1}{4} \big(\epsilon^2 \|  \xi\|^2+ |  \zeta|^2\big) 
\end{aligned}
\end{equation*}
for $\epsilon_0>0$ maybe smaller. %
Similarly, we get for all $  \eta = (  \xi,   \zeta) \in \R^{n-1}\times \R$ that
\[
g_{\epsilon,(x,z)} \left(  \eta,   \eta \right) \leq \frac{4}{1-\norm b^2} \big(\epsilon^2 \|  \xi\|^2+ |  \zeta|^2\big).  \qedhere
\]
\end{proof}
%
%
%
\subsection{Linearise The Equations}
We linearise our equations \eqref{eq:eq4} at the limit solution $\gamma_0$ defined in \eqref{eq:adia3}.

We set $\gamma = \gamma_0+ \eta$, where $\eta =(\xi, \zeta) \in \X$. We plug this into \eqref{eq:eq4} while using the defining equation for the limit solution in \eqref{eq:adia1}
\begin{equation}\label{eq:exis6}
\begin{cases}
\epsilon \, \dot \xi &= -   {A  \xi} - b  ( 2 (\gamma_0)_z   \zeta +   \zeta^2) + \epsilon R_\epsilon^x(\gamma_0+\eta) - \epsilon (\dot \gamma_0)_x,\\
\dot \zeta&=  \scal b {A \xi} - (2 (\gamma_0)_z   \zeta +   \zeta^2)  + R_\epsilon^z(\gamma_0+\eta).
\end{cases}
\end{equation}
Note that such a $\gamma$ has automatically the right boundary conditions\footnote{Cf. \cite[Corollary VIII.8.]{Bre1}.} in \eqref{eq:eq4}.  
Setting $h=0$ in \eqref{eq:exis6} and linearising the equation at $\eta=0$, we obtain the first part of our linearisation as $D_\epsilon: W^{1,2}(\R, \R^n)  \to L^2(\R, \R^n)$ given by
\begin{equation}\label{eq:exis7}
\begin{aligned}
\hspace{-0.5em} D_\epsilon \eta:= \dot \eta + Q_\epsilon  {\eta}, \text{ where } Q_\epsilon= \begin{pmatrix}   {A} / \epsilon & 2 (\gamma_0)_z \,b  /\epsilon \\  {b }^\top   {A } & 2 (\gamma_0)_z \end{pmatrix}. 
\end{aligned}
\end{equation}
The second part of the linearisation $E_\epsilon$ (defined below) has to do with $h$ not being zero in general. 
We define $\mathcal F_\epsilon: W^{1,2}(\R, \R^n) \to L^2(\R, \R^n)$ by
\begin{align} 
\mathcal F_\epsilon (\eta):= D_\epsilon \eta + \begin{pmatrix} b/\epsilon \\ 1 \end{pmatrix}   \zeta^2 - R_\epsilon (\gamma_0+\eta)-\vect{(\dot \gamma_0)_x}{0}. \label{eq:exis9}
\end{align}
Then we have for $\eta \in \X$
\begin{equation}\label{eq:exis10}
\begin{aligned}
&\gamma_{0}+\eta \text{ is solution of } \eqref{eq:eq4} \iff \mathcal F_\epsilon (\eta)=0.
\end{aligned}
\end{equation}
The linearisation of $\mathcal F_\epsilon$ at $\eta=0$ is given by
\[
d \mathcal F_\epsilon(0) = D_\epsilon + E_\epsilon,
\]
where $E_\epsilon: \X \to \Y$ is the error term given by
\begin{equation}\label{eq:exis11}
E_\epsilon \eta := -d R_\epsilon (\gamma_0) \eta.
\end{equation}
%
%
%
\subsection{Linear Estimates}\label{sec:lin}
In this section, we establish linear estimates for $D_\epsilon$, $E_\epsilon$ and $d \mathcal F_\epsilon(0):=D_\epsilon+E_\epsilon$ which are needed for the existence result and local uniqueness.
\subsubsection{Linear Estimates For $D_\epsilon$ Defined In \eqref{eq:exis7}}
The triple $(A,b,h)$ is as in Assumptions \ref{sta:trans1}. The limit solution $\gamma_0$ was defined in \eqref{eq:adia3}. We will be interested in estimates related to the operator $D_\epsilon$ defined in \eqref{eq:exis7}. 
We start by calculating the adjoint operator $D_\epsilon^\ast$ in \ref{lem:lin1}, followed by a proof of surjectivity of $D_\epsilon$ in Proposition \ref{prop:surj}. With some more work, we end up with the crucial estimate in Proposition \ref{prop:lin2}.  
%
%
\begin{Lemma}\label{lem:adjoint}
The operator $D_\epsilon^\ast: \X \to \Y$ given by
\begin{equation}\label{eq:lin8}
\begin{aligned}
D_\epsilon^\ast \eta:= -\dot \eta + Q_\epsilon \eta,\qquad \text{ where }\qquad Q_\epsilon= \begin{pmatrix}   {A} / \epsilon & 2 (\gamma_0)_z \,b /\epsilon \\  b^\top   A & 2 (\gamma_0)_z \end{pmatrix},
\end{aligned}
\end{equation}
is the adjoint operator of $D_\epsilon$ with respect to the $L^2$-norm \eqref{eq:funcnorm1}, seen as an unbounded operator with dense domain from $L^2 \to L^2$.
\end{Lemma}
\begin{proof}
We can write $Q_\epsilon= \left( g_\epsilon^0 \right)^{-1} B_\epsilon$ where $B_\epsilon=\begin{pmatrix} \epsilon A & 0 \\ 0 & 2 (\gamma_0)_z \end{pmatrix}$. Thus we have
\[
\geps{\eta'}{Q_\epsilon \eta} = \eta'^\top B \zeta  = \eta'^\top  B^\top \eta =  \geps{Q_\epsilon  \eta'}{\eta} 
\]
for all $\eta', \eta \in \R^n$. 
So one verifies by integration by part that 
\begin{align*}
\int_\R   \geps{\eta'}{ D_\epsilon \eta} \, dt &= \int_\R (\geps{\eta'}{ \dot \eta} + \geps{\eta'}{ Q_\epsilon \eta}) \, dt = \int_\R ( \geps{- \dot \eta'}{\eta}+\geps{Q_\epsilon \eta'}{\eta}) \, dt \\
&=  \int_\R   \geps{D^{\ast}_\epsilon \eta'}{\eta} \, dt. \qedhere
\end{align*}
\end{proof}
Here is a Lemma obtained by doubling the operators.
\begin{Lemma}\label{lem:lin1}
For $\eta =(\xi, \zeta)\in \X$, we have the following equalities. 
\begin{align*}
\norm {D_\epsilon \eta}^2_{L^2_\epsilon} = \norm {\dot \eta} ^2_{L^2_\epsilon} - 2\int_\R (\dot \gamma_0)_z\, \zeta^2 + \norm {Q_\epsilon \eta}^2_{L^2_\epsilon},\\
\norm {D_\epsilon^\ast \eta}^2_{L^2_\epsilon} = \norm {\dot \eta} ^2_{L^2_\epsilon} + 2\int_\R (\dot \gamma_0)_z\, \zeta^2 + \norm {Q_\epsilon \eta}^2_{L^2_\epsilon}.
\end{align*}
\end{Lemma}
\begin{proof}
We start by doubling the operator for $\eta \in C^\infty_0(\R, \R^n)$
\begin{equation}\label{eq:double}
D_\epsilon^\ast D_\epsilon \eta= \left(  -\frac{d}{d t} + Q_\epsilon \right) (\dot \eta + Q_\epsilon \eta)= -\ddot \eta - \dot Q_\epsilon \eta + Q_\epsilon^2 \eta
\end{equation}
and pairing this expression \eqref{eq:double} with $\eta$. We obtain the identity
\[
\norm {D_\epsilon \eta}^2_{L^2_\epsilon} = \norm {\dot \eta} ^2_{L^2_\epsilon} - \scal {\dot Q_\epsilon \eta}{\eta}_{L^2_\epsilon} + \norm {Q_\epsilon \eta}^2_{L^2_\epsilon}
\]
which is also true for all $\eta \in \X$ by a density argument.
Furthermore,
\begin{equation}\label{eq:lin2}
g^0_\epsilon ( \dot Q_\epsilon \eta, \eta)= g^0_\epsilon \left( \left[ 
\begin{array}{cc}
0		& 2(\dot \gamma_0)_z b/\epsilon\\
0 		& 2 (\dot \gamma_0)_z 
\end{array} \right] \eta, \eta \right) = 2 (\dot \gamma_0)_z\, \zeta \, \geps {\vect{b/\epsilon}{1}}{ \eta} = 2 (\dot \gamma_0)_z \zeta^2.
\end{equation}
where the last equality follows from $(g_\epsilon^0)^{-1}= \begin{pmatrix}\mathbbm 1/\epsilon^2 & b/\epsilon \\ b/\epsilon & 1 \end{pmatrix}$.
Similar calculations can be done for $D_\epsilon^\ast$.
\end{proof}
The following Proposition proves that $D_\epsilon$ is surjective.
\begin{Proposition}\label{prop:surj}
There exist $C, \epsilon_0>0$ such that for $0<\epsilon\leq \epsilon_0$, we have 
\begin{equation}\label{eq:surj}
\norm \eta_{W^{1,2}_\epsilon} \leq C \norm {D_\epsilon^\ast \eta} _{L^2_\epsilon},
\end{equation}
for all $\eta \in \X$ and where the norms are defined in \eqref{eq:funcnorm1}. %
This means in particular, that $D_\epsilon$ is surjective, and that $D_\epsilon^\ast$ is injective with closed range. %
\end{Proposition}
\begin{proof}
By Lemma \ref{lem:lin1}, we have for $\eta = (\xi, \zeta) \in \X$
\begin{equation}\label{eq:lin9}
\normw{\eta}^2 =  \norm {\eta}^2_{L^2_\epsilon} +  \norm {\dot \eta}^2_{L^2_\epsilon} = \left(\norm {D_\epsilon^\ast \eta}_{L^2_\epsilon}^2-2 \int_\R (\dot \gamma_0)_z \zeta^2 -\norm {Q_\epsilon\eta}_{L^2_\epsilon}^2 \right)+\norm {\eta}_{L^2_\epsilon}^2.
\end{equation}
We use that $Q_\epsilon= \left( g_\epsilon^0 \right)^{-1} B_\epsilon$ where $B_\epsilon=\begin{pmatrix} \epsilon A & 0 \\ 0 & 2 (\gamma_0)_z \end{pmatrix}$, and get
\begin{equation}\label{eq:lin13}
\begin{aligned}
\norm {Q_\epsilon\eta}_{L^2_\epsilon}^2 	&= \int_\R \eta^\top B_\epsilon (g_\epsilon^0)^{-1} B_\epsilon \eta\, dt = \norm{A\xi}^2_{L^2} + 4 \norm{( \gamma_0)_z \zeta}^2_{L^2} + 4 \int_\R ( \gamma_0)_z \zeta \scal{A\xi} b \\
					&\geq (1-\norm b^2) \norm {A\xi}^2_{L^2} \geq (1-\norm b^2) \kappa^2 \norm \xi^2_{L^2},
\end{aligned}
\end{equation}
where we used Cauchy-Schwarz and $\kappa>0$ is such that
$
\norm {A x} \geq \kappa \norm x
$
for $x\in \R^{n-1}$ which exists as $A$ is symmetric and invertible.
Furthermore, we have
\begin{align*}
\norm \eta^2_{L^2_\epsilon} 		&\leq \frac{2}{1-\norm b^2} \left( \epsilon^2 \norm{\xi}^2_{L^2} + \norm{\zeta}^2_{L^2} \right) 
\end{align*}
by \eqref{eq:equiv}. %
Putting these two estimates back into \eqref{eq:lin9}, we get
\begin{align*}
&\norm {\eta}^2_{L^2_\epsilon} +  \norm {\dot \eta}^2_{L^2_\epsilon}\\ 		
&\leq \left(-(1-\norm b^2) \kappa + \epsilon^2 \frac{2}{1-\norm b^2}  \right) \norm \xi^2_{L^2} + \frac{2}{1-\norm b^2} \norm\zeta^2_{L^2} + \norm {D^\ast_\epsilon \eta}^2_{L^2_\epsilon} - 2 \int_\R (\dot \gamma_0)_z \zeta^2 \\
&\leq \frac{2}{1-\norm b^2} \norm\zeta^2_{L^2} + \norm {D^\ast_\epsilon \eta}^2_{L^2_\epsilon}
\end{align*}
where we chose $\epsilon_0 >0$ smaller such that $-(1-\norm b^2) \kappa \norm \xi^2_{L^2}$ dominates the other term in $\normltwo{\xi}^2$ and where we used the fact that $(\dot \gamma_0)_z>0$. 
So it only remains to find an estimate for $\norm\zeta^2_{L^2}$. Towards this goal, we use Lemma \ref{lem:lin1} and the definition of the norm in \eqref{eq:norm3}
\begin{align*}
&\norm {D^\ast_\epsilon \eta}_{L^2_\epsilon}^2 	= \norm {\dot \eta}_{L^2_\epsilon}^2 + 2 \int_\R (\dot \gamma_0)_z \zeta^2 + \norm {Q_\epsilon\eta}_{L^2_\epsilon}^2\\
							&\geq 2 \int (\dot \gamma_0)_z \zeta^2 + \norm {A \xi + 2( \gamma_0)_z b \zeta}^2_{L^2} + \fracb \int_\R (\scal {A\xi} b +\! 2 (1-\norm b^2) ( \gamma_0)_z \zeta -\! \scal{A\xi} b )^2	\\
							&\geq 2(1-\norm b^2) \int_\R \zeta^2 - 2 (1-\norm b^2) \int_\R ( \gamma_0)_z^2 \zeta^2 + 4 (1-\norm b^2) \int_\R ( \gamma_0)_z^2 \zeta^2 					
							\geq 2 (1-\norm b^2) \int_\R \zeta^2
\end{align*}
where we used the defining differential equation 
$
(\dot \gamma_0)_z =(1-\norm b ^2) (1-( \gamma_0)_z^2)
$
in line three.
Hence combining all of the above, we get
\[
\norm {\eta}_{W^{1,2}_\epsilon} \leq \sqrt{1 + \left(\fracb\right)^2} \norm {D_\epsilon^\ast \eta}_{L^2_\epsilon}
\]
for all $0<\epsilon\leq\epsilon_0$ and $\eta \in \X$.

The statements about $D_\epsilon$ being surjective and $D_\epsilon^\ast$ being injective with closed image follow by the closed image Theorem which can be found in Theorem 6.2.3 of \cite{DietmarFA}.
\end{proof}
A weaker estimate holds for $D_\epsilon$. This is expected since invariance under time-shift of \eqref{eq:eq4} means that its linearisation $D_\epsilon$ has a non-trivial kernel.
\begin{Proposition} \label{prop:deps} There is a constant $C, \epsilon_0>0$, such that for all $0<\epsilon<\epsilon_0$, we have
\begin{equation}\label{eq:lin47}
\begin{aligned}
\norm {\eta}_{W^{1,2}_\epsilon} \leq C ( \norm {D_\epsilon \eta}_{L^2_\epsilon} + \norm {\zeta}_{L^2})
\end{aligned}
\end{equation}
 and for all $\eta=(\xi, \zeta) \in \X$ where the norms are defined in \eqref{eq:funcnorm1}.
\end{Proposition}
\begin{proof}
Use Lemma \ref{lem:lin1} and similar estimates as in the proof of Proposition \ref{prop:surj}, to get
\begin{align*}
&\norm {\eta}^2_{L^2_\epsilon} +  \norm {\dot \eta}^2_{L^2_\epsilon} 	= \left(\norm {D_\epsilon \eta}_{L^2_\epsilon}^2 + 2 \int_\R ( \dot \gamma_0)_z \zeta^2 -\norm {Q_\epsilon\eta}_{L^2_\epsilon}^2 \right)+\norm {\eta}_{L^2_\epsilon}^2\\
													&\leq \norm {D_\epsilon \eta}_{L^2_\epsilon}^2+\hspace{-0.2em}\bigg(\frac{2 \epsilon^2}{1-\norm b^2} -(1-\norm b^2) \kappa \bigg)\hspace{-0.2em} \norm \xi ^2_{L^2}
													+\hspace{-0.2em} \bigg(2(1-\norm b^2)+ \frac{2}{1-\norm b^2}\bigg) \hspace{-0.2em}\norm \zeta_{L^2}^2 \\
													&\leq  \norm {D_\epsilon \eta}_{L^2_\epsilon}^2+ C \norm \zeta^2_{L^2}										
\end{align*}
where $\epsilon_0>0$ small and we use $\norm {( \dot \gamma_0)_z}_{L^\infty} \leq (1-\norm b^2)$ in the second inequality.			
\end{proof}
%
The observation that $D_\epsilon$ does not have a kernel on $\image(D_\epsilon^\ast)$ leads to the following result.  
\begin{Proposition}\label{prop:lin2}
There are $C, \epsilon_0>0$, such that for all $0<\epsilon<\epsilon_0$, we have 
\begin{equation}\label{eq:crucial}
\begin{split}
\norm {D_\epsilon^\ast \upsilon}_{W^{1,2}_\epsilon} \leq C  \norm {D_\epsilon D_\epsilon^\ast \upsilon }_{L^2_\epsilon} 
\end{split}
\end{equation}
for all $\upsilon \in \W$ and where the norms are defined in \eqref{eq:funcnorm1}.
\end{Proposition}
\begin{proof}
We start by plugging in $\eta=D_\epsilon^\ast \upsilon$ into equation \eqref{eq:lin47}, to get
\[
\norm{D_\epsilon^\ast \upsilon}_{W^{1,2}_\epsilon}\leq C_1 \left(  \norm {D_\epsilon D_\epsilon^\ast \upsilon} _{L^2_\epsilon} +\norm{\pi_\zeta(D_\epsilon^{\ast} \upsilon)}_{L^2} \right)
\]
where $\pi_\zeta(\eta)=\zeta$ is the projection on the $\zeta$ part.
Next we estimate using \eqref{eq:equiv} and Lemma \ref{lem:adjoint} that for $D_\epsilon^\ast \upsilon \neq 0$
\begin{align*}
\norm{\pi_\zeta(D_\epsilon^{\ast} \upsilon)}_{L^2}			&	\leq C_2 \norm{D_\epsilon^\ast\upsilon}_{L^2_\epsilon} %
														= \frac{C_2}{\norm{D_\epsilon^\ast\upsilon}_{L^2_\epsilon}} \scal{D_\epsilon^\ast \upsilon}{D_\epsilon^\ast \upsilon}_{L^2_\epsilon} 
														= \frac {C_2} {\norm{D_\epsilon^\ast\upsilon}_{L^2_\epsilon}} \scal{\upsilon}{D_\epsilon D_\epsilon^\ast \upsilon}_{L^2_\epsilon}\\											
													&\leq \frac{C_2  \norml{\upsilon}}{\norml{D_\epsilon^\ast\upsilon}} \norm{D_\epsilon D_\epsilon^\ast \upsilon}_{L^2_\epsilon}\leq C_3 \norm{D_\epsilon D_\epsilon^\ast \upsilon}_{L^2_\epsilon}
\end{align*}
where the last line follows from Proposition \ref{prop:surj}. \qedhere

\end{proof}
%
%
\subsubsection{Linear Estimates For $E_\epsilon$ Defined In \eqref{eq:exis11} And $d \mathcal F_\epsilon(0):= D_\epsilon + E_\epsilon$}\label{sec:err}
In this section, we study the part of the linearisation $E_\epsilon$ which stems from the general metric $(h \neq 0)$ and prove that the behaviour of $d \mathcal F_\epsilon(0)=D_\epsilon + E_\epsilon$ and $D_\epsilon$ do not differ much on the image of the adjoint operator $D_\epsilon^\ast$. We start by expressing $E_\epsilon$ in terms of $h$. We use $\mathcal M_\epsilon = \diag(\epsilon, \ldots, \epsilon, 1)$ to express $(\epsilon^2 x, \epsilon z) = \epsilon \mathcal M_\epsilon (x,z)$ for $(x,z) \in \R^{n-1}\times \R$. Also in this notation, we express $R_\epsilon$ in \eqref{eq:eq4} as $\mathcal M_\epsilon R_\epsilon(x,z)  = - h_{\epsilon^2, \epsilon \mathcal M_\epsilon (x,z)} (A x, z^2-1)$. Thus
\begin{equation}\label{eq:lin10}
\begin{aligned}
\mathcal M_\epsilon E_\epsilon \eta 	&= -d(\mathcal M_\epsilon R_\epsilon(\gamma_0))\eta \\
								&=  (\epsilon \mathcal M_\epsilon \eta)^\top (dh_{\epsilon^2})_{\epsilon \mathcal M_\epsilon\gamma_0} \vect{(A \gamma_0)_x}{(\gamma_0)_z^2-1} + h_{\epsilon^2, \epsilon \mathcal M_\epsilon \gamma_0} \vect{A \xi}{ 2 (\gamma_0)_z \zeta}
\end{aligned}
\end{equation}
for all $\eta \in \X$ and where $\gamma_0$ is the limit solution \eqref{eq:adia3}. Also in this notation, we can rewrite the equivalence of norms \eqref{eq:equiv} for $\eta \in \R^n$ as 
\begin{equation}\label{eq:equiv2}
\frac{1}{2} \norm{\mathcal M_\epsilon \eta}^2 \leq  \norm{\eta}_\epsilon^2 \leq \frac 2 {1-\norm b^2} \norm{\mathcal M_\epsilon \eta}^2.
\end{equation}
We will now show that on $\image D_\epsilon^\ast$ the operator $E_\epsilon$ is small compared to $D_\epsilon$. 
%
%
\begin{Proposition}\label{prop:err1}
There is $C, \epsilon_0>0$ such that for all $0<\epsilon<\epsilon_0$, we have
\[
\norml{E_\epsilon D_\epsilon^\ast \upsilon} \leq \epsilon C \norml{D_\epsilon D_\epsilon^\ast \upsilon}
\]
for all $\upsilon \in \W$ where $D_\epsilon$ and $D_\epsilon^\ast$ is defined in \eqref{eq:exis7} resp. \eqref{eq:lin8} and  the norm is defined in \eqref{eq:funcnorm1}. 
\end{Proposition}
\begin{proof}
We use $h_{0,(0,0)}=0$ from the Assumptions \ref{sta:trans1} to get by Taylor expansion the estimate
\begin{equation}\label{eq:lin11}
\norm{h_{\epsilon^2, \epsilon \mathcal M_\epsilon \gamma_0}} \leq 2 \epsilon \norm{h}_{C^1(B_2(0))}(1+\normlinfty{\gamma_0}) \leq C_1 \epsilon
\end{equation}
for all $0<\epsilon<\epsilon_0\leq 1$ and $\epsilon_0>0$ small such that $\norm{\epsilon_0 \mathcal M_{\epsilon_0} \gamma_0}_{L^\infty}<1$. So all together \eqref{eq:lin10}, \eqref{eq:equiv2} and \eqref{eq:lin11} give for all $\eta=(\xi, \zeta) \in \X$ 
\begin{equation}\label{eq:lin12}
\begin{aligned}
\norml{E_\epsilon \eta} \leq C_2 \epsilon \norml \eta + C_2 \epsilon \normltwo \xi.
\end{aligned}
\end{equation}
This estimate is not yet sufficient. We thus go on for $\eta = (\xi, \zeta) \in \X$
\begin{align*}
\epsilon^2 \norml \xi^2 	&\leq \epsilon^2/\kappa^2 \norml{A \xi}^2 = \epsilon^2/\kappa^2 \norml{Q_\epsilon \vect{\xi}{0}}^2 = \epsilon^2/\kappa^2 \norml{Q_\epsilon \eta - Q_\epsilon \vect{0}{\zeta}}^2\\
					&\leq C_3 \epsilon^2 (\norml{D_\epsilon \eta}^2 + \norm{\zeta}_{L^2}^2).
\end{align*}
Here we used that there is $\kappa>0$ such that $\kappa \norm{\xi} \leq \norm{A \xi}$, Lemma \ref{lem:lin1}, and the identities $\norm{A x} = \norm{Q_\epsilon (\xi,0)}_\epsilon$ for $\xi\in \R^{n-1}$ and $\norml{Q_\epsilon (0,\zeta)}^2 = 4 \norm{(\gamma_0)_z \zeta}_{L^2}^2 \leq 4 \normltwo{\zeta}$ for $\zeta \in \R$ which are special cases of \eqref{eq:lin13}. Plugging this back into \eqref{eq:lin12} with $\eta = D_\epsilon^\ast \upsilon$ for $\upsilon \in \W$ and using Proposition \ref{prop:lin2}, we get by equivalence of norms \eqref{eq:equiv}
\[
\norml{E_\epsilon \upsilon} \leq C_4 \epsilon ( \norml{D_\epsilon D_\epsilon^\ast \upsilon}+\norml{D_\epsilon^\ast \upsilon}) \leq C_5 \epsilon \norml{D_\epsilon D_\epsilon^\ast \upsilon}. \qedhere
\]
\end{proof}
Next, we prove a lemma that gives conditions under which a bounded linear operator $A: \X \to \Y$ restricted to $\image D_\epsilon^\ast$ is surjective. These conditions are met by the full linearisation $d \mathcal F_\epsilon(0):= D_\epsilon+E_\epsilon$ as shown in Corollary \ref{cor:lin1} below.
\begin{Lemma}\label{lem:lin2}
Given $A:\X \to \Y$ a bounded operator and assume there is $\beta, C, \epsilon_1>0$ such that for all $0<\epsilon<\epsilon_1$ and all $\upsilon \in \W$
\begin{equation}\label{eq:lin20}
\begin{aligned}
\norm{A D_\epsilon^\ast \upsilon  - D_\epsilon D^\ast_\epsilon \upsilon}_{L^2_\epsilon} \leq C\epsilon^\beta \norm{D_\epsilon D_\epsilon^\ast \upsilon}_{L^2_\epsilon}.
\end{aligned}
\end{equation}
Then there is $0<\epsilon_0<\epsilon_1$ such that $A$ restricted to $\image D_\epsilon^\ast$ is surjective. 
\end{Lemma}
%
%
\begin{proof}
We know by Proposition \ref{prop:lin2} that there is $C>0$ such that,
\begin{equation}\label{eq:lin21}
\begin{aligned}
\norm{D_\epsilon^\ast \upsilon}_{W^{1,2}_\epsilon} \leq C \norm{D_\epsilon D_\epsilon^\ast \upsilon}_{L^2_\epsilon}.
\end{aligned}
\end{equation}
This proposition requires $\epsilon_0>0$ sufficiently small.
Thus for $\mu \in \Y$, define inductively a sequence $(\upsilon_n,\mu_n)_n$ for $n\in \N$
\begin{align*}
\mu=: D_\epsilon D_\epsilon^\ast \upsilon_0, \quad \mu_{n}:= A D_\epsilon^\ast \upsilon_{n}, \quad \mu-\sum_{k=0}^n \mu_k=: D_\epsilon D_\epsilon^\ast \upsilon_{n+1},
\end{align*}
which is possible, as $D_\epsilon D_\epsilon^\ast$ is surjective by Proposition \ref{prop:surj}.
One can prove inductively that for $w_n:= \sum_{k=0}^n \mu_k \in \image{A}$ and $x_n:=\sum_{k=0}^n D^\ast_\epsilon \upsilon_k$, we have by \eqref{eq:lin20} and \eqref{eq:lin21}
\begin{align*}
&\norm{\mu-w_n}_{L^2_\epsilon} \leq (C\epsilon^\beta)^{n+1}\norm \mu_{L^2_\epsilon},  \\
&\norm{x_n-x_{n-1}}_{W^{1,2}_\epsilon} \leq C \norm{D_\epsilon D_\epsilon^\ast \upsilon_n}_{L^2_\epsilon}= C \norm{\mu-w_{n-1}}_{L^2_\epsilon} \leq C (C\epsilon^\beta)^{n}\norm \mu_{L^2_\epsilon}
\end{align*}
for $n\in \N$.
Now we choose $\epsilon_0>0$ maybe even smaller, such that $C \epsilon_0^\beta<1$, and we get $(x_n)_n \subset \image D_\epsilon^\ast \cap \X$ is Cauchy and so there is $x \in \image D_\epsilon^\ast \cap \X$ with $x= \lim_{n\to \infty} x_n$ and $A x = \lim_{n\to \infty} w_n =\mu$, where we used that $A$ is bounded and that $\image D_\epsilon^\ast$ is closed in $\Y$ by Proposition \ref{prop:surj}. As $\mu$ was arbitrary, $A$ restricted to $\image D_\epsilon^\ast$ is surjective.
\end{proof}
%
%
\begin{Corollary}\label{cor:lin1}
There is $\epsilon_0>0$ such that for all $0<\epsilon<\epsilon_0$, the linearised operator at zero $d \mathcal F_\epsilon(0)=D_\epsilon+E_\epsilon$ restricted to $\image D_\epsilon^\ast$ is surjective. 
\end{Corollary}
%
%
\begin{proof} 
We know by Proposition \ref{prop:err1} that
\begin{align*}
\norm{d \mathcal F_\epsilon(0) D_\epsilon^\ast \upsilon  - D_\epsilon D^\ast_\epsilon \upsilon}_{L^2_\epsilon} =  \norm{E_\epsilon D_\epsilon^\ast \upsilon}_{L^2_\epsilon} \leq \epsilon C \norm{D_\epsilon D^\ast_\epsilon \upsilon}_{L^2_\epsilon},
\end{align*}
for all $\upsilon \in \W$ and $0<\epsilon<\epsilon_0$. So we can apply Lemma \ref{lem:lin2} with $A=d \mathcal F_\epsilon(0)$ and $\beta=1$, to finish the proof.
\end{proof}
%
Finally, we extend the crucial estimate in Proposition \ref{prop:lin2} to the full linearisation $d \mathcal F_\epsilon(0)$.
\begin{Corollary}\label{cor:lin2}
There is $\epsilon_0>0$ and a constant $\hat C>0$ such that for all $0<\epsilon<\epsilon_0$,
\begin{align}\label{eq:crux}
\norm{D_\epsilon^\ast \upsilon}_{W^{1,2}_\epsilon} + \epsilon^{\frac 1 2} \norm{D_\epsilon^\ast \upsilon}_{L^\infty_\epsilon} \leq \hat C \norm{d \mathcal F_\epsilon(0) D_\epsilon^\ast \upsilon}_{L^2_\epsilon} 
\end{align}
for all $\upsilon \in \W$, where $d \mathcal F_\epsilon(0):=D_\epsilon+E_\epsilon$ and where the norms are defined in \eqref{eq:funcnorm1} and \eqref{eq:funcnorm3}.
\end{Corollary}
\begin{proof}
By Sobolev embedding \eqref{eq:funcnorm2}, we reduce to only bounding the $W^{1,2}_\epsilon$-norm. Furthermore, by Proposition \ref{prop:lin2}, we only need for $\epsilon_0>0$ small enough to bound the expression
$
\norm{D_\epsilon D_\epsilon^\ast \upsilon}_{L^2_\epsilon}.
$
To this end, we use Proposition \ref{prop:err1} and estimate
\begin{align*}
\norm{D_\epsilon D_\epsilon^\ast \upsilon}_{L^2_\epsilon} \leq  \norm{d \mathcal F_\epsilon(0) D_\epsilon^\ast \upsilon}_{L^2_\epsilon} + \norm{E_\epsilon D_\epsilon^\ast \upsilon}_{L^2_\epsilon} 
\leq  \norm{d \mathcal F_\epsilon(0) D_\epsilon^\ast \upsilon}_{L^2_\epsilon} + C_1 \epsilon \norm{D_\epsilon D_\epsilon^\ast \upsilon}_{L^2_\epsilon}.
\end{align*} 
Thus as soon as $C_1 \epsilon_0< \frac{1}{2}$, we get
$
\norm{D_\epsilon D_\epsilon^\ast \upsilon}_{L^2_\epsilon}\leq 2 \norm{d \mathcal F_\epsilon(0) D_\epsilon^\ast \upsilon}_{L^2_\epsilon}. 
$
\end{proof}
%
\subsection{Quadratic Estimates}\label{sec:quad}
This section is dedicated to quadratic estimates, first for $R_\epsilon$ defined in \eqref{eq:eq4} and then for $\mathcal F_\epsilon$ defined in \eqref{eq:exis9}. These stem from the Hessian of $\frac{1}{3} z^3 - z$ being constant equal to~$2$. These estimates will be the last step before we can establish existence in the next section.
\begin{Lemma}[\textbf{Quadratic Estimates for $R_\epsilon$}]\label{lem:quad1}
There are $M, \epsilon_0>0$ such that for all $0<\epsilon< \epsilon_0$ and for all $\hat \eta, \, \Delta := \gamma - \gamma_0 \in \X$ with $\normlinfty{\hat \eta }, \normlinfty{\Delta} \leq 1$, we have
\begin{align*}
&\norml{(dR_\epsilon(\gamma) - dR_\epsilon(\gamma_0))\hat\eta} 
\leq M\big( \normlinfty{\hat \eta} \! (\epsilon \normlinfty{ \Delta} \! +\! \norml{\Delta} \big) \!+ \! \norml{\hat \eta} ( \normlinfty{\Delta}\!+\!  \epsilon \norml{\Delta} )  \big), \\
&\norml{R_\epsilon(\gamma+\hat \eta) - R_\epsilon(\gamma) - d R_\epsilon(\gamma) \hat \eta}
\leq M \bigg( \normlinfty{\hat \eta} \big (\epsilon \normlinfty{ \hat \eta} + \norml{\hat \eta} \big)+ \epsilon \norml{\hat \eta}^2 \bigg),
\end{align*}
where $R_\epsilon$ is defined in \eqref{eq:eq4} and the norms are defined in  \eqref{eq:funcnorm1} and  \eqref{eq:funcnorm3}.
\end{Lemma}
\begin{proof}
Denote by $L(\eta):=R_\epsilon(\gamma_0+\eta)$ and note that 
\[
dR_\epsilon(\gamma)\hat \eta - dR_\epsilon(\gamma_0)\hat\eta= (dL(\Delta)-dL(0))\hat \eta.
\]
We recall the notation $\mathcal M_\epsilon=\diag(\epsilon, \ldots, \epsilon, 1)$ from the beginning of Subsection \ref{sec:err}. Then \eqref{eq:eq4} reads $\mathcal M_\epsilon R_\epsilon(x,z)  = - h_{\epsilon^2, \epsilon \mathcal M_\epsilon (x,z)} (A x, z^2-1)$ for all $(x,z) \in \R^n$. By Taylor expansion, we have with $\Delta=(\Delta^\xi, \Delta^\zeta)$ that
\begin{align*}
&(dL(\Delta)-dL(0))\hat \eta = \int_0^1 d^2 L(\sigma\Delta) \big(\Delta,\hat \eta \big) \, d\sigma = \mathcal M_\epsilon^{-1} \sum_{i=1}^4 S_i,\\
&\text{where }
S_1:= \int_0^1   \epsilon^2\, (\mathcal M_\epsilon \Delta)^\top (d^2h_{\epsilon^2})_{(\epsilon \sigma \mathcal M_\epsilon (\gamma))} \mathcal M_\epsilon \hat \eta   \vect{ {A((\gamma_0)_x + \sigma \Delta^\xi)}}{( {(\gamma_0)_z+\sigma \Delta^\zeta})^2-1} \, d\sigma ,\\
&\phantom{\text{where }}
S_2:= \int_0^1   \epsilon\, (dh_{\epsilon^2})_{(\epsilon \sigma \mathcal M_\epsilon (\gamma))} \mathcal M_\epsilon \hat \eta \vect{ {A \Delta^\xi}}{2  {((\gamma_0)_z+\sigma \Delta^\zeta) \Delta^\zeta}} \, d\sigma ,\\
&\phantom{\text{where }}
S_3:= \int_0^1  \epsilon\, (dh_{\epsilon^2})_{(\epsilon \sigma \mathcal M_\epsilon (\gamma))}  \mathcal M_\epsilon \Delta  \vect{ {A \hat \xi}}{2  {((\gamma_0)_z+\sigma \Delta^\zeta) \hat \zeta}} \, d\sigma , \\
&\phantom{\text{where }}
S_4:=\int_0^1   h_{\epsilon^2,(\epsilon \sigma \mathcal M_\epsilon (\gamma))} \vect{0}{2  { \Delta^\zeta\, \hat \zeta}} \, d\sigma.
\end{align*}
We recall that, by \eqref{eq:equiv2} for some $c>0$,\,  $c^{-1} \norm{\mathcal M_\epsilon \eta} \leq \norm{\eta}_\epsilon\leq c \norm{\mathcal M_\epsilon \eta }$ for all $\eta \in \R^n$. Thus
$
\norml{(dL(\Delta)-dL(0))\hat \eta}\leq c \sum_{i=1}^4 \normltwo{S_i }
$
and as in \eqref{eq:lin11} we have \[
\epsilon^{-1} \normlinftyy{h_{\epsilon^2, (\epsilon \sigma \mathcal M_\epsilon (\gamma))}}+ \normlinftyy{(dh_{\epsilon^2})_{(\epsilon \sigma \mathcal M_\epsilon (\gamma))}} + \normlinftyy{(d^2h_{\epsilon^2})_{(\epsilon \sigma \mathcal M_\epsilon (\gamma))}}\leq C_1
\]  
where $C_1>0$ only depends on $\norm{h}_{C^2(B_3(0))}$, $\norm b$ and $\norm A$  , and where we used the assumptions $\normlinfty{\Delta}\leq 1$ and $\epsilon_0$ small. Furthermore for $C_3>0$, we have for $0<\sigma<1$
\begin{align*}
& \normltwo{\left(A\left((\gamma_0)_x+\sigma\Delta^\xi\right), \left((\gamma_0)_z+\sigma \Delta^\zeta\right)^2-1 \right) }   \\
&\leq C_2\left( \normltwo{(\gamma_0)_x}+ \normltwo{\Delta^\xi} + \normltwo{(\gamma_0)_z^2-1} + 2 \normlinftyy{(\gamma_0)_z}\normltwo{\Delta^\zeta}+ \normltwo{(\Delta^\zeta)^2} \right) \\
&\leq C_3 (1+ \epsilon^{-1} \norml{\Delta}), \quad \normltwo{(A \hat \xi, 2 ((\gamma_0)_z+\sigma \Delta^\zeta) \hat \eta)} \leq \epsilon^{-1} C_3 \norml{\hat \eta},\\
& \normltwo{(A \Delta^\xi, 2 ((\gamma_0)_z+\sigma \Delta^\zeta) \Delta^\zeta)}\leq \epsilon^{-1} C_3 \norml{\Delta}.
\end{align*}
where we used \eqref{eq:equiv}.
Therefore, we see that
\begin{align*}
\hspace{-1.5em}\sum_{i=1}^4 \normltwo{S_i } \leq C_4\bigg( \normlinfty{\hat \eta} \big (\epsilon \normlinfty{ \Delta} + \norml{\Delta} \big) +  \norml{\hat \eta} \big( \normlinfty{\Delta}+  \epsilon \norml{\Delta} \big)  \bigg). 
\end{align*}
For the other expression, we have again by Taylor 
\begin{align*}
&R_\epsilon(\gamma+\hat \eta) - R_\epsilon(\gamma) - d R_\epsilon(\gamma) \hat \eta = L(\Delta+\hat \eta) - L(\Delta)-d L( \Delta) \hat \eta \\
&=\int_0^1 \int_0^1 d^2 L (\Delta + \sigma \tau \hat \eta) ( \hat \eta, \hat \eta) \, d \tau \, d \sigma =: S.
\end{align*}
By a similar argument as above, we can conclude that
\[
\norml{S} \leq C_5 \bigg( \normlinfty{\hat \eta} \big (\epsilon \normlinfty{ \hat \eta} + \norml{\hat \eta} \big)+ \epsilon \norml{\hat \eta}^2 \bigg).\qedhere
\]
\end{proof}
\begin{Corollary}[\textbf{Quadratic Estimates for $\mathcal F_\epsilon$}]\label{cor:quad1}
There are $M, \epsilon_0>0$ such that for all $0<\epsilon< \epsilon_0$ and for all $\hat \eta, \, \Delta \in \X$ with $\normlinfty{\hat \eta }, \normlinfty{\Delta} \leq 1$, we have
\begin{align*}
&\norml{(d\mathcal F_\epsilon(\Delta) - d\mathcal F_\epsilon(0))\hat\eta} 
\leq M\big( \normlinfty{\hat \eta} \! (\epsilon \normlinfty{ \Delta} \! +\! \norml{\Delta} \big) \!+ \! \norml{\hat \eta} ( \normlinfty{\Delta}\!+\!  \epsilon \norml{\Delta} )  \big), \\
&\norml{\mathcal F_\epsilon(\Delta+\hat \eta) - \mathcal F_\epsilon(\Delta) - d \mathcal F_\epsilon(\Delta) \hat \eta}
\leq M \bigg( \normlinfty{\hat \eta} \big (\epsilon \normlinfty{ \hat \eta} + \norml{\hat \eta} \big)+ \epsilon \norml{\hat \eta}^2 \bigg),
\end{align*}
where $\mathcal F_\epsilon$ is defined in \eqref{eq:exis9} and the norms are defined in  \eqref{eq:funcnorm1} and  \eqref{eq:funcnorm3}.
\end{Corollary}
\begin{proof}
We recall the formula for $\mathcal F_\epsilon$ with $\eta \in \X$
\[
\mathcal F_\epsilon (\eta):= D_\epsilon \eta + \begin{pmatrix} b/\epsilon \\ 1 \end{pmatrix}   \zeta^2 - R_\epsilon (\gamma_0+\eta)-\vect{(\dot \gamma_0)_x}{0},
\]
and see that the first and last term will not appear in the quadratic expressions above. Therefore, we get for $\hat \eta = (\hat \xi, \hat \zeta),\,  \Delta=(\Delta^\xi, \Delta^\zeta)$ and $\gamma:= \gamma_0+\Delta$ that 
\begin{align*}
&\left( d \mathcal F_\epsilon (\Delta) - d \mathcal F_\epsilon (0) \right) \hat \eta= 2 \Delta^\zeta \, \vect{b/\epsilon}{1} \hat \zeta + \left( d R_\epsilon(\gamma)-d R_\epsilon(\gamma_0) \right) \hat \eta,\\
&\mathcal F_\epsilon(\Delta + \hat \eta) - \mathcal F_\epsilon(\Delta) - d \mathcal F_\epsilon(\Delta)\hat \eta= \vect{b/\epsilon}{1} \hat \zeta^2+ (R_\epsilon(\gamma+ \hat \eta) - R_\epsilon(\gamma) - d R_\epsilon(\gamma) \hat \eta ).
\end{align*}
Now we can use \eqref{eq:equiv} and $\norm{(b/\epsilon,1)}_\epsilon=1$ from \eqref{eq:lin2}, to conclude that
\[
\norml{2 \Delta^\zeta \, \vect{b/\epsilon}{1} \hat \zeta} = 2 \normltwo{\Delta^\zeta \, \hat \zeta} \leq 4 \normlinfty{\Delta} \norml{\hat \eta}, \quad \norml{\vect{b/\epsilon}{1} \hat \zeta^2} \leq 2 \normlinfty{\hat \eta} \norml{\hat \eta},
\]
and so the result follows from the quadratic estimates on $R_\epsilon$ in Lemma \ref{lem:quad1}.
\end{proof}
%
\subsection{Existence Via Newton}

After establishing linear estimates in section \ref{sec:lin} and quadratic estimates in section \ref{sec:quad}, we can now prove the existence of a solution to \eqref{eq:eq4} near the limit solution $\gamma_0$, which was defined in \eqref{eq:adia3}. We recall from \eqref{eq:exis10}, that $\gamma=\gamma_0+\eta$ is a solution to \eqref{eq:eq4} if and only if $\eta$ is a zero of the functional $\mathcal F_\epsilon$ defined in \eqref{eq:exis9}. Therefore we are interested in finding a zero of a functional, which is exactly what a Newton iteration method does. This method requires a bijective linearisation $d \mathcal F_\epsilon(0)$. Our problem does not meet this condition. So to circumvent this obstacle, we iterate on the slice $\image D_\epsilon^\ast$ which is orthogonal to the kernel. Such a Newton iteration method is typical for adiabatic limit analysis, see e.g. \cite{Sal1}. 
As a corollary, we can prove that the found solution is transverse as shown in Corollary~\ref{cor:existence}.
\begin{Theorem}[\textbf{Existence}] \label{thm:existence}
There are $C,\epsilon_0>0$ such that for all $0<\epsilon<\epsilon_0$, there exists $\eta_\epsilon \in \X$ such that $\gamma_\epsilon:= \gamma_0+ \eta_\epsilon$ is a smooth solution of  \eqref{eq:eq4}. Furthermore, we have
\[
\eta_\epsilon \in \image D_\epsilon^\ast, \qquad \normw {\eta_\epsilon} + \epsilon^{\frac 1 2} \normlinfty{\eta_\epsilon} \leq C \epsilon
\]
where the norms are defined in \eqref{eq:funcnorm1}, \eqref{eq:funcnorm3}, the operator $D_\epsilon^\ast$ is defined in \eqref{eq:lin8} and $\gamma_0$ is the limit solution in \eqref{eq:adia3}.
\end{Theorem}
\begin{proof}
From \eqref{eq:exis10}, we know that $\gamma=\gamma_0+\eta$ for $\eta \in \X$ is a solution of \eqref{eq:eq4} exactly if $\mathcal F_\epsilon (\eta) = 0$. Our starting point is the limit solution $\gamma_0=\gamma_0+0$. We recall that $d \mathcal F_\epsilon(0):=D_\epsilon + E_\epsilon$, where $D_\epsilon, E_\epsilon$ were defined in \eqref{eq:exis7} resp. \eqref{eq:exis11}.

So let $\upsilon_0 \in \W$ be a solution to $d \mathcal F_\epsilon(0) D_\epsilon^\ast \upsilon_0=-\mathcal{F}_\epsilon (0)=((\dot \gamma_0)_x,0) - R_\epsilon(\gamma_0)$. Such a $\upsilon_0$ exists by Corollary \ref{cor:lin1} for all $0<\epsilon<\epsilon_0$. Set
$
\eta_0:= D_\epsilon^\ast \upsilon_0, \, \eta_0=:(\xi_0,\zeta_0).
$
We have
\begin{align*}
\norml {\mathcal F_\epsilon(0)}^2&= \int_\R \norm{\vect{(\dot \gamma_0)_x}{0}+R_\epsilon(\gamma_0) }_\epsilon^2 \ d t \leq \frac{4 \epsilon^2 \norm{ (\dot \gamma_0)_x}_{L^2}^2}{1-\norm b^2}   + 2  \norm{R _\epsilon(\gamma_0)}_{L^2_\epsilon}^2  \leq C_1 \epsilon^2
\end{align*}
where we used \eqref{eq:equiv}, the definition of $R_\epsilon$ in \eqref{eq:eq4} and $\norm{h_{\epsilon^2,(\epsilon^2 ( \gamma_0)_x,\epsilon ( \gamma_0)_z)}}_{L^\infty}\leq C_2 \epsilon$ as in~\eqref{eq:lin11}. 
We get from \eqref{eq:crux} that
\begin{equation}\label{eq:crux2}
\begin{split}
\norm {\eta_0}_{W^{1,2}_\epsilon} + \epsilon^{\frac{1}{2}}\norm {\eta_0}_{L^\infty_\epsilon} \leq \hat C  \norm {-\mathcal F_\epsilon(0)}_{L^2_\epsilon}  \leq C_3 \epsilon. 
\end{split}
\end{equation}
Now we need to estimate the value of $\mathcal F_\epsilon$ on our new solution $\eta_0$. As $\mathcal F_\epsilon(0) = - d \mathcal F_\epsilon(0) \eta_0$, we get
$
\mathcal F_\epsilon (\eta_0) 	= \mathcal F_\epsilon (\eta_0)-\mathcal F_\epsilon(0) - d \mathcal F_\epsilon(0) \eta_0.				
$
Therefore, we get from the quadratic estimate of Corollary \ref{cor:quad1} and the Sobolev embedding \eqref{eq:funcnorm2} for $0<\epsilon<\epsilon_0$
\begin{align}\label{eq:crux3}
&\norm {\mathcal F_\epsilon (\eta_0)}_{L^2_\epsilon}	\leq C_4 \epsilon^{\frac 1 2} \norm {\eta_0}_{W^{1,2}_\epsilon}
\end{align}
where $\epsilon_0>0$ small, such that $\norm {\eta_0}_{L^\infty_\epsilon} \leq C_3 \epsilon_0^{1/2} \leq 1$. 
%
Now continue by defining inductively 
\[
\begin{array}{lll}
\gamma_{k+1}:= \gamma_0 + \sum_{l=0}^{k} \eta_l, 	&d \mathcal F_\epsilon(0) D_\epsilon^\ast \upsilon_k:=-\mathcal F_\epsilon(\gamma_k-\gamma_0), & \eta_k:= D_\epsilon^\ast \upsilon_k, \\
\eta_k=:\vect{\xi_k}{\zeta_k},  			&\Delta_k:= \gamma_k -\gamma_0,			&\text{for } k \in \N.
\end{array}
\]
We will prove by induction that we have the following inequalities for some $C_5\geq C_4$
 \begin{align*}
&\norm {\eta_k}_{W^{1,2}_\epsilon} + \epsilon^\frac{1}{2} \norm{\eta_k}_{L^\infty_\epsilon} \leq 2^{-k} C_3 \epsilon, \tag{$\square_k$} \\
&\norm {\mathcal F_\epsilon(\Delta_{k+1})}_{L^2_\epsilon} \leq C_5 \epsilon^{\frac{1}{2}}\norm {\eta_{k}}_{W^{1,2}_\epsilon}, \tag{$\lozenge_k$}
\end{align*}
for all $k\in \N$, all $0<\epsilon<\epsilon_0$ with $\epsilon_0$ small.
We start our induction, by noting that ($\square_0$) and ($\lozenge_0$) have already been established in \eqref{eq:crux2} and \eqref{eq:crux3}. Assume now for $k\geq 1$ that we proved ($\square_{l}$) and ($\lozenge_l$) for $l=0, \ldots, k-1$. Let us start by proving $(\square_k)$. 
\begin{align*}
\norm {\eta_k}_{W^{1,2}_\epsilon} + \epsilon^\frac{1}{2} \norm{\eta_k}_{L^\infty_\epsilon} 		&\leq \hat C \norm {\mathcal F_\epsilon(\Delta_k)}_{L^2_\epsilon}
																				\leq \hat C C_5 \epsilon^\frac{1}{2} \norm {\eta_{k-1}}_{W^{1,2}_\epsilon}\\
																				&\leq 2^{-1}  \norm {\eta_{k-1}}_{W^{1,2}_\epsilon} \leq \ldots 
																				\leq 2^{-k} \norm {\eta_0}_{W^{1,2}_\epsilon}
																				\leq 2^{-k} C_3 \epsilon
\end{align*}
for all $0<\epsilon<\epsilon_0$ where $\epsilon_0>0$ such that $\hat C  C_5 \epsilon_0^{1/2} \leq 2^{-1}$. Here, we used \eqref{eq:crux} in the first inequality and we used $(\lozenge_{k-1})$ in the second inequality. Using the same argument repeatedly with $(\lozenge_l)$ for $l\leq k-2$, we get to the penultimate inequality. The last one follows from \eqref{eq:crux2}.
Hence, we get from $(\square_l)$ for $l=0,\ldots, k$ that
\begin{equation}\label{eq:difference}
\begin{aligned}
&\norm {\Delta_{k+1}}_{W^{1,2}_\epsilon} \leq \sum_{l=0}^k \norm {\eta_l}_{W^{1,2}_\epsilon} \leq 2 C_3 \epsilon, \; %
\norm {\Delta_{k+1}}_{L^\infty_\epsilon} \leq \epsilon^{-\frac{1}{2}} \sum_{l=0}^k \epsilon^{\frac{1}{2}}\norm {\eta_l}_{L^\infty_\epsilon} \leq 2 C_3 \epsilon^{\frac{1}{2}}.
\end{aligned}
\end{equation}
For proving $(\square_k)$, we observe that due to $d \mathcal F_\epsilon(0) \eta_k=-\mathcal F_\epsilon(\Delta_k)$, 
\begin{align*}
\mathcal F_\epsilon(\Delta_{k+1}) = \left(  \mathcal F_\epsilon(\Delta_{k+1}) - \mathcal F_\epsilon(\Delta_k)- d \mathcal F_\epsilon(\Delta_k) \eta_k \right) + \left( d \mathcal F_\epsilon(\Delta_{k}) - d \mathcal F_\epsilon (0) \right) \eta_k.
\end{align*}
Taking $0<\epsilon_0<1$ small such that $\normlinfty{\Delta_k} + \normlinfty{\eta_k} \leq 3 C_3 \epsilon_0^{1/2} \leq 1$, we can apply our quadratic estimates from Corollary \ref{cor:quad1} with constant $M>0$. Thus
\begin{align*}
\norm {\mathcal F_\epsilon(\Delta_{k+1})}_{L^2_\epsilon} 	
&\leq \norm {\mathcal F_\epsilon(\Delta_{k+1}) - \mathcal F_\epsilon(\Delta_k)- d \mathcal F_\epsilon(\Delta_k) \eta_k}_{L^2_\epsilon} 
												+\norm{(d \mathcal F_\epsilon(\Delta_k) - d \mathcal F_\epsilon(0)) \eta_k}_{L^2_\epsilon}\\
												&\leq  6 C_3 (C_6+1) M \epsilon^{\frac 12} \normw{\eta_k}=: C_5 \epsilon^{\frac 12} \normw{\eta_k}
\end{align*}
where $C_6>0$ stems from Sobolev embedding \eqref{eq:funcnorm2}. We also used  \eqref{eq:difference} and $(\square_k)$. 

Thus, we get by $(\lozenge_k)$ that for fixed $0<\epsilon<\epsilon_0$, $\gamma_k$ is converging in $\X$ to some $\gamma_\epsilon=\gamma_0+\eta_\epsilon$, such that by \eqref{eq:difference}, $(\lozenge_k)$ and $(\square_k)$, 
\[
\norm {\gamma_\epsilon-\gamma_0}_{W^{1,2}_{\epsilon}} \leq 2 C_3 \epsilon \qquad \text{ and }  \mathcal F_\epsilon(\eta_\epsilon)=0.
\] 
In addition, given that $\image D_\epsilon^\ast$ is closed by Proposition \ref{prop:surj}, $\gamma_\epsilon-\gamma_0 \in \image D_\epsilon^\ast$. By bootstrapping, solutions $\eta_\epsilon$  in $\X$ of $\mathcal F_\epsilon(\eta_\epsilon)=0$ are automatically smooth. So we have established the existence result. 
\end{proof}
%
%
\begin{Corollary}[Transversality]\label{cor:existence}
There is $\epsilon_0>0$, such that for all $0<\epsilon<\epsilon_0$ the gradient trajectory $\gamma_{\epsilon}$ of Theorem \ref{thm:existence} is transverse.
\end{Corollary}
\begin{proof}
Put $\eta_\epsilon:= \gamma_{\epsilon}-\gamma_0$. By Theorem \ref{thm:existence}, $\norm{\eta_\epsilon}_{W^{1,2}_\epsilon}\leq C \epsilon$.  Transversality is equivalent to proving that $d \mathcal F_\epsilon(\eta_\epsilon)$ is surjective. 
By the quadratic estimates in Corollary \ref{cor:quad1}, we have, similarly as in the proof above, $C,\epsilon_0 >0$ such that
\[
\norm{\left( d \mathcal F_\epsilon(\eta_\epsilon)- d \mathcal F_\epsilon(0) \right) \hat \eta}_{L^2_\epsilon} \leq C \epsilon^{1/2} \norm{\hat \eta}_{W^{1,2}_\epsilon}
\]
for all $\hat \eta \in W^{1,2}(\R,\R^n)$ and all $0<\epsilon<\epsilon_0$. 
Thus by Proposition \ref{prop:lin2} and Proposition \ref{prop:err1}, we get for $\upsilon \in \W$ and $0<\epsilon_0<1$ small
\begin{align*}
\norm{(d \mathcal F_\epsilon(\eta_\epsilon)-D_\epsilon) D_\epsilon^\ast \upsilon}_{L^2_\epsilon}\hspace{-0.2em} &\leq \hspace{-0.2em}\norm{(d \mathcal F_\epsilon(\eta_\epsilon)-d \mathcal F_\epsilon(0)) D_\epsilon^\ast \upsilon}_{L^2_\epsilon} + \norm{(d \mathcal F_\epsilon(0)-D_\epsilon) D_\epsilon^\ast \upsilon}_{L^2_\epsilon}\\
&\leq C \epsilon^{1/2} \norm{D_\epsilon^\ast \upsilon}_{W^{1,2}_\epsilon} + \norml{E_\epsilon D_\epsilon^\ast \upsilon} \leq   C \epsilon^{1/2} \norm{D_\epsilon D_\epsilon^\ast \upsilon}_{L^2_\epsilon}.
\end{align*}
Now we simply apply Lemma \ref{lem:lin2} with $A=d \mathcal F_\epsilon(\eta_\epsilon)$ and $\beta=\frac{1}{2}$ to finish the proof.
\end{proof}
%
%
%
\subsection{Local Uniqueness}

An immediate consequence of the leg work done so far will be the following local uniqueness result (with constraints). In the next section, we will strengthen this result in Theorem $\ref{thm:loc1}$ by removing these constraints.
\begin{Proposition}[Local Uniqueness]\label{prop:loc1}
There is $0<\epsilon_0<1$ and $\mu_1>0$ such that for all $0<\epsilon<\epsilon_0$ and for any solution $\gamma$ of \eqref{eq:eq4}, with  
\[
\textbf{(C1)} \; \norm{\gamma-\gamma_0}_{L^\infty_\epsilon}+ \epsilon \norml{\gamma-\gamma_0}<\mu_1 \qquad \text{ and }\qquad \textbf{(C2)}\; \gamma-\gamma_0 \in \image D_\epsilon^\ast\cap \X,
\]
we have $\gamma=\gamma_\epsilon$, where $\gamma_\epsilon$ is the solution from Theorem \ref{thm:existence}, $D_\epsilon^\ast$ is the adjoint in \eqref{eq:lin8} and the norms are defined in \eqref{eq:funcnorm1} and \eqref{eq:funcnorm3}.
\end{Proposition}
\begin{proof}
Denote $\eta:=\gamma-\gamma_0$ and $\eta_\epsilon:=\gamma_\epsilon-\gamma_0$. Then by assumptions and Theorem \ref{thm:existence}, we have $\eta, \eta_\epsilon \in \image D_\epsilon^\ast\cap \X$, and the estimates $\normw{\eta_\epsilon} + \epsilon^{\frac{1}{2}} \normlinfty{\eta_\epsilon} \leq C_1 \epsilon$, $\normlinfty{\eta_\epsilon}\leq \mu_1$. Since $\gamma$ and $\gamma_\epsilon$ solve \eqref{eq:eq4}, $\mathcal F_\epsilon(\eta)=\mathcal F_\epsilon(\eta_\epsilon)=0$ by \eqref{eq:exis10}. So we can read off \eqref{eq:exis9} and \eqref{eq:exis11} that 
\[
F_\epsilon(\hat \eta)=D_\epsilon(\hat \eta) + E_\epsilon(\hat \eta) = - \vect{b/\epsilon \hat \zeta^2-(\dot \gamma_0)_x}{\hat\zeta^2} - (R_\epsilon(\gamma_0+\hat \eta) - d R_\epsilon(\gamma_0) \hat \eta)
\]
holds for $\hat \eta = \eta, \eta_\epsilon$. The difference of these identity gives for $\eta=(\xi, \zeta)$ and $\eta_\epsilon=(\xi_\epsilon, \zeta_\epsilon)$
\begin{align*}
F_\epsilon(\eta_\epsilon-\eta)	=  \vect{b/\epsilon}{1}(\zeta+\zeta_\epsilon)(\zeta-\zeta_\epsilon) 
								+ (R_\epsilon(\gamma) - d R_\epsilon(\gamma_0) \eta) - (R_\epsilon(\gamma_\epsilon) - d R_\epsilon(\gamma_0) \eta_\epsilon).
\end{align*}
The first term can be estimated using \eqref{eq:equiv} with $C_2>0$ equal $C_2^2 = 2/(1-\norm b^2)$,
\begin{align}
&\norm { \vect{b/\epsilon}{1}(\zeta+\zeta_\epsilon)(\zeta-\zeta_\epsilon) }_{L^2_\epsilon}	=\norm {(\zeta+\zeta_\epsilon)(\zeta-\zeta_\epsilon)}_{L^2} 
																		\leq (\norm{\zeta}_{L^\infty}+\norm{\zeta_\epsilon}_{L^\infty}) \norm{\zeta-\zeta_\epsilon}_{L^2}\nonumber\\
																		&\leq C_2 (\norm{\eta}_{L^\infty_\epsilon}+\norm{\eta_\epsilon}_{L^\infty_\epsilon}) \norm{\eta-\eta_\epsilon}_{W^{1,2}_\epsilon}
																		\leq C_2 ( C_1 \epsilon^{\frac{1}{2}} +\mu_1) \norm{\eta_\epsilon-\eta}_{W^{1,2}_\epsilon} \label{eq:loc1}
\end{align}
where the first equality follows by $\norm{(b/\epsilon,1)}_\epsilon=1$ from \eqref{eq:lin2} and the last inequality uses the properties of $\eta_\epsilon$ and $\eta$. For the terms in $R_\epsilon$, we can use the quadratic estimates in Lemma \ref{lem:quad1} as soon as $\epsilon_0, \mu_1 >0$ small such that $\normlinfty{\eta_\epsilon - \eta}\leq C_1 \epsilon_0^{1/2} +\mu_1 \leq 1$.
\begin{align}\label{eq:loc2}
&\norm{(R_\epsilon(\gamma) - d R_\epsilon(\gamma_0) \eta) - (R_\epsilon(\gamma_\epsilon) - d R_\epsilon(\gamma_0) \eta_\epsilon)}_{L^2_\epsilon} \nonumber\\
&\leq \norml{R_\epsilon(\gamma_\epsilon+\eta-\eta_\epsilon) - R_\epsilon(\gamma_\epsilon) - dR(\gamma_\epsilon) (\eta-\eta_\epsilon) } + \norml{(dR(\gamma_\epsilon)-dR(\gamma_0)) (\eta-\eta_\epsilon)}\nonumber\\
&\leq M \bigg( \normlinfty{\eta-\eta_\epsilon} \big (\epsilon \normlinfty{ \eta-\eta_\epsilon} + \norml{\eta-\eta_\epsilon} \big)+ \epsilon \norml{\eta-\eta_\epsilon}^2 \bigg) \nonumber\\
&\phantom{\leq}+M\big( \normlinfty{\eta-\eta_\epsilon} \! (\epsilon \normlinfty{ \eta_\epsilon} \! +\! \norml{\eta_\epsilon} \big) \!+ \! \norml{\eta-\eta_\epsilon} ( \normlinfty{\eta_\epsilon}\!+\!  \epsilon \norml{\eta_\epsilon} )  \big)\nonumber\\
&\leq  C_3( \epsilon^{\frac 1 2}+\mu_1) \norm{\eta_\epsilon-\eta}_{W^{1,2}_\epsilon},
\end{align}
where the third term in line three required the $L^2_\epsilon$ bound on $\eta=\gamma-\gamma_\epsilon$.  
So in total, we get by \eqref{eq:loc1} and \eqref{eq:loc2} a $C_4>0$ such that
\begin{equation}\label{eq:loc3}
\begin{aligned}
\norm{F_\epsilon(\eta_\epsilon-\eta)}_{L^2_\epsilon} 	\leq C_4( \epsilon^{\frac 1 2}+ \mu_1) \norm{\eta_\epsilon-\eta}_{W^{1,2}_\epsilon}.
\end{aligned}
\end{equation}
As $\eta_\epsilon-\eta= (\gamma_\epsilon-\gamma_0)-(\gamma-\gamma_0) \in \image D_\epsilon^\ast$, we have by \eqref{eq:crux} in Corollary \ref{cor:lin2} and \eqref{eq:loc3}
\begin{align*}
&\norm{\eta_\epsilon-\eta}_{W^{1,2}_\epsilon}				\leq \hat C \norm{F_\epsilon(\eta_\epsilon-\eta)}_{L^2_\epsilon}
													\leq \hat C  C_4( \epsilon^{\frac 1 2}+ \mu_1) \norm{\eta_\epsilon-\eta}_{W^{1,2}_\epsilon}
													\leq \frac{1}{2} \norm{\eta_\epsilon-\eta}_{W^{1,2}_\epsilon}
\end{align*}
once $\mu_1>0$ and $0<\epsilon_0<1$ such that 
$
\hat C  C_4 \epsilon_0^{\frac{1}{2}}\leq \frac{1}{4} \text{ and  } \hat C C_4 \mu_1\leq \frac{1}{4}.
$
So we conclude that $\norm{\eta_\epsilon-\eta}_{W^{1,2}_\epsilon}=0$, that is $\gamma=\gamma_\epsilon$.
\end{proof}
%
%
\section{Strong Local Uniqueness}\label{sec:four}
This section contains the steps (X) - (XII) from the introduction, each in a separate subsection. Namely, we will strengthen the local uniqueness in Proposition \ref{prop:loc1} to the strong local uniqueness in Theorem \ref{thm:loc1}. The proof needs two analytic ingredients: time-shift in section \ref{sec:shift}, and exponential decay in Section \ref{sec:unif}. $(A,b,h)$ will always satisfy Assumptions \ref{sta:trans1}. The limit solution $\gamma_0$ was defined in \eqref{eq:adia3}.
%
%
%
\subsection{Time-shift} \label{sec:shift}
We will show that if $\gamma$ is close enough to $\gamma_0$, then for a time-shift $\gamma_\tau(t)= \gamma(t+\tau)$ of $\gamma$,  $\gamma_\tau-\gamma_0$ is in the codimension one subspace $\image D_\epsilon^\ast$. 
\begin{Proposition}[Time-shift]\label{prop:shift1}
There are $\mu_2, \mu_3, \epsilon_0 >0$ such that for all $0<\epsilon<\epsilon_0$ and for all $\gamma= \gamma_0 + \eta$ with $\eta \in \X$ and
$
\norm{\eta}_{L^2_\epsilon} <\mu_2,
$
there is $\tau\in \R$ such that
\begin{align}
\gamma_\tau-\gamma_0\in \image D_\epsilon^\ast,  \quad
\abs{\tau} \leq \mu_3 \norml{\eta}  \quad \text{ and } \quad \norm{\gamma_\tau-\gamma_0}_{L^2_\epsilon} \leq \mu_3 \norm{\eta}_{L^2_\epsilon},
\label{eq:shift11}
\end{align}
where $\gamma_0$ is defined in \eqref{eq:adia3} and the norms are defined in \eqref{eq:funcnorm1}.
\end{Proposition}
\begin{proof}
We start with some preliminary definitions.

Put $w_0:=\dot \gamma_0 \in \X$. Applying time derivative on equation \eqref{eq:adia1}, we find
\begin{equation}\label{eq:shift1}
\begin{split}
D_\epsilon w_0= \vect{( \ddot \gamma_0)_x}{( \ddot \gamma_0)_z} + \vect {\frac{1}{\epsilon} A (\dot \gamma_0)_x+ \frac{2b}{\epsilon} ( \gamma_0)_z ( \dot \gamma_0)_z}{\scal{A( \gamma_0)_x}{b} + 2 ( \gamma_0)_z ( \dot \gamma_0)_z } = \vect{(\ddot \gamma_0)_x}{0}.
\end{split}
\end{equation}
Now denote by $P_\epsilon:=D_\epsilon^\ast(D_\epsilon D_\epsilon^\ast)^{-1}D_\epsilon: \X\to \X \subset \Y$ and notice that
$
P_\epsilon^2=P_\epsilon=P_\epsilon^\ast, 
$
i.e. $P_\epsilon$ is an orthogonal projection. We also have
\[
\ker P_\epsilon=\ker D_\epsilon, \ \image P_\epsilon=(\ker D_\epsilon)^{\perp_{L^2_\epsilon}}= \image D_\epsilon^\ast \cap \X.
\]
Therefore, $w_\epsilon:=w_0-P_\epsilon w_0 \in \ker P_\epsilon=\ker D_\epsilon$. By \eqref{eq:shift1}, $w_\epsilon\neq 0$. By \cite[Proposition~2.16]{Sch1}, we know that the Fredholm index $\ind(D_\epsilon)=1$ and by Proposition \ref{prop:surj}, we know that $D_\epsilon$ is surjective for $0<\epsilon<\epsilon_0$, therefore we conclude that 
\begin{equation}\label{eq:shift8}
\begin{aligned}
\spann(w_\epsilon)=\ker D_\epsilon.
\end{aligned}
\end{equation}
Therefore, given $\gamma = \gamma_0+ \eta$ for $\eta \in \X$ and $\tau \in \R$, we define the time-shift
$
\gamma_\tau:\R\to \R^n \text{ by } \gamma_\tau(t)=\gamma(t+\tau) \text{ for } t\in \R.
$
Then $\gamma_\tau - \gamma_0 \in \X$.
Furthermore, we define $\rho:\R \to \R$ by 
\begin{equation}\label{eq:shift9}
\rho(\tau)=\scal{w_\epsilon}{\gamma_\tau-\gamma_0}_{L^2_\epsilon}.
\end{equation}
By construction and \eqref{eq:shift8}, we have
\begin{equation}\label{eq:shift10}
\gamma_\tau-\gamma_0 \in \image D_\epsilon^\ast \quad \Leftrightarrow \quad \rho(\tau)=0.
\end{equation}
We now prove this result in 5 steps. 

\textbf{Step 1:} \textbf{Estimate for $\normw{w_\epsilon-w_0}$ where $w_\epsilon$ is defined above \eqref{eq:shift8}.} 
We have $w_\epsilon-w_0=-D_\epsilon^\ast (D_\epsilon D_\epsilon^\ast)^{-1}D_\epsilon w_0 \in \image D_\epsilon^\ast$ 
and using $w_\epsilon\in\ker D_\epsilon$ and \eqref{eq:shift1}, we get
\begin{equation}\label{eq:shift2}
\begin{split}
D_\epsilon(w_\epsilon-w_0)=-D_\epsilon w_0=-\vect{( \ddot \gamma_0)_x}{0}.
\end{split}
\end{equation}
Therefore, we obtain by Proposition \ref{prop:lin2}, \eqref{eq:equiv} and \eqref{eq:shift2} that
\begin{equation}\label{eq:shift3}
\begin{split}
\norm{w_\epsilon-w_0}_{W^{1,2}_\epsilon}	&\leq C_1 \left( \norm{D_\epsilon (w_\epsilon-w_0)}_{L^2_\epsilon} \right) 
										\leq C_2 \epsilon\normltwo{(\ddot \gamma_0)_x} 
										\leq C_3 \epsilon.
\end{split}
\end{equation}
%
%
\textbf{Step 2:} \textbf{Looking at $\rho(0)$ where $\rho$ was defined in \eqref{eq:shift9}.} %
We have by \eqref{eq:shift3} and by the assumption $\norml{\gamma-\gamma_0}= \norml{\eta}\leq \mu_2$ that
\begin{align}
\abs{\rho(0)}		&=\abs{\scal{w_\epsilon}{\gamma-\gamma_0}_{L^2_\epsilon}}
				\leq \norm{w_\epsilon}_{L^2_\epsilon} \norm{\gamma-\gamma_0}_{L^2_\epsilon}
				\leq \left(  \norm{w_0}_{L^2_\epsilon} + \norm{w_\epsilon-w_0}_{L^2_\epsilon}\right)\norm{\gamma-\gamma_0}_{L^2_\epsilon} \nonumber\\
				&\leq C_4(1+\epsilon) \norm{\gamma-\gamma_0}_{L^2_\epsilon}
				\leq C_5 \mu_2. \label{eq:shift4}
\end{align}
%
\textbf{Step 3:} \textbf{Looking at $\rho'(\tau)$ for $\abs{\tau}$ small.} %
\begin{align*}
\rho'(\tau)	&= \frac{d}{d \tau} \scal{w_\epsilon}{\gamma_\tau-\gamma_0}_{L^2_\epsilon} 
			= \scal{w_\epsilon}{\partial_t \gamma_\tau}_{L^2_\epsilon}\\
			&= \int_\R \partial_t(g_\epsilon(w_\epsilon, \gamma_\tau-\gamma_0)) \ d t - \scal{\dot w_\epsilon}{\gamma_\tau-\gamma_0}_{L^2_\epsilon}+\scal{w_\epsilon}{\dot \gamma_0}_{L^2_\epsilon}
			=:S_1+S_2+S_3
\end{align*}
where the second equality uses $\partial_t \gamma_\tau= \partial_\tau \gamma_\tau$ and $S_i$ stands in for the $i^{th}$ summand in the preceding expression. We now estimate each summand. 
\begin{align*}
S_1 &=\lim_{R\to \infty} \int_{-R}^R \partial_t g_\epsilon(w_\epsilon, \gamma_\tau-\gamma_0) \ d t =\lim_{R\to \infty} [ g_\epsilon (w_\epsilon, \gamma_\tau-\gamma_0)]_{-R}^R =0
\end{align*}
where we used a continuous representative of $\gamma_\tau-\gamma_0$. Cf. \cite[Corollary~VIII.8.]{Bre1}.
\begin{align}
\abs{S_2}&= \abs{\scal{\dot w_\epsilon}{\gamma_\tau-\gamma_0}_{L^2_\epsilon}} 	\leq \left( \norm{w_\epsilon-w_0}_{W^{1,2}_\epsilon} + \norm{w_0}_{W^{1,2}_\epsilon} \right) \norm{\gamma_\tau-\gamma_0}_{L^2_\epsilon} \nonumber\\
																&\leq C_6 (1+ \epsilon ) \norm{\gamma_\tau-\gamma_0}_{L^2_\epsilon} \nonumber 
																\leq C_6 (1+\epsilon ) ( \norm{\gamma-\gamma_0}_{L^2_\epsilon} +\norm{(\gamma_0)_\tau-\gamma_0}_{L^2_\epsilon}) \\\																&\leq C_7 ( \norm{\gamma-\gamma_0}_{L^2_\epsilon} + \abs{\tau}) \leq C_7(\mu_2+\abs{\tau}) \label{eq:shift6}
\end{align}
for $\epsilon_0<1$ and $C_6 \geq 1$ where the first line uses Cauchy--Schwarz, the second one uses \eqref{eq:shift3}, the third one uses the invariance of the $L^2_\epsilon$ norm under time-shift and the last line uses the assumption and the exponential decay to $(0,\pm 1)$ of $\gamma_0$ at its ends.
\begin{align*}
\abs{S_3} &= \scal{w_\epsilon}{w_0}_{L^2_\epsilon}	\geq \norm{w_0}^2_{L^2_\epsilon}-\abs{\scal{w_\epsilon-w_0}{w_0}_{L^2_\epsilon}} 
										\geq  \norm{w_0}_{L^2_\epsilon} ( \norm{w_0}_{L^2_\epsilon}- \norm{w_\epsilon-w_0}_{L^2_\epsilon}) \\
										&\geq \norm{w_0}_{L^2_\epsilon} ( \norm{w_0}_{L^2_\epsilon} - C_3 \epsilon ) 
										\geq C_8 \norm{(w_0)_z}_{L^2}\left( \norm{(w_0)_z}_{L^2}-C_3 \epsilon  \right) 
										\geq 2 c_0 >0
\end{align*}
where $w_0:= \dot \gamma_0$ and $c_0>0$ can be chosen independent of $\epsilon$ if $ 0<\epsilon<\epsilon_0< \norm{(w_0)_z}_{L^2}/C_3$. 
In these inequalities, we used \eqref{eq:shift3} and \eqref{eq:equiv}.
Combining these estimates, we end up with
\begin{equation}\label{eq:shift5}
\abs{\rho'(\tau)} \geq 2 c_0 - C_7(\abs{\tau}+\mu_2).
\end{equation}
%
\textbf{Step 4:} \textbf{Finding $\tau$ such that $\rho(\tau)=0$.} %
Assume that $\mu_2<\frac{1}{2 C_7} c_0$ and $\abs{\tau} \leq \frac{1}{2 C_7} c_0$, then by \eqref{eq:shift5}, we have
\begin{align*}
\abs{\rho'(\tau)} \geq 2 c_0 - C_7(\abs{\tau}+\mu_2) \geq c_0.
\end{align*}
Restricting $\mu_2$ further to $\mu_2\leq \frac{1}{2C_7C_5} c_0^2$, we have that $\frac{\mu_2 C_5}{c_0} \leq  \frac{1}{2 C_7} c_0$, which is enough by \eqref{eq:shift4}, and the intermediate value theorem to guarantee a zero of $\rho$ in $[-\frac{\mu_2 C_5}{c_0}, \frac{\mu_2 C_5}{c_0}]$.

\textbf{Step 5:} \textbf{Conclusion.} %
By \eqref{eq:shift10}, we found $\tau$ such that $\gamma_\tau-\gamma_0\in \image D_\epsilon^\ast$.
By \eqref{eq:shift6}, we get
\begin{align*}
\norm{\gamma_\tau-\gamma_0}_{L^2_\epsilon} &\leq C_7( \norm{\gamma-\gamma_0}_{L^2_\epsilon} + \abs{\tau})
\end{align*}
and furthermore the zero $\tau$ of $\gamma$ was found to have $\abs{\tau} \leq \mu_2 C_5/c_0$. However the penultimate line of \eqref{eq:shift4} gives a better estimate of
\begin{equation*}
\begin{split}
\abs{\tau} \leq \frac{2C_4 \norm{\gamma-\gamma_0}_{L^2_\epsilon}}{c_0}.
\end{split}
\end{equation*}
Meaning that $\mu_3:= \max \left( C_7\left( 1+\frac{2 C_4}{c_0} \right), \frac{2C_4}{c_0} \right) $ does the job.
\end{proof}
%
%
%
\subsection{Exponential Decay}\label{sec:unif}

This section is dedicated to proving that solutions of our equation have exponential decay (uniform in $\epsilon$) to the critical points at the ends as stated in Proposition \ref{prop:unif} below. The proof is adapted from \cite[Proposition 2.10]{Sch1}.
\begin{Proposition} \label{prop:unif}
There is $\epsilon_0>0$ such that for all $0<\epsilon<\epsilon_0$,
a solution $\gamma$ of \eqref{eq:eq4} has exponential decay at both ends with uniform rate $\eta=(1-\norm b^2)$ i.e. there is $t_0 \in \R$ and $C>0$ such that for all $t\geq t_0$
\begin{align}
\epsilon^2 \norm {\gamma_x(t)}^2 + \abs{\gamma_z(t)-1}^2 &\leq C \left(  \epsilon^2 \norm {\gamma_x(t_0)}^2 + \abs{\gamma_z(t_0)-1}^2 \right) e^{-\eta(t-t_0)}, \notag \\
\epsilon^2 \norm {\gamma_x(-t)}^2 + \abs{\gamma_z(-t)+1}^2 &\leq C \left(  \epsilon^2 \norm {\gamma_x(-t_0)}^2 + \abs{\gamma_z(-t_0)+1}^2 \right) e^{-\eta(t-t_0)}. \label{eq:unif4}
\raisetag{2.5\baselineskip}
\end{align}
Furthermore, \eqref{eq:unif4} holds whenever $t_0$ was chosen such that for every $\abs t \geq t_0$ estimate \eqref{eq:unif2} for $\delta=\frac{1}{2}$ below holds.
\end{Proposition}
\begin{proof}
Take $0<\epsilon<1$. Let $\gamma=(\gamma_x, \gamma_z)$ be a solution of the differential equations in \eqref{eq:eq4}. We recall that $(f_\epsilon, g_\epsilon)$ in \eqref{eq:eq5} is a gradient pair for \eqref{eq:eq4}.  

\vspace{-\parskip}Now also assume that $\lim_{t\to \infty} \gamma(t) = (0,1)=p_+$. Then there is $t_0>0$ such that for all $t\geq t_0$, $\gamma(t)~\in~B_{1}(0)~\times~B_2(0)\subset \R^{n-1}\times \R$. Thus for $t\geq t_0$, we get from $h_{0,(0,0)}=0$
\begin{equation}\label{eq:unif5}
\norm{h_{\epsilon^2, (\epsilon^2 \gamma_x(t), \epsilon \gamma_z(t))}} \leq C_1 \epsilon \norm{h}_{C^1(B_4(0))} \leq M \epsilon.
\end{equation}
We look at the function $\alpha:\R \to \R_{\geq 0}$ defined by
\begin{equation}\label{eq:unif3}
\alpha(t):=1/2 g_{\epsilon,\gamma(t)} \left(\gamma(t) -p_+, \gamma(t) - p_+ \right).
\end{equation}
As $\gamma$ is a gradient trajectory, the derivative of $\alpha$ is given by 
\[
\alpha'=-(d f_\epsilon\circ \gamma) \left( \gamma - p_+ \right) =-  \epsilon {(\gamma_x)}^\top A \gamma_x -(\gamma_z)^3+(\gamma_z)^2+(\gamma_z)-1.
\]
By the equations in \eqref{eq:eq4}, its second derivative is given by
\begin{align*}
&\alpha''	=-2 (\gamma_x)^\top A ( \epsilon \, \dot \gamma_x) - (3(\gamma_z)^2 -2\gamma_z - 1) \, \dot \gamma_z\\
& = 2 (\gamma_x)^\top A ( A \gamma_x + b(\gamma_z^2-1) - \epsilon R^x_\epsilon(\gamma)) + (3(\gamma_z)^2 -2\gamma_z - 1) ( \scal{A\gamma_x}{b}+\gamma_z^2-1 - R^z_\epsilon(\gamma))\\
&=2 \norm{A\gamma_x}^2 + 8 \abs{\gamma_z-1}^2 \frac{(3\gamma_z+1)(\gamma_z+1)}{8} + (5\gamma_z +3)(\gamma_z-1) \scal b {A\gamma_x} \\
&\phantom{=}- (2 (\gamma_x)^\top A, 3 \gamma_z^2 - 2 \gamma_z -1) \vect{\epsilon R^x_\epsilon(\gamma)}{R^z_\epsilon(\gamma)}
\end{align*}
Now we can estimate using the definition of $R_\epsilon$ in \eqref{eq:eq4} and \eqref{eq:unif5} for $t\geq t_0$ that
\begin{align*}
&\abs{(2 (\gamma_x)^\top A, 3 \gamma_z^2 - 2 \gamma_z -1) \vect{\epsilon R^x_\epsilon(\gamma)}{R^z_\epsilon(\gamma)}} = \abs{(2 (\gamma_x)^\top A, 3 \gamma_z^2 - 2 \gamma_z -1) h_{\epsilon^2, (\epsilon^2 \gamma_x, \epsilon \gamma_z)}\! \vect{A \gamma_x}{\gamma_z^2-1}}\\
&\leq M\epsilon \left( \norm{A\gamma_x} + \abs{(3\gamma_z+1)(\gamma_z-1)} \right)\left( \abs{A \gamma_x} + \abs{(\gamma_z+1)(\gamma_z-1)}\right)\\
&\leq M\epsilon \left(4 \norm{A\gamma_x}^2 + \left( \abs{(\gamma_z+1)(3 \gamma_z+1)} +(2 \abs{\gamma_z+1} +\abs{3 \gamma_z+1})^2\right) \abs{\gamma_z-1}^2 \right)  \\
&\leq M \epsilon^{1/2} ( \norm{\gamma_x}^2 + \abs{\gamma_z-1}^2),
\end{align*}
for all $0<\epsilon<\epsilon_0$ for $\epsilon_0$ small. On the other hand, we have
\[
8\abs{\frac{5\gamma_z +3}{8}} \abs{\gamma_z-1} \abs{\scal b {A\gamma_x}} \leq 4 \frac{\norm b^2 \norm{A \gamma_x}^2}{1+\norm b^2} + 4 (1+\norm b^2) \abs{\frac{5\gamma_z +3}{8}}^2 \abs{\gamma_z-1}^2.
\]
Let $\kappa>0$ be such that $\norm{Ax}\geq \kappa \norm x$, for all $x\in \R^{n-1}$ and choose $0<\delta<1$. We have
\begin{equation}\label{eq:unif2}
\abs{\frac{(5\gamma_z +3)^2}{8^2} -1} < \frac{\delta}{2}, \quad \abs{\frac{(3\gamma_z+1)(\gamma_z+1)}{8}-1} < \frac{\delta}{2}, \quad (\text{and }\norm{\gamma_x(t)} \leq 1) 
\end{equation}
for all $t \geq t_0$ for $t_0 = t_0(\delta)$ big. This choice of $t_0$ is possible, since $\lim_{t\to \infty} \gamma_z(t) = 1$. Putting all these estimates together, we get for $t\geq t_0$
\begin{align*}
\alpha''	& \geq \left(2 \kappa^2 \frac{1-\norm b^2}{1+ \norm b^2} -M \epsilon^{1/2} \right)\norm {\gamma_x}^2 + (4 (1-\norm b^2) (1-\frac{\delta}{2}) - M \epsilon^{1/2}) (\gamma_z-1)^2 \\
		&\geq (1-\norm b^2) \left(\frac{\kappa^2}{(1+\norm b^2) \epsilon^2} (\epsilon^2 \norm{\gamma_x}^2) + 4(1-\delta) (\gamma_z-1)^2 \right)
\end{align*}
whenever $0<\epsilon<\epsilon_0$ for $\epsilon_0$ small. For $\epsilon_0$ even smaller, $\kappa^2\geq 4(1-\delta)(1+\norm b^2) \epsilon_0^2$ and Lemma \ref{lem:equiv1} comparing norms holds with $M$ as in \eqref{eq:unif5} with factor $4/(1-\norm b^2)$. Thus we get for $t\geq t_0$, by definition of $\alpha$ in \eqref{eq:unif3},
\begin{align*}
\alpha''(t) 	&\geq 4 (1-\norm b^2) (1-\delta) (\epsilon^2 \norm{\gamma_x(t)}^2 + \abs{\gamma_z(t) -1}^2) \\
		&\geq (1-\norm b^2)^2 (1-\delta)\; g_{\epsilon,\gamma(t)} \left(\gamma(t) -p_+, \gamma(t) - p_+ \right) \geq 2 (1-\norm b^2)^2(1-\delta) \alpha(t).
\end{align*}
We now set $\eta_\delta := (1-\norm b^2) \sqrt{ 2 \, (1-\delta)}$ and define $\alpha_0(t):=\alpha(t_0) \exp(-\eta_\delta(t-t_0))$ and $\Delta(t) := \alpha(t)-\alpha_0(t)$. $\Delta$ has the following properties
\begin{align*}
\Delta''(t) \geq \eta_\delta^2 \Delta(t) & \text{ for all } t \geq t_0, && \Delta(t_0)=0, && \lim_{t\to \infty} \Delta(t) =0.
\end{align*}
Due to these properties, no positive maximum can be attained on $[t_0, \infty)$. Therefore, we get $\Delta(t)\leq 0$ for all $t\geq t_0$ i.e. by Lemma \ref{lem:equiv1} there are $C_1, C_2>0$, such that for $t\geq t_0$
\begin{align*}
C_2(\epsilon^2 \norm {\gamma_x(t)}^2 + \abs{\gamma_z(t)-1}^2 )	&\leq \alpha(t) \leq  \alpha(t_0) e^{-\eta_\delta(t-t_0)} \\& \leq C_3 \left(  \epsilon^2 \norm {\gamma_x(t_0)}^2 + \abs{\gamma_z(t_0)-1}^2 \right) e^{-\eta_\delta(t-t_0)}.
\end{align*}
Repeating the argument with $\tilde A:=-A$, $\tilde h_{\lambda,(x,z)} := h_{\lambda,(-x,-z)}$ and $\tilde \gamma(t):=-\gamma(-t)$, which again fulfil $\eqref{eq:eq4}$ with $(\tilde A, b, \tilde h)$, we get $t_0$ and $\epsilon_0>0$ as in the argument before such that for all $t\geq t_0$,
\[
\epsilon^2 \norm {\gamma_x(-t)}^2 + \abs{\gamma_z(-t)+1}^2 \leq C_4 \left(  \epsilon^2 \norm {\gamma_x(-t_0)}^2 + \abs{\gamma_z(-t_0)+1}^2 \right) e^{-\eta_\delta(t-t_0)}. \qedhere
\]
\end{proof}
%
\subsection{Strong Local Uniqueness}

By the existence Theorem \ref{thm:existence}, there is a solution $\gamma_\epsilon$ of \eqref{eq:eq4} for every $0<\epsilon<\epsilon_0$. We now prove that if any solution $\gamma$ of \eqref{eq:eq4} for $\epsilon>0$ small is contained in a neighbourhood of the limit solution $\gamma_0$ in \eqref{eq:adia3}, then $\gamma$ has to be already the solution $\gamma_\epsilon$ (up to time-shift). This strengthens the local uniqueness of Proposition~\ref{prop:loc1}.  
\begin{Theorem}[Strong Local Uniqueness]\label{thm:loc1} Fix $\nu \in (0, 1/2)$ and $R>1$. 
There is $\epsilon_0>0$, such that for all $0<\epsilon<\epsilon_0$, a solution $\gamma$ of \eqref{eq:eq4} with
\begin{equation}\label{eq:uniq13}
\normlinftyy{A \gamma_x + b(\gamma_z^2-1)} \leq \epsilon^\nu, \quad \normlinftyy{\gamma_z} \leq R,
\end{equation}
then $\gamma$ is equal to $\gamma_\epsilon$ up to a time-shift, where $\gamma_\epsilon$ is the solution from Theorem \ref{thm:existence}.
\end{Theorem}

\textbf{Note: } Equation \eqref{eq:adia1} reads $A(\gamma_0)_x + b ((\gamma_0)_z^2-1)=0$. So solutions $\gamma$ of \eqref{eq:eq4} with \eqref{eq:uniq13} are understood to be in a neighbourhood of $\gamma_0$. At this point, it is not clear that there is a single solution that fulfils \eqref{eq:uniq13}. However, the topological result in Section \ref{sec:five} will prove that every solution of \eqref{eq:eq4} satisfies \eqref{eq:uniq13} which leads to uniqueness in Theorem \ref{thm:goal2} 
\begin{proof}
Suppose we are given solutions $\gamma_i$ of \eqref{eq:eq4} for a sequence $\epsilon_i >0$ with $\epsilon_i \to 0$ and
\begin{equation}\label{eq:uniq2}
\normlinftyy{A (\gamma_i)_x + b((\gamma_i)_z^2-1)} \leq (\epsilon_i)^\nu, \quad \normlinftyy{(\gamma_i)_z} \leq R.
\end{equation}
We want to show that there is $I>0$ depending only on the sequence $\epsilon_i$ such that for all $i\geq I$, there is $\tau_i$ such that 
$
(\gamma_i)_{\tau_i} = \gamma_{\epsilon_i},
$ 
where $\gamma_\epsilon$ is the solution in Theorem \ref{thm:existence}.
Since $\lim_{t\to \pm \infty} (\gamma_i)_z(t) = \pm 1$ and $\gamma_i$ is continuous, we can assume up to time-shift that
\begin{equation}\label{eq:uniq1}
\begin{split}
(\gamma_i)_z(0)=0
\end{split}
\end{equation}
for all $i \in \N$. We note that the second line of \eqref{eq:eq4} reads
\begin{equation}\label{eq:uniq12}
\begin{aligned}
(\dot \gamma_i)_z = - \scal{b}{A (\gamma_i)_x + b((\gamma_i)_z^2-1)}+ (1-\norm b^2)(1-(\gamma_i)_z^2) + R_{\epsilon_i}^z(\gamma_i).
\end{aligned}
\end{equation}
We also have $\normlinftyy{(\gamma_i)_x} \leq C_1$ and \eqref{eq:unif5} gives $\norm{\epsilon_i R^x_{\epsilon_i}(\gamma_i)}+\abs{R^z_{\epsilon_i}(\gamma_i)}\leq C_2 \epsilon_i$.  

%
%
\textbf{Step 1: For $T>0$, $\gamma_i$ converges uniformly on $[-T,T]$ to the limit solution $\gamma_0$.}

Fix $T>0$. We have $\norm b<1$ and for $I$ big that $\epsilon_i<1$. So on $[-T,T]$
\begin{align}
&\abs{ {(\dot \gamma_i)_z}-(\dot \gamma_0)_z} 	\nonumber\\
				&= \abs{R_{\epsilon_i}^z(\gamma_i)-  \scal{b}{A (\gamma_i)_x + b((\gamma_i)_z^2-1)}  +(1-\norm b^2)(1-(\gamma_i)_z ^2) -(1-\norm b^2)(1-(\gamma_0)_z^2 )} \notag\\
				&\leq (C_2+1) \epsilon_i^\nu +  (1-\norm b^2)\abs{(\gamma_i)_z +(\gamma_0)_z } \abs{(\gamma_i)_z -(\gamma_0)_z }\nonumber\\
				& \leq (C_2+ 1) \epsilon_i^\nu + 2R \abs{(\gamma_i)_z -(\gamma_0)_z },\label{eq:uniq10}
\end{align}
where the bound on $\gamma$, \eqref{eq:uniq12}  for $(\dot \gamma_i)_z$ and \eqref{eq:adia2} for $(\dot \gamma_0)_z$.
Next, let $N\in \N$ be such that $\frac{N}{4R}>T$. We will estimate the distance between $(\gamma_i)_z$ and $(\gamma_0)_z$ on subintervals that cover $[-T,T]$.  

\vspace{-\parskip} Define $Q_n^i:= \norm{(\gamma_i)_z-(\gamma_0)_z}_{L^\infty([\frac{n-1}{4R}, \frac{n}{4R}])}$ for $n\geq 1$ and $Q^i_0=\abs{(\gamma_i)_z(0) - (\gamma_0)_z(0)}$. By \eqref{eq:uniq1}, we have that $\abs{Q^i_0} = 0$. We also have for $n\geq 1$ that
\begin{align*}
Q_n^i &\leq \abs{(\gamma_i)_z\left( \frac{n-1}{4R} \right)-(\gamma_0)_z\left( \frac{n-1}{4R} \right)} + \int_{\frac{n-1}{4R}}^{\frac{n}{4R}}\abs{(\dot \gamma_i)_z(t)-(\dot \gamma_0)_z(t)} \, d t\\
	& \leq Q^i_{n-1} +\frac{1}{4R} (C_2+1)\epsilon_i^\nu + \frac{1}{2} Q_n^i  \quad\Rightarrow Q_n^i \leq 2 Q^i_{n-1} +\frac{1}{2R} (C_2+1) \epsilon_i^\nu	
\end{align*}
where we used \eqref{eq:uniq10} in the second line.
Now setting $Q_n^i = a_n \epsilon_i^{\nu}$, we get that $a_n$ verifies $a_0\leq 1$, $a_n \leq 2a_{n-1}+1$ as soon as $i \geq I$ such that $(C_2+1) \epsilon_i^{\nu}/(2R)<1$. Such a recursive inequality implies that $a_n \leq \sum_{i=0}^n 2^i = 2^{n+1}-1$ and so
$
Q_n^i \leq  2^{n+1} \epsilon_i^\nu .
$
The same argument works for $t<0$, and so
\[
\norm{(\gamma_i)_z-(\gamma_0)_z}_{L^\infty([-T,T])} \leq 2^{N+1} \epsilon_i^\nu.
\]
As $\epsilon_i \to 0$, this inequality reads $(\gamma_i)_z$ converges uniformly to $(\gamma_0)_z$ on $[-T,T]$. 
Combining the uniform convergence of $(\gamma_i)_z$ to $(\gamma_0)_z$ and that of $A (\gamma_i)_x-  b (1- (\gamma_i)_z^2)$ to zero on $[-T,T]$, we get that $(\gamma_i)_x$ converges uniformly to $( \gamma_0)_x:=A^{-1}   b (1-(\gamma_0)_z^2)$. 
Since all the $\gamma_i$ solve \eqref{eq:eq4}, we also get that
\begin{equation}\label{eq:uniq9}
\begin{aligned}
\vect{\epsilon (\dot \gamma_i)_x}{(\dot \gamma_i)_z} \text{ converges uniformly on }[-T,T] \text{ to }\vect{0}{(1-\norm b^2) (1-(\gamma_0)_z^2)}. 
\end{aligned}
\end{equation}
%
%
\textbf{Step 2: Energy convergence and exponential decay on the ends.} 

Define the one dimensional gradient pair $f_0(v) := \frac{1}{3}v^3-v$, $g_0:= \fracb $ and its associated energy functional,
\[
\energy_0(\sigma):= \frac{1}{2} \int_\R \fracb \abs{\dot \sigma}^2 \, d t + \frac 12 \int_\R (1-\norm b^2) \abs{1-  \sigma^2}^2 \, d t,
\]
where $\sigma:\R \to \R$ is smooth with $\lim_{t\to\pm \infty} \sigma(t)=\pm 1$.
A minimum of this energy functional is $(\gamma_0)_z(t)=\tanh((1-\norm b^2)t)$, and $\energy_0(z_0)=\frac{4}{3}= f_{\epsilon_i}(p_-) - f_{\epsilon_i}(p_+) =  \energy_{\epsilon_i}(\gamma_i)$, where $(f_\epsilon, g_\epsilon)$ is the gradient pair in \eqref{eq:eq5}.

Since \eqref{eq:eq4} is the negative gradient field of the pair $(f_{\epsilon_i},g_{\epsilon_i})$, we get the energy functional
\begin{align}
\energy_{\epsilon_i}\left( \left. \gamma_i\right|_{[-T,T]} \right)	& = \int_{-T}^{T} g_{\epsilon_i, \gamma_i}(\dot \gamma_i, \dot \gamma_i) \; dt= - \int_{-T}^{T} \left( d f_{\epsilon_i}(\gamma_i) \dot \gamma_i \right)\, d t  \label{eq:uniq11}\\
													&= - \int_{-T}^{T}   \left({(\gamma_i)_x}^\top A (\epsilon (\dot \gamma_i)_x) + ((\gamma_i)_z^2-1) (\dot \gamma_i)_z \right)\,d t 
													\, \xrightarrow{i \to \infty}  \energy_0\left(\left. (\gamma_0)_z\right|_{[-T,T]}\right) \nonumber
\end{align}
where we used uniform convergence \eqref{eq:uniq9} established in Step 1.

Next fix $0<\rho<\frac{2}{9}$. Take $T:=T(\rho)>0$ such that $\energy_0\big(\left. (\gamma_0)_z\right|_{[-T,T]}\big)\geq  \frac{4}{3}-\rho$. Then by \eqref{eq:uniq11}, there is $I:=I(T,\rho)>0$ such that for all $i\geq I$, we get 
\[
\abs{\energy_{\epsilon_i}\left( \left. \gamma_i\right|_{[-T,T]} \right)  - \energy_0\left( \left. (\gamma_0)_z\right|_{[-T,T]}  \right)} < \rho.
\]
Therefore $\energy_{\epsilon_i}(\left. \gamma_i\right|_{[-T,T]} ) \geq \frac{4}{3}- 2\rho$.
Thus, we have
\begin{align*}
\frac{4}{3}=\energy_{\epsilon_i}(\gamma_i)	&= \energy_{\epsilon_i}\left(\res{\gamma_i}{(-\infty, -T)}\right) +\energy_{\epsilon_i}\left(\res{\gamma_i}{[-T, T]}\right)+ \energy_{\epsilon_i}\left(\res{\gamma_i}{(T, \infty)}\right)\\
									&\geq \frac{4}{3}-2\rho  + \energy_{\epsilon_i}\left(\res{\gamma_i}{(-\infty, -T)}\right) + \energy_{\epsilon_i}\left(\res{\gamma_i}{(T, \infty)}\right)
\end{align*} 
and so
\begin{align*}
\energy_{\epsilon_i}(\res{\gamma_i}{(-\infty, -T)}) < 2\rho &\Rightarrow f_{\epsilon_i}(\gamma_i(-t)) > \frac{2}{3} -2\rho\\
\energy_{\epsilon_i}(\res{\gamma_i}{(T, \infty)}) < 2\rho &\Rightarrow f_{\epsilon_i}(\gamma_i(t)) < -\frac{2}{3} +2\rho
\end{align*}
for all $t\geq T$.
Next, since $\normlinftyy{(\gamma_i)_x}\leq C_1$, we have the first summand of $f_{\epsilon_i}$ bounded by 
\begin{align*}
\abs{\frac{1}{2} \epsilon_i (\gamma_i)_x^\top A (\gamma_i)_x} \leq  \frac 12 \epsilon_i \norm A C_1^2  <\rho,  
\end{align*}
whenever $i\geq I$ for $I$ maybe bigger. Thus by \eqref{eq:eq5}, we have for $t\geq T$ 
\[
-\frac{2}{3} +2\rho \geq f_{\epsilon_i}(\gamma_i ) \geq  \frac{1}{3} (\gamma_i)_z^3  - (\gamma_i)_z  -\rho.
\]
Thus, for all $t\geq T$,
\begin{equation}\label{eq:uniq3}
\begin{split}
\frac{1}{3}(\gamma_i)_z^3 -(\gamma_i)_z  \leq - \frac{2}{3}+3 \rho \iff ((\gamma_i)_z -1)^2((\gamma_i)_z +2)\leq 9 \rho.
\end{split}
\end{equation}
From the first inequality in \eqref{eq:uniq3} and $\rho <\frac{2}{9} \iff - \frac{2}{3}+3 \rho<0$, we conclude that $\frac{1}{3}(\gamma_i)_z  ((\gamma_i)_z^2 -3) <0$ for $t\geq T$. This implies
$
(\gamma_i)_z <-\sqrt{3},\text{ or  }0\leq (\gamma_i)_z  \leq \sqrt{3}
$
for $t\geq T$. 
As $\lim_{t\to \infty} (\gamma_i)_z =1$, we conclude by continuity that for $t\geq T$,   
\[
0\leq (\gamma_i)_z  \leq \sqrt{3}.
\]
Therefore by the second inequality of \eqref{eq:uniq3} due to $(\gamma_i)_z +2\geq 2$, we conclude that 
\[
((\gamma_i)_z -1)^2 \leq \frac{9 \rho}{2} \iff \abs{(\gamma_i)_z -1} \leq \frac{3\sqrt2}{2} \sqrt \rho
\]
for $t\geq T$. We combine this with $\norm b <1$ and \eqref{eq:uniq2} to get
\begin{align*}
\norm{\gamma_z -1 }_{L^\infty([T,\infty))}\leq \frac{3\sqrt2}{2} \sqrt \rho, \; \norm{\gamma_x }_{L^\infty([T,\infty))} \leq \norm A^{-1}\!\bigg( \epsilon^{\nu}+ (1+ R) \bigg(\frac{3\sqrt 2}{2} \sqrt \rho \bigg).\bigg) 
\end{align*}
Thus for $i \geq I$ for $I$ maybe bigger, we end up with some constant $C_3>0$ such that
$
\norm{\gamma  - p_+}_{L^\infty([T,\infty))} \leq C_3 \sqrt \rho,
$ where $p_+=(0,1)$ as always.

Hence for $\rho>0$ small, \eqref{eq:unif2} is verified for $t\geq T(\rho)$ and all $i\geq I(T,\rho)$. Therefore, we have by Proposition \ref{prop:unif} and \eqref{eq:equiv} that there is $C_4>0$ such that for all $t \geq T$,
\begin{equation}\label{eq:uniq6}
\begin{split}
\abs{\gamma_i(t)-p_+}_{\epsilon_i} 	& \leq C_4 e^{-(1-\norm b^2)(t-T)}.
\end{split}
\end{equation}
On par with this, we may choose $C_4>0$ even bigger, such that for $t \geq T$
\begin{equation}\label{eq:uniq7}
\begin{split}
\abs{\gamma_0(t)-p_+}_{\epsilon_i} 		&\leq C_4  e^{-(1-\norm b^2)(t-T)}.
\end{split}
\end{equation}
Similar estimates as \eqref{eq:uniq6} and \eqref{eq:uniq7} hold for $t \leq -T$ by the same arguments.

%
\newpage
\textbf{Step 3: Conclusion.} %

Take $\tilde T>T$, we can estimate by using \eqref{eq:uniq6} and \eqref{eq:uniq7}
\begin{align}
\norm{\gamma_i-\gamma_0}_{L_{\epsilon_i}^\infty\left( [\tilde T,\infty) \right)}  	&\leq \norm{\gamma_i-p_+}_{L_{\epsilon_i}^\infty\left( [\tilde T,\infty) \right)}+\norm{p_+-\gamma_0}_{L_{\epsilon_i}^\infty\left( [\tilde T,\infty) \right)} \notag\\
																&\leq 2C_4  e^{-(1-\norm b^2)(\tilde T-T)}, \label{eq:uniq4}\\
\norm{\gamma_i-\gamma_0}_{L_{\epsilon_i}^2\left( [\tilde T,\infty) \right)}  &\leq \norm{\gamma_i-p_+}_{L_{\epsilon_i}^2\left( [\tilde T,\infty) \right)}+\norm{p_+-\gamma_0}_{L_{\epsilon_i}^2\left( [\tilde T,\infty) \right)}\notag\\
		&\leq \frac{2 C_4}{1-\norm b^2}  e^{-(1-\norm b^2)(\tilde T-T)} \notag,				
\end{align}
with norms as in \eqref{eq:funcnorm1}.
Similar estimates as \eqref{eq:uniq4} can be established on $(-\infty, \tilde T]$.

Now recall the constants $\delta_1$ from Proposition \ref{prop:loc1}, $\delta_2, \delta_3$  from Proposition \ref{prop:shift1}. Then by \eqref{eq:uniq4}, we may choose $\tilde T$ such that
\begin{align*}
\norm{\gamma_i-\gamma_0}_{L_{\epsilon_i}^\infty\left( (-\infty,-\tilde T] \cup [\tilde T,\infty) \right)} < \frac{\delta_1}{8}, \quad
\norm{\gamma_i-\gamma_0}_{L_{\epsilon_i}^2\left( (-\infty,-\tilde T] \cup [\tilde T,\infty) \right)}  < \frac{\tilde \delta_2}{2} 
\end{align*}
where $\tilde \delta_2 := \min \left( \delta_2 , \frac{\delta_1}{ 8 \delta_3} \right)$.
Also by norm equivalence \eqref{eq:equiv}, there is $C_5>0$ such that 
\begin{align*}
\norm{\gamma_i-\gamma_0}_{L_{\epsilon_i}^2\left(  [-\tilde T, \tilde T] \right)} \leq C_5 \tilde T  \norm{\gamma_i-\gamma_0}_{L_{\epsilon_i}^\infty\left(  [-\tilde T, \tilde T] \right)}.
\end{align*}
We know by Step 1, that we have uniform convergence of $\gamma_i$ to $\gamma_0$ on $[-\tilde T, \tilde T]$ and so for $I$ maybe even bigger, we get for all $i\geq I$
\begin{align*}
 \norm{\gamma_i-\gamma_0}_{L_{\epsilon_i}^\infty\left(  [-\tilde T, \tilde T] \right)} < \min\left( \frac{\delta_1}{8}, \frac{\tilde \delta_2}{2 C_5  \tilde  T } \right).
\end{align*}
Combining these estimates, we get for $i\geq I$ that
\begin{equation}\label{eq:uniq8}
\begin{split}
 \norm{\gamma_i-\gamma_0}_{L_{\epsilon_i}^\infty} \leq \frac{\delta_1}{4}, \quad
 \norm{\gamma_i-\gamma_0}_{L_{\epsilon_i}^2} \leq \min\left( \delta_2, \frac{\delta_1}{8 \delta_3}\right).
\end{split}
\end{equation}
By Lemma \ref{lem:equiv1} and \eqref{eq:equiv}, we have $\norml{\dot \gamma_i}^2 \leq C \energy_{\epsilon_i}(\gamma_i) = C 4/3$ and so $\gamma_i-\gamma_0 \in \X$.
With this in hand, we can apply Proposition \ref{prop:shift1} due to the second inequality in \eqref{eq:uniq8}. This gives us $\tau_i>0$ for $i\geq I$ such that the shifted solutions $(\gamma_i)_{\tau_i}$ have the property that $(\gamma_i)_{\tau_i}- \gamma_0 \in \image D_{\epsilon_i}^\ast$. Proposition \ref{prop:shift1} also gives us the bound
\[
\abs{\tau_i} \leq \delta_3 \norm{\gamma_i-\gamma_0}_{L_{\epsilon_i}^2} \leq \frac{\delta_1}{8}.
\]
To finish the proof, we would like to apply the local uniqueness in Proposition \ref{prop:loc1}. For this we need to estimate $ (\gamma_i)_{\tau_i}- \gamma_0$ in the $L^\infty_{\epsilon_i}$ and $L^2_{\epsilon_i}$ norm. Namely,
\begin{align*}
\norm{(\gamma_i)_{\tau_i}- \gamma_0}_{L^\infty_{\epsilon_i}} &\leq \norm{(\gamma_i)_{\tau_i}- (\gamma_0)_{\tau_i}}_{L^\infty_{\epsilon_i}}+\norm{(\gamma_0)_{\tau_i}- \gamma_0}_{L^\infty_{\epsilon_i}}\\
& \leq \frac{\delta_1}{2} + \sup_{\R} \norm{\dot \gamma_0}_{\epsilon_i} \abs{\tau_i} \leq  \frac{\delta_1}{4}+ 2  \frac{\delta_1}{ 8} =\frac{\delta_1}{2}\\
\epsilon_i \norm{(\gamma_i)_{\tau_i}- \gamma_0}_{L^2_{\epsilon_i}} & \leq \epsilon_i \norm{(\gamma_i)_{\tau_i}- (\gamma_0)_{\tau_i}}_{L^2_{\epsilon_i}}+ \epsilon_i \norm{(\gamma_0)_{\tau_i}- \gamma_0}_{L^2_{\epsilon_i}}\\
& \leq \epsilon_i  \delta_2 + C_6 \epsilon^2_i \normltwo{( \gamma_0)_x} + C_6 \epsilon_i \left(\cosh\left((1-\norm b^2)\frac{\delta_1}{8} \right) +1\right) \leq \frac{\delta_1}{2}
\end{align*}
where we used that $\sup_{\R} \norm{\dot \gamma_0}_{\epsilon_i}\leq 2$ for all $i\geq I$ in line two for $I$ maybe bigger. For line four, we used \eqref{eq:equiv}, $( \gamma_0)_x \in L^2(\R, \R^{n-1})$ and estimated 
\[
\normltwo{((\gamma_0)_z)_\tau - \sgn} \leq 2\sqrt{2} \cosh((1-\norm b^2) \tau)/ \sqrt{1-\norm b^2}
\] 
with $\sgn$ the sign function. This results in a constant $C_6>0$ and the last inequality holds as soon as $i\geq I$ for $I$ maybe bigger.
Now applying Proposition \ref{prop:loc1}, we finally get
$
(\gamma_i)_{\tau_i}=\gamma_{\epsilon_i}
$
for all $i\geq I$. This proves the existence of an $\epsilon_0>0$ such that strong uniqueness holds for $0<\epsilon<\epsilon_0$. 
\end{proof}
%
%
%
\section{Global Uniqueness}\label{sec:five}

This section contains the steps (XIII) - (XV) from the introduction, each in a separate subsection. Namely, we get a priori estimates in Corollary \ref{cor:isol1} and prove the global to local result in Proposition \ref{prop:globloc}. This will finish the proof of uniqueness of the Main Theorem. $(A,b,h)$ will always satisfy Assumptions \ref{sta:trans1}.
%
%
%
%
%
\subsection{Conley Index Pair}\label{sec:isol}

In this section, we will construct a pair of sets $(N_\epsilon, L_\epsilon)$ with $L_\epsilon \subset N_\epsilon$ such that all solutions of \eqref{eq:eq4} are contained in $N_\epsilon \setminus L_\epsilon$. (Cf. Theorem \ref{thm:isol1}) The proof of the isolating statement will require an argument on the bigger set $N_\epsilon$. As a direct application of this theorem, we will get global uniqueness in subsections \ref{sec:globloc} and \ref{sec:fivethree}. One crucial ingredient in the proof of this theorem will be the energy-length inequality in Proposition \ref{prop:isol1}, which states that all solutions with bounded energy have bounded length.

\begin{Theorem}\label{thm:isol1}
Let $\nu \in (0,1/2)$. There is $K>0$ as in Proposition \ref{prop:isol1} and $\epsilon_0>0$ such that for $0<\epsilon<\epsilon_0$, there are compact sets $L_\epsilon \subset N_\epsilon \subset \R^n$ (as in Definition \ref{def:isol1}) with the following properties:
\begin{enumerate}[label=(\roman*),leftmargin=*]
\item every solution $\gamma:\R \to \R^n$ of the ODE in \eqref{eq:eq4} with $\lim_{t\to \infty} \gamma(t) = (0,-1)$ which leaves $N_\epsilon$ at time $T \in \R$ has $f_\epsilon(\gamma(T)) < -\frac 23$, where $f_\epsilon$ was defined in \eqref{eq:eq5},
\item every solution $\gamma:\R \to \R^n$ of the ODE in \eqref{eq:eq4} with $\lim_{t\to \pm \infty} \gamma(t) = (0,\pm 1)$ is contained in $N_\epsilon \setminus L_\epsilon$,
\item if $(x,z) \in N_\epsilon$, then $\norm{Ax+b(z^2-1)} \leq \sqrt 2 \epsilon^{-1}$ and $\abs{z} \leq K$,
\item if $(x,z) \in N_\epsilon \setminus L_\epsilon$, then $\norm{Ax+b(z^2-1)} \leq \sqrt 2 \epsilon^\nu$ and $\abs{z} \leq K$.
\end{enumerate}
\end{Theorem}

\begin{proof}
See page \pageref{proof:isol1}.
\end{proof}

\subsubsection{Preliminary definitions and results.}

We start things off by a change of coordinates centred around the limit solution $\gamma_0$. This can be done by taking part of the right hand side of \eqref{eq:eq4} as new variable 
\[
w(x,z) :=Ax+ b (z^2-1).
\] 
So we can rewrite \eqref{eq:eq4} as 
\[
\dot x = - w(x,z) /\epsilon + R_\epsilon^x(x,z), \text{ and } \dot z = -   \scal{b}{w(x,y)} + (1-\norm b^2)(1-  z^2) + R_\epsilon^z(x,z). 
\] 
This leads to the new equations for differentiable $w:\R \to \R^{n-1}$ and  $z:\R \to \R$
\begin{equation}\label{eq:isol1}
\begin{cases}
&\dot w = - {A_\epsilon(z) w}/\epsilon + \tilde R_\epsilon^w(w,z), \\
&\dot z = -   \scal b w + (1-\norm b^2)(1-  z^2) + \tilde R_\epsilon^z(w,z),\\	
& \lim_{t\to \pm \infty} (w(t), z(t)) = (0, \pm 1),
\end{cases}
\end{equation}
\vspace{-1em}\begin{align}
&\text{where } {\tilde R_\epsilon^w(w,z)} = AR_\epsilon^x(x(w,z),z) + 2  b z (1-\norm b) (1-  z^2) + 2   b z R_\epsilon^z(x(w,z),z) , \label{eq:isol2}\\
&  {\tilde R_\epsilon^z(w,z)} = R_\epsilon^z(x(w,z), z), \,  A_\epsilon(z) := A + 2 \epsilon z b b^\top   \text{ and } x(w,z) := A^{-1}   w +  A^{-1}   b (1- z^2). \notag
\end{align}
Then, since $A$ is invertible, $(x,z)$ is a solution for \eqref{eq:eq4} if and only if $(w(x,z),z)$ is a solution for \eqref{eq:isol1}. Furthermore, $w(\gamma_0)= 0$ by \eqref{eq:adia1}.  Now if we define the diffeomorphism $\Psi(x,z) = (w(x,z),z)$, then we have that $(\Psi_\ast f_\epsilon, \Psi_\ast g_\epsilon)$ is a gradient pair for \eqref{eq:isol1} where $(f_\epsilon, g_\epsilon)$ are as in \eqref{eq:eq5}. This means in particular that
\begin{align}\label{eq:isol3}
(\Psi_\ast f_\epsilon) (w,z) = \frac{\epsilon}{2}w^\top A^{-1} w +  \frac{\epsilon}{2} (b(1-  z^2))^\top A^{-1} (2 w+b(1-  z^2)) + \frac 13 z^3 - z
\end{align}
decreases along solutions of \eqref{eq:isol1}. For future use, we have $(\Psi_\ast f_\epsilon) (0,\pm 1) = \mp \frac{2}{3}$.

The advantage of this new set of variables is that the dynamics of $w$ has a hyperbolic term ${A_\epsilon(z) w}/\epsilon$ which is big compared to the terms in $\tilde R_\epsilon^w(w,z)$ as $\epsilon$ approaches zero. The sets $L_\epsilon$ and $N_\epsilon$ will be chosen to reflect the dominating hyperbolic part. 

\begin{figure}[htb]
\centering
\includegraphics[scale=0.9]{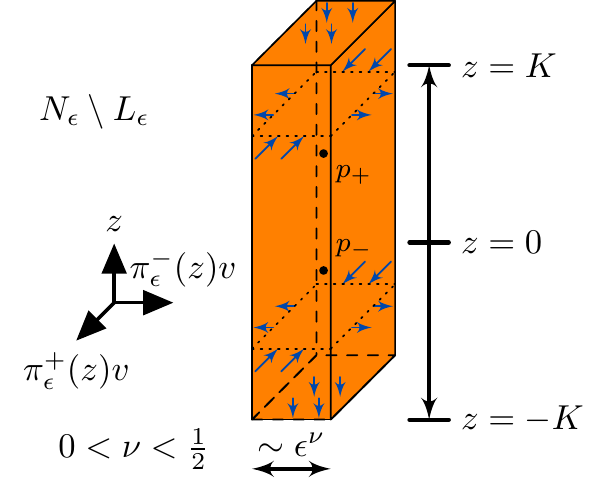}
\caption{The set $N_\epsilon\setminus L_\epsilon$ with direction of the vector field on its boundary. }\label{fig:isol1}
\end{figure}

\begin{figure}[htb]
\centering
\includegraphics[scale=0.9]{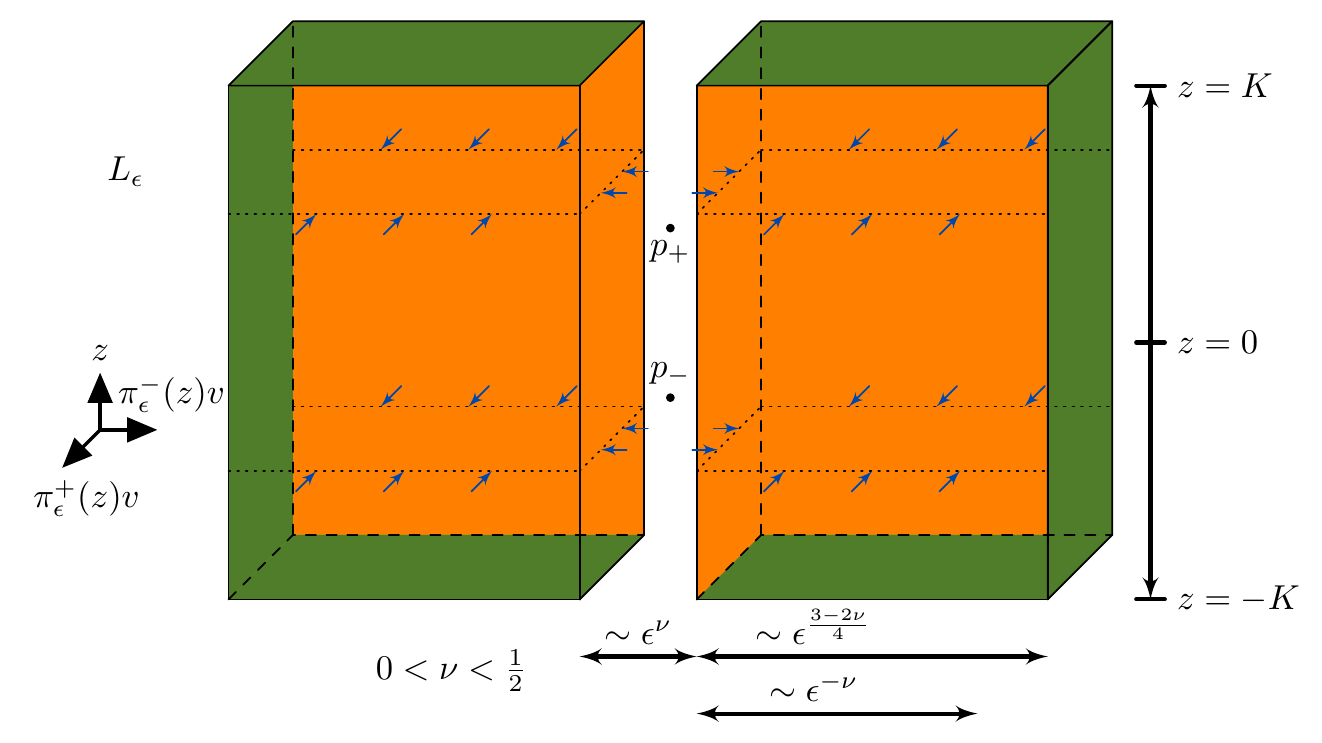}
\caption{The set $L_\epsilon$ with direction of the vector field on some boundary faces. Green faces indicate where $f_\epsilon$ will be too small, once we exit this set.}\label{fig:isol2}
\end{figure}

The matrix $A_\epsilon(z) := A + 2 \epsilon z  b  b^\top $ is still symmetric in $\R^{(n-1)\times (n-1)}$ and so is diagonalisable with real eigenvalues. In addition, for fixed $R>0$, there is $\epsilon_0>0$ such that $A_\epsilon(z)$ is invertible for all $0<\epsilon<\epsilon_0$ and $z\in [-R,R]$. 
From \eqref{eq:isol1}, we are inspired to split $w= \pi_\epsilon^{+}(z)w+\pi_\epsilon^{-}(z)w$ where
\begin{equation}\label{eq:isol5}
\begin{aligned}
\pi_\epsilon^{\pm}(z) \text{ are the orthogonal positive resp. negative eigenspace projection of $A_\epsilon(z)$.}
\end{aligned}
\end{equation}
%
%
This lets us define the pair of sets $(N_\epsilon , L_\epsilon)$ for Theorem \ref{thm:isol1}. 
\begin{Definition}\label{def:isol1}
Fix $0<\nu<\frac{1}{2}$ and take $K$ as in Proposition \ref{prop:isol1} below. 
\begin{equation*}
\begin{aligned}
N_\epsilon &:= \left\{(x,z) \in \R^n \, : \, \abs{z} \leq K,  \norm{\pi_\epsilon^+(z) w(x,z)} \leq \epsilon^\nu, \,\phantom{\epsilon^\nu \leq}\norm{\pi_\epsilon^-(z) w(x,z)} \leq \epsilon^{\frac{3-2\nu}{2}}\right\},\\
L_\epsilon &:= \left \{(x,z) \in \R^n \, : \, \abs{z} \leq K,   \norm{\pi_\epsilon^+(z) w(x,z)} \leq \epsilon^\nu, \,\epsilon^\nu \leq \norm{\pi_\epsilon^-(z) w(x,z)} \leq \epsilon^{\frac{3-2\nu}{2}}\right\}.
\end{aligned}
\end{equation*}
See also Figure \ref{fig:isol1} and \ref{fig:isol2}. 
\end{Definition}
%
%
\begin{Notation} \label{not:isol1}
For vectors $\gamma \in \R^n$, we use the following coordinate notation: $\gamma = (\gamma_x, \gamma_z)$ for its $x$- resp. $z$-component as always. We define $\gamma_w:=A\gamma_x+ b (\gamma_z^2-1)$. Hence a solution $\gamma:\R \to \R^n$ of \eqref{eq:eq4} will fulfil the following differential equations
\begin{equation}\label{eq:isol8}
\begin{aligned}
&\epsilon \dot \gamma_x = -A\gamma_x +b (1- \gamma_z^2) + \epsilon R_\epsilon^x(\gamma_x,\gamma_z), \; \dot \gamma_z = -   b^\top  {A\gamma_x} + 1-  {\gamma_z}^2 + R_\epsilon^z(\gamma_x,\gamma_z),\\
&\text{which is equivalent to }\\
&\epsilon \dot \gamma_w = \!-A_\epsilon(z) \gamma_w + \epsilon \tilde R_\epsilon^w(\gamma_w,\gamma_z), \dot \gamma_z = \!-  \scal b {\gamma_w} + (1-\norm b^2)(1- { \gamma_z}^2) + \tilde R_\epsilon^z(\gamma_w,\gamma_z).
\end{aligned}
\end{equation}
Also in this notation, we have $f_\epsilon(\gamma_x, \gamma_z) = (\Psi_\ast f_\epsilon)(\gamma_w, \gamma_z)$, where $f_\epsilon$ was defined in \eqref{eq:eq5}.
\end{Notation}
%
\begin{Remark}
$(N_\epsilon, L_\epsilon)$ can be made into a Conley index pair \cite{Con1, Sal3} for the invariant set 
\[
S_\epsilon := \left\{ \gamma(\R) \,:\, \gamma: \R \to \R^n, \, \gamma(\pm \infty) = (0,\pm 1),\, \gamma \text{ solves } \eqref{eq:eq4} \right\} \cup \{(0,\pm 1)\},\]
by adding the points $(x,z) \in \R^n$ with $z=-K,\,  \norm{\pi_\epsilon^\pm(z) w(x,z)} \leq \epsilon^\nu$ to the exit set $L_\epsilon$.
We note that since $p_\pm$ are the only critical points of $ f_\epsilon$ in $N_\epsilon \setminus L_\epsilon$, the set $S_\epsilon$ is the biggest invariant set that could possibly be contained in $N_\epsilon \setminus L_\epsilon$ and Theorem \ref{thm:isol1} (ii) states that this set is indeed contained therein. The proof of Theorem \ref{thm:isol1} will imitate some features of the Conley property, but omits its full proof.
\end{Remark}
We now establish the following technical lemma with useful estimates for later on.
\begin{Lemma}\label{lem:isol1}
Fix $R>0$. There is $\epsilon_0>0$ and $M>0$ such that for all $0<\epsilon<\epsilon_0$, $w \in \R^{(n-1)}$ with $\norm w \leq \sqrt2\epsilon^{-1}$, $z \in [-R,R]$, $(\hat w, \hat z) \in \R^{(n-1)}\times \R$,  we have
\begin{equation*}
\begin{array}{ll}
\textbf{(i)}\; \norm{ x(w,z) } \leq 2 \norm{A^{-1}} \epsilon^{-1},				&\textbf{(ii)} \,\epsilon \norm{\tilde R_\epsilon^w(w,z)}  + \abs{\tilde R_\epsilon^z(w,z)}\leq  M (1+ \norm w) \epsilon, \\[0.5em]
\textbf{(iii)}  \norm{h_{\epsilon^2,(\epsilon^2 x(w,z), \epsilon z)}} \leq M \epsilon, 
&\textbf{(iv)}\,  \begin{cases} \norm{(d \pi_\epsilon^\pm(z) \hat z) \hat w} \leq M\epsilon \norm{\hat z} \norm{\hat w}, \\[0.5em]
						\norm{\pi^\pm_\epsilon(z) - \pi_\epsilon^\pm(0)} \leq M \epsilon, \end{cases}
\end{array}
\end{equation*}
where $\tilde R_\epsilon^w, \tilde R_\epsilon^z, x(w,z)$ are defined in \eqref{eq:isol2}, $\pi_\epsilon^{\pm}(z)$ in \eqref{eq:isol5}, $h_{\epsilon^2}$ as in Assumption \ref{sta:trans1}.
\end{Lemma}
\begin{proof}
\textbf{(i)} is true by the definition in  \eqref{eq:isol2} for $\epsilon_0>0$ small enough.  \textbf{(iii)} follows similarly to \eqref{eq:unif5}. \textbf{(ii)} is a consequence of \textbf{(i)} and \textbf{(iii)} combined with the definitions in \eqref{eq:eq4} and \eqref{eq:isol2}. Let us prove \textbf{(iv)}. 
For $z\in [-R,R]$,
$
\| \epsilon 2zbb^\top \|  \leq 2R \epsilon.
$
So for $\epsilon_0>0$ small, we have that for $0<\epsilon<\epsilon_0$, $A_\epsilon(z)$ is still invertible. Indeed for $\epsilon_0>0$ sufficiently small, there are $0<\kappa<\mathcal K$, such that
\begin{equation}\label{eq:isol10}
\begin{split}
0<\kappa \norm{\hat w}= \frac{1}{2\norm {A^{-1}}} \norm{\hat w}\leq \frac{1}{\norm{A^{-1}}}\norm{\hat w} - 2R \epsilon  \norm{\hat w} \leq  \norm{ {A_\epsilon(z)  \hat w}}  \leq \mathcal K \norm{\hat w},
\end{split}
\end{equation}
for all $\hat w \in \R^{(n-1)}$.
As $A_\epsilon(z)$ is symmetric and fulfils \eqref{eq:isol10}, the positive eigenvalues are in the real segment $[\kappa, \mathcal K]$ and the negative eigenvalues are in the real segment $[-\mathcal K,-\kappa]$.

Thus, choose a simple loop $\gamma^+$ in $\C$ encircling the real segment $[\kappa, \mathcal K]$ and $\gamma^-$ encircling $[-\mathcal K,-\kappa]$. This lets us define the projections on the positive respectively negative eigenvalues of $A_\epsilon(z)$ by%
\footnote{These formulae can either be computed directly using the residue theorem of complex analysis or one uses functional calculus of a bounded linear operator from functional analysis which can be found in \cite[5.2.10]{DietmarFA}. }
$
\pi^\pm_\epsilon(v) \hat w:= \frac{1}{2\pi i} \int_{\gamma^\pm} \left( z \mathbbm 1- A_\epsilon(z) \right)^{-1} \hat w \,d z
$
for $\hat w \in \R^{(n-1)}$. So differentiating yields for $\hat z \in \R$ and $\hat w \in \R^{(n-1)}$, \,
\[
(d \pi_\epsilon^\pm(z) \hat z) \hat w = \frac{2\epsilon \hat z}{2 \pi i} \int_{\gamma^\pm} (z \mathbbm 1 - A_\epsilon(z))^{-1}   {b b^\top}  (z \mathbbm 1 - A_\epsilon(z))^{-1} \hat w  \, d z.
\]
Now define for $\mu \in \R$ the linear operator $B(\mu) =A+2\mu b b^\top$ and
\[
C:=\sup_{\hat w \in B_1(0),\, \abs{\mu} \leq R \epsilon_0, \, j=\pm}   \abs { \frac{1}{ \pi i} \int_{\gamma^j} (z \mathbbm 1 - B(\mu))^{-1} {b b^\top} (z \mathbbm 1 - B(\mu))^{-1} \hat w \, d z} <\infty.
\]
Thus, we get
$
\norm{(\pi_\epsilon^\pm(z) - \pi_\epsilon^\pm(0)) \hat w} = \norm{ \int_0^1(d \pi_\epsilon^\pm (tz) z) \hat w \, dt} \leq \epsilon C R \norm{\hat w}.
$ 
\end{proof}
\subsubsection{Energy-Length inequality}
A key ingredient in the proof of Theorem \ref{thm:isol1} is that trajectories with bounded energy have bounded length. We formulate this in the following result. We use Notation \ref{not:isol1}.
\begin{Proposition}[Length-Energy Inequality]\label{prop:isol1}
There is $\epsilon_0>0$ and $K>0$ such that any solution $\gamma$ of the ODE in \eqref{eq:isol8} for $0<\epsilon< \epsilon_0$ with $\lim_{t\to -\infty}\gamma(t) = (0,-1)$ and
\[
\abs{\gamma_z(T)} > K, \qquad \norm{\gamma_w(t)} \leq \sqrt 2 \epsilon^{-1}  \text{ for }t\leq T 
\]
for some $T \in \R$ has energy 
\[
E_\epsilon\big(\gamma|_{(-\infty, T]}\big) := \int_{-\infty}^T g_{\epsilon, \gamma} (\dot \gamma, \dot \gamma) > \frac{4}{3} = f_\epsilon(0,-1) -f_\epsilon(0,1).
\]
In formulae, $K = \sqrt{\frac{2-\norm b^2}{1-\norm b^2}} + \frac{4}{3} \sqrt {32}$.
\end{Proposition}
%
%
\begin{proof}
We will prove the converse. Up to shifting time, we may assume 
\[
\abs{\gamma_z(0)} \geq l:= \sqrt{\fracb +1}, \; \abs{\gamma_z(T)}=K>l \text{ and }\abs{\gamma_z(t)} \geq l \text{ on }[0,T] 
\]
with $E_\epsilon\big(\gamma|_{(-\infty, T]}\big)\leq \frac{4}{3}$ and $\norm{\gamma_w(t)} \leq \sqrt 2 \epsilon^{-1}$ for $t\leq T$. We get from \eqref{eq:isol8} that 
\[
\dot \gamma_z - \epsilon  \scal b {\dot \gamma _x} = (1-\norm b^2) (1- \gamma_z^2).
\]
So we can estimate that for $t \in [0,T]$, $\abs{\dot \gamma_z - \epsilon  \scal b {\dot \gamma_x}} \geq (1-\norm b^2) (l^2-1) = 1$. Now we can use $\norm b <1$ and apply Lemma \ref{lem:isol1} (iii) with $R=K$ in combination with equivalence of norms from Lemma \ref{lem:equiv1} to get on $[0,T]$
\begin{equation}\label{eq:isol7}
\begin{aligned}
1 \leq \abs{\dot \gamma_z - \epsilon  \scal b {\dot \gamma_x}}^2 \leq 2(\epsilon^2 \abs{\dot \gamma _x}^2 + \abs{\dot \gamma_z}^2) \leq 8 \abs{\dot \gamma}_{g_\epsilon,\gamma}^2,
\end{aligned}
\end{equation}
for $\epsilon_0$ small and where $\abs{\dot \gamma(t)}_{g_\epsilon,\gamma(t)} := \sqrt{g_{\epsilon,\gamma(t)} (\dot \gamma(t), \dot \gamma(t))}$.
We can integrate \eqref{eq:isol7} to get
\[
T \leq 8 \int_0^T \abs{\dot \gamma}_{g_\epsilon,\gamma}^2 \leq 8 E_\epsilon\big(\gamma|_{(-\infty, T]}\big).
\]
Now we can use Cauchy-Schwartz inequality to conclude that
\begin{align*}
K-l & \leq \abs{\gamma_z(T)-\gamma_z(0)} \leq \int_0^T \abs{\dot \gamma_z} \leq 2 \int_0^T \abs{\dot \gamma}_{g_\epsilon,\gamma} \leq 2 \sqrt{\int_0^T 1} \; \sqrt{ \int_0^T \abs{\dot \gamma}_{g_\epsilon,\gamma}^2}\\
& \leq 2 \sqrt T \sqrt{E_\epsilon\big(\gamma|_{(-\infty, T]}\big)} \leq \sqrt{32}\, E_\epsilon\big(\gamma|_{(-\infty, T]}\big)\leq \frac 4 3 \sqrt{32}. \qedhere
\end{align*}
\end{proof}
%
\subsubsection{Proof of Theorem \ref{thm:isol1}.}
\begin{proof}[Proof of Theorem \ref{thm:isol1}]\label{proof:isol1}
We have \textbf{(iii)} and \textbf{(iv)} by Definition \ref{def:isol1}. Take $R=K$ and let $M>0$ be the constant from Lemma \ref{lem:isol1} with this $R$. We start proving \textbf{(i)}. Let $\gamma$ be a solution of the ODE in \eqref{eq:isol8}.

We want to investigate the direction along the boundary of $L_\epsilon$ (inward or outward pointing) of the gradient vector field of which $\gamma$ is by definition an integral curve. Any boundary point $(x,z) \in \partial L_\epsilon$ fulfils at least one of the following conditions. 
\begin{multicols}{2}
\begin{enumerate}[label=(\alph*)]
\item $\abs{\piminus w(x,z)} = \epsilon^\nu,$
\item $\abs{\piminus w(x,z)} = \epsilon^{\frac{2\nu-3}{4}},$
\item $\abs{\piplus w(x,z)} = \epsilon^\nu,$
\item $z = \pm K$.
\end{enumerate}
\end{multicols}
\vspace{-\parskip} Any boundary point $(x,z) \in \partial N_\epsilon$ fulfils at least one of (b), (c) or (d). 

\textbf{Case 1:} Vector field direction for cases (a) and (b).
 
To prove that the vector field along this boundary points in for (a) and out for (b), we need to prove that $\rho_1: \R \to \R: t \mapsto \abs{\pi_\epsilon^-(\gamma_z(t)) \gamma_w(t)}^2$ is strictly increasing at $T\in \R$, whenever $\gamma(T)$ fulfils (a) or (b). Indeed, we estimate the derivative of $\rho_1$ from below.
\begin{align*} 
\rho_1' 	&= 2 \langle \pi_\epsilon^-(\gamma_z) \gamma_w, (d \pi_\epsilon^-(\gamma_z) \dot \gamma_z) \gamma_w\rangle  +  2 \langle \pi_\epsilon^-(\gamma_z) \gamma_w, \pi_\epsilon^-(\gamma_z) \dot \gamma_w \rangle \\
	& \geq  -2 \norm{\gamma_w} M \epsilon \abs{\dot \gamma_z} \norm{\gamma_w} - \frac{2}{\epsilon} \langle \pi_\epsilon^-(\gamma_z) \gamma_w, A_\epsilon(\gamma_z)(\pi_\epsilon^-(\gamma_z) \gamma_w)\rangle - \norm{\gamma_w} M (1+\norm{\gamma_w}) \\
	&\geq 2 \epsilon^{-1} \kappa \norm {\pi^-_\epsilon(\gamma_z)\gamma_w}^2 - 2\epsilon M C   \norm{\gamma_w}^2(1+ \norm{\gamma_w}) - M \norm{\gamma_w} (1+\norm{\gamma_w})      
\end{align*}
where we estimated quantities using (ii), (iv) in Lemma \ref{lem:isol1}, the equations for $\dot \gamma$ from \eqref{eq:isol8}, some constant $C>0$ and the definition of $\pi^-_\epsilon$ in \eqref{eq:isol5} with $\kappa$ being the lower bound for the absolute value of eigenvalues for $A_\epsilon(z)$ as in \eqref{eq:isol10}. The first term dominates the other two as long as $\norm{\pi^-_\epsilon(\gamma_z)\gamma_w} \leq \epsilon^{\frac{2\nu-3}{2}}<\epsilon^{-1}$, $\norm{\pi^+_\epsilon(\gamma_z)\gamma_w}\leq \epsilon^\nu$  and $\epsilon_0>0$ maybe smaller. Thus $\rho_1'(T)>0$.

\textbf{Case 2:} Vector field direction for case (c).  

To prove that the vector field along this boundary points in for (c), we need to prove that $\rho_2: \R \to \R: t \mapsto \abs{\pi_\epsilon^+(\gamma_z(t)) \gamma_w(t)}^2$ is strictly decreasing at $T\in \R$, whenever $\gamma(T)$ fulfils (c). Indeed, we estimate the derivative of $\rho_2$ from above.
\begin{align*}
\rho_2' 	&= 2 \langle \pi_\epsilon^+(\gamma_z) \gamma_w, (d \pi_\epsilon^+(\gamma_z) \dot \gamma_z) \gamma_w\rangle  \!+\!  2 \langle \pi_\epsilon^+(\gamma_z) \gamma_w, \pi_\epsilon^+(\gamma_z) \dot \gamma_w \rangle \leq 2 \norm{\pi^+_\epsilon(\gamma_z)\gamma_w} M \epsilon \abs{\dot \gamma_z} \norm{\gamma_w} \\
&- \frac{2}{\epsilon} \langle \pi_\epsilon^+(\gamma_z) \gamma_w, A_\epsilon(\gamma_z)(\pi_\epsilon^+(\gamma_z) \gamma_w)\rangle +  \norm{\pi^+_\epsilon(\gamma_z)\gamma_w} M (1+\norm{\gamma_w})\\
		&\leq -2 \epsilon^{-1} \kappa \norm{\pi^+_\epsilon(\gamma_z)\gamma_w}^2  +  \norm{\pi^+_\epsilon(\gamma_z)\gamma_w} \big(2M C \epsilon (1+\norm{\gamma_w}) \norm{\gamma_w} + M (1+\norm{\gamma_w})\big) 
\end{align*}
where we used the same sort of estimates as in Case 1. Now plug in that $\norm{\pi^+_\epsilon(\gamma_z)\gamma_w}= \epsilon^\nu$ and $\norm{\pi^-_\epsilon(\gamma_z)\gamma_w} \leq \epsilon^{\frac{2\nu-3}{4}}$, to get
\[
\rho_2' \leq -2 \kappa \epsilon^{2\nu-1} + 4 M C \epsilon^{2\nu - \frac{1}{2}} + 2 M \epsilon^{6\frac{ \nu -3}{4}}. 
\]
We have that $2\nu-1 < \frac{6 \nu -3}{4} \iff \nu <\frac{1}{2}$ and so the first term dominates the other two as long as $\epsilon_0>0$ is maybe smaller. Thus $\rho_2'(T)<0$.

\textbf{Case 3:} Vector field direction for case (d) cannot be controlled for large $\gamma_w$. This is where the whole necessity of having Proposition \ref{prop:isol1} comes in. 

Now assume in addition that $\lim_{t\to-\infty}\gamma(t) = (0,-1) \in N_\epsilon$. Since $(0,-1)$ is an interior point, $\gamma(T) \notin N_\epsilon$ at time $T \in \R$ implies that there is $\tau<T$ such that $\gamma((-\infty, \tau]) \subset N_\epsilon$ and $\gamma(\tau) \in \partial N_\epsilon$, i.e. $\gamma(\tau)$ fulfils (b), (c) or (d). 

First we assume condition (b) holds. We need to prove that $f_\epsilon( \gamma(\tau)) < -\frac{2}{3}$. For this, we start making estimates on 
$\gamma^{\pm}_w(\tau):= \pi_\epsilon^\pm(0)  \gamma_w(\tau)$. 
\begin{align*}
\abs{\gamma_w(\tau)} &\leq 2 e^{\frac{2\nu-3}{4}}, \\
\abs{\gamma^{-}_w(\tau)} &\geq \abs{\pi_\epsilon^-(\gamma_z(\tau)) \gamma_w(\tau)} +\abs{\gamma^{-}_w(\tau) - \pi_\epsilon^-(\gamma_z(\tau)) \gamma_w(\tau)} \geq \epsilon^{\frac{2\nu-3}{4}} - \epsilon M \abs{\gamma_w(\tau)}  \\
&\geq \epsilon^{\frac{2\nu-3}{4}} (1 - 2\epsilon M)  \geq \frac{1}{2}  \epsilon^{\frac{2\nu-3}{4}}, \\
\abs{\gamma^{+}_w(\tau)} &\leq \abs{\pi_\epsilon^+(\gamma_z(\tau)) \gamma_w(\tau)} +\abs{\gamma^{+}_w(\tau) - \pi_\epsilon^+(\gamma_z(\tau)) \gamma_w(\tau)} \leq \epsilon^\nu +2 M \epsilon^{1+\frac{2\nu-3}{4}} \leq 2 \epsilon^\nu,
\end{align*}
where we used (iv) from Lemma \ref{lem:isol1}, $\abs{\pi_\epsilon^-(\gamma_z(\tau)) \gamma_w(\tau)} = \epsilon^{\frac{2\nu-3}{4}}$, $\abs{\pi_\epsilon^+(\gamma_z(\tau)) \gamma_w(\tau)} \leq \epsilon^{\nu}$, $\nu < \frac{2\nu+1}{4} \iff \nu < \frac{1}{2}$ and $\epsilon_0>0$ maybe smaller. We recall the expression for $\Psi_\ast f_\epsilon$ in \eqref{eq:isol3} and $A$ is a symmetric, invertible matrix whose eigenvalues are in $[-\mathcal K, -\kappa] \cup [\kappa, \mathcal K]$ as in \eqref{eq:isol10}. By definitions, $\gamma^{\pm}_w(\tau)$ is the positive/negative eigenspace projection with respect to $A$. So we can use the estimates for $\gamma^{\pm}_w(\tau) $ above to get
\begin{align*}
\Psi_\ast f_\epsilon(\gamma_w(\tau),\gamma_z(\tau)) &\leq \frac{\epsilon}{2} \gamma^{-}_w(\tau)^\top  A^{-1} \gamma^{-}_w(\tau) +  \frac{\epsilon}{2} \gamma^{+}_w(\tau)^\top A^{-1} \gamma^{+}_w(\tau) + C (1+ \epsilon^{\frac{2 \nu +1}{4}}) \\
&\leq -\frac{\epsilon^{1+\frac{2\nu-3}{2}}}{8 \mathcal K} + \frac{2 \epsilon^{1+2\nu}}{\kappa}+ C(1+\epsilon^{\frac{2 \nu +1}{4}})
\end{align*}
where we used $\abs{\gamma_z} \leq K$ to get $C>0$ independent of $\epsilon$ and $1+\frac{2\nu-3}{2} <0 \iff \nu <\frac{1}{2}$. Hence, the first term dominates the other two as long as $\epsilon_0>0$ maybe even smaller. So by the equation in the line below \eqref{eq:isol8},  $f_\epsilon(\gamma(\tau)) < -\frac{2}{3} $.

Due to the direction of the vector field in case (c) (Case 2 above), this possibility is excluded.
If case (d) happens, then $\norm{\gamma_w(t)} \leq \epsilon^{-1} $ for $t \leq \tau$ and $|\gamma_z(\tau)|\geq K$. So the Length-Energy Inequality \ref{prop:isol1} gives%
\footnote{Fact on gradient flows: $\dot \gamma = -\nabla_g f \Rightarrow \int_{t_0}^{t_1} g_\gamma(\dot \gamma, \dot \gamma) = \int_{t_0}^{t_1} g_\gamma ( -\nabla_{g} f(\gamma), \dot \gamma) \, d t = -f(\gamma(t_0)) + f(\gamma(t_1))$ for $t_0<t_1$.}
\[
\frac2 3 - f_\epsilon(\gamma(\tau)) = \lim_{t\to-\infty} f_\epsilon(\gamma(t)) -f_\epsilon(\gamma(\tau)) =  \int_{-\infty}^{\tau}g_{\epsilon, \gamma} (\dot \gamma, \dot \gamma) > \frac{4}{3}
\]  
which implies $f_\epsilon(\gamma(\tau)) <-\frac{2}{3}$. This proves \textbf{(i)}.

We now prove how \textbf{(i)} and Case 1 imply \textbf{(ii)}. So assume that $\lim_{t \to \pm \infty} \gamma(t) = (0,\pm 1)$. As $f_\epsilon(0,\pm 1)=\mp \frac2 3$, we have $\abs{f_\epsilon(\gamma(t))} \leq \frac{2}{3}$ for $t\in \R$. By (i), this implies that $\gamma(\R) \subset N_\epsilon$. Now assume that there is $T \in \R$ such that $\gamma(T) \in L_\epsilon$. By Case 1, this implies that $\abs{\pi_\epsilon^-(\gamma_z(t)) \gamma_w(t)} \geq \epsilon^\nu$ for all $t \geq T$. This is a contradiction since $\lim_{t\to \infty} \abs{\pi_\epsilon^-(\gamma_z(t)) \gamma_w(t)} = 0$. So $\gamma(\R) \subset N_\epsilon \setminus L_\epsilon$.     
\end{proof}
%
%
\subsection{A Priori Estimates} \label{sec:fivethree}

\begin{Corollary}\label{cor:isol1} Fix $\nu\in (0, 1/2)$.
There is $R>1$ and $\epsilon_0>0$, such that for all $0<\epsilon<\epsilon_0$, a solution $\gamma$ of \eqref{eq:eq4} has the a priori bound
\begin{equation}
\normlinftyy{A \gamma_x + b(\gamma_z^2-1)} \leq \epsilon^\nu, \quad \text{and} \quad \normlinftyy{\gamma_z} \leq R.
\end{equation}  
\end{Corollary}
\begin{proof}
Take $\nu<\nu'<1/2$. By Theorem \ref{thm:isol1} (ii), we have for $\epsilon_0$ small, that a solution $\gamma$ of \eqref{eq:eq4} for $0<\epsilon<\epsilon_0$ has the a priori bound
\[
\normlinftyy{A \gamma_x + b(\gamma_z^2-1)} \leq \sqrt 2 \epsilon^{\nu'} \leq \epsilon^\nu, \quad  \normlinftyy{\gamma_z} \leq K=:R. \qedhere
\]
\end{proof}

\begin{Remark}
In particular, the estimates in Corollary \ref{cor:isol1} hold for the solution $\gamma_\epsilon$ from the Existence Theorem \ref{thm:existence}. Analysing the proof of the Strong Uniqueness Theorem \ref{thm:loc1}, we get for $\gamma_0$ the limit solution in \eqref{eq:adia3} that
\[
\normlinfty{\gamma_\epsilon - \gamma_0} = \epsilon \normlinftyy{(\gamma-\gamma_0)_x} + \normlinftyy{(\gamma-\gamma_0)_z} \xrightarrow{\epsilon\to 0} 0.
\]
This means that while the $z$-component converges uniformly to zero, the $x$-component may only converge after multiplying it by $\epsilon>0$. 
This is why we call the limit solution $\gamma_0$ the adiabatic limit.
\end{Remark}

\begin{Corollary}[Uniqueness in Theorem \ref{thm:goal2}]\label{thm:uniq2}
There exists $\epsilon_0>0$ such that for all $0<\epsilon<\epsilon_0$, any solution $\gamma$ of \eqref{eq:eq4} is up to time shift equal to $\gamma_\epsilon$ where $\gamma_\epsilon$ is the solution from Theorem \ref{thm:existence}.
\end{Corollary}
\begin{proof}
By Corollary \ref{cor:isol1}, there is $\epsilon_0>0$ such that for all $0<\epsilon<\epsilon_0$, any solution $\gamma$ of \eqref{eq:eq4} fulfils the condition in the Strong Uniqueness Theorem \ref{thm:loc1}. So uniqueness follows for $\epsilon_0>0$ maybe smaller. This finishes the proof of uniqueness and of Theorem \ref{thm:goal2}.
\end{proof}

%
%
\subsection{Global To Local}\label{sec:globloc}

We can now use Theorem \ref{thm:isol1} to reduce the global problem on the manifold to the local theorem in the chart. This reduces the Main Theorem to its local version, Theorem~\ref{thm:goal2}.
\begin{Proposition}[Global to Local]\label{prop:globloc}
Let $(F_\lambda, G_\lambda)$ be as in the Main Theorem from the Introduction
and let $\varphi_\lambda:U \subset M \to \R^m$ be the charts from Theorem \ref{thm:normal2} for $0<\lambda<\epsilon_1^2$. 

Then there is $0<\epsilon_0<\epsilon_1$ such that for all $0<\lambda<\epsilon_0^2$, every gradient trajectory $\Gamma_\lambda:\R\to M$ solving $\dot \Gamma_\lambda~=~-\nabla_{G_\lambda} F_\lambda~\circ~ \Gamma_\lambda$ and connecting the critical points $p_-(\lambda)$ to $p_+(\lambda)$ is contained in the chart i.e. 
\[
\Gamma_\lambda (\R) \subset U.
\] 
Furthermore, there is a bijection (up to time-shift)%
\footnote{By this we mean that $\Gamma_1$ differs from $\Gamma_2$ by time-shift exactly if their corresponding $\gamma_1$ and $\gamma_2$ also do.}
 between such gradient trajectories $\Gamma_\lambda$ and solutions $\gamma_{\epsilon}$ of equation \eqref{eq:eq4} for $\epsilon=\sqrt \lambda$. Transversality is preserved under this bijection.
\end{Proposition}

\begin{proof}[Proof that Proposition \ref{prop:globloc} and Theorem \ref{thm:goal2} imply the Main Theorem] \label{proof:goal1}
By the statement of Proposition \ref{prop:globloc}, there is $\epsilon_0>0$ such that for all $0<\lambda<\epsilon_0^2$, we have a bijection (up to time-shift) between
\begin{itemize}
\item gradient trajectories $\Gamma_\lambda:\R\to M$ solving $\dot \Gamma_\lambda~=~-\nabla_{G_\lambda} F_\lambda~\circ~ \Gamma_\lambda$ and connecting the critical points $p_-(\lambda)$ to $p_+(\lambda)$, and
\item solutions $\gamma_\epsilon$ of \eqref{eq:eq4} with $\epsilon=\sqrt{\lambda}$.
\end{itemize}
By Theorem \ref{thm:goal2}, there exists a unique solution $\gamma_\epsilon$ of \eqref{eq:eq4} up to time-shift. Since transversality is preserved under the bijection, this proves the Main Theorem. 
\end{proof}

\begin{proof}
We use $\epsilon= \sqrt \lambda$. We recall from \eqref{eq:trans1} or Theorem \ref{thm:normal2} that there are
 $\mathfrak c, \rho, \epsilon_1>0$, a family of charts $\varphi_{\epsilon^2}: U\subset M \to \R^n$ with $\varphi_0(p_0)=0$ and a family of affine invertible maps $\chi_{{\epsilon^2}}:\R \to \R$ such that
\begin{equation}\label{eq:globloc1}
\begin{aligned}
\chi_{{\epsilon^2}} \circ F_{{\epsilon^2}} \circ \varphi_{{\epsilon^2}}^{-1}(x,z) 	&=\frac 1 {\mathfrak c}\bigg(\frac{1}{2} x^\top A x + \frac{1}{3} z^3 - {\epsilon^2} z\bigg),\\
((\varphi_{{\epsilon^2}})_\ast G_{\epsilon^2})_{(x,z)}^{-1} 					&= \mathfrak c\bigg(\begin{pmatrix} \mathbbm 1 & b \\   b^\top & 1 \end{pmatrix} + h_{{\epsilon^2},(x,z)}\bigg), 
\end{aligned}
\quad \text{ and }\quad B_{\sqrt2 \rho}(0) \subset \varphi_{\epsilon^2}(U),
\end{equation}
for all $(x,z)\in (\R^{n-1}\times \R)\cap \varphi_{\epsilon^2}(U)$ and all $0<\epsilon <\epsilon_1$. Here $(A,b,h)$ are as in Assumptions~\ref{sta:trans1}. Now as in the proof of Lemma \ref{lem:trans1}, we can further define $\varphi_\epsilon^{(3)}(x,z)~=~( x/ \epsilon^2, z/ \epsilon)$ and $\chi_\epsilon^{(3)}(w) = z/\epsilon$. We get for $0< \epsilon < \epsilon_1$,  $\varphi_\epsilon^{(4)}:= \varphi_\epsilon^{(3)} \circ \varphi_{\epsilon^2}$, and $\chi_\epsilon^{(4)}:=\chi_\epsilon^{(3)} \circ \chi_{\epsilon^2}$ that
\[
(\chi_\epsilon^{(4)} ) \circ F_{\epsilon^2} \circ (\varphi_\epsilon^{(4)} )^{-1}= \frac{ \epsilon^2}{\mathfrak c} f_\epsilon \text{ and }((\varphi_\epsilon^{(4)})_\ast G_{\epsilon^2})^{-1}= \frac {\mathfrak c} {\epsilon^2} g_\epsilon^{-1},
\]
where $(f_\epsilon, g_\epsilon)$ is the gradient pair in \eqref{eq:eq5}  for \eqref{eq:eq4}. By \eqref{eq:globloc1}, we have that
\[
V_\epsilon := \{(x,z) \in \R^{n-1} \times \R \, : \,  \norm x < \rho \epsilon^{-2}, \abs z < \rho \epsilon^{-1}\} \subset \varphi_{\epsilon^2}^{(3)}(B_{\sqrt 2\rho}(0)) = (\varphi_\epsilon^{(4)} )(U). 
\]
Now picking $0<\epsilon_0 < \epsilon_1$ small, we get by Theorem \ref{thm:isol1} (iii) and Lemma \ref{lem:isol1} (i), that 
\[
N_\epsilon \subset  \left\{(x,z) \in \R^{n-1} \times \R \, : \,  \norm x \leq \sqrt 2 \norm{A^{-1}} \epsilon^{-1}, \abs z \leq K \right\} \subsetneq  V_\epsilon \subset (\varphi_\epsilon^{(4)} )(U).
\]
Fix $0<{\epsilon}<\epsilon_0$. Let $\Gamma$ be a gradient trajectory solving $\dot \Gamma~=~-\nabla_{G_{\epsilon^2}} F_{\epsilon^2}~\circ~ \Gamma$ and having $\lim_{t \to -\infty} \Gamma(t) = p_-({\epsilon^2})$. Assume there is $T\in \R$ such that $\Gamma(T) \notin U$. Thus define 
\[
\gamma: (-\infty, T'] \to \R^n\text{ given by }\gamma(t) :=(\varphi_\epsilon^{(4)} ) \circ \Gamma_{\epsilon^2} ( a_{\epsilon^2} t),
\] 
for $t \leq T'$ and where $(\chi_\epsilon^{(4)} )(w) =: a_{\epsilon} w + b_\epsilon, \, a_{\epsilon} \neq 0$. Here $T'$ is such that $\gamma(T') \notin N_\epsilon$ and $\Gamma(a_{\epsilon^2} t) \in U$ for all $t \leq T'$, which exists by assumption. Then $\gamma$ solves the ODE in \eqref{eq:eq4} by construction and by Theorem \ref{thm:isol1} (i), we get 
\[
(\chi_\epsilon^{(4)} ) \circ F_{\epsilon^2}(\Gamma(a_{\epsilon^2} T')) =\frac{ \epsilon^2}{\mathfrak c} f_\epsilon(\gamma(T')) < -\frac{2}{3} \frac{ \epsilon^2}{\mathfrak c} = (\chi_\epsilon^{(4)} ) \circ F_{\epsilon^2}(p_+(\epsilon^2)) 
\] 
Hence $F_{\epsilon^2}(\Gamma(a_{\epsilon^2} T')) < F_{\epsilon^2}(p_+(\epsilon^2))$ and so $\lim_{t \to \infty} \Gamma(t) \neq p_+(\epsilon^2)$. 

This proves that any gradient trajectory $\Gamma$ connecting the critical points $p_-(\epsilon^2)$ to $p_+(\epsilon^2)$ is contained in $U$. Hence there is a bijection $\Gamma \mapsto \gamma_\Gamma$ defined by $\gamma_\Gamma(t) :=(\varphi_\epsilon^{(4)} ) \circ \Gamma_{\epsilon^2} ( a_{\epsilon^2} t)$. Then $\Gamma$ is a connecting gradient trajectory precisely if $\gamma_\Gamma$ is a solution to \eqref{eq:eq4}. We also see from that formula that $\Gamma$ is transverse exactly if $\gamma_\Gamma$ is. Furthermore, $\Gamma_1$ and $\Gamma_2$ differ by a time-shift if and only if $\gamma_{\Gamma_1}$ and $\gamma_{\Gamma_2}$ differ by time-shift.  
\end{proof}

%
%
%
\appendix
\section{Birth-Death Critical Points And Their Local Form}\label{sec:app}

In this appendix, we define birth-death critical points and introduce their normal form in local charts.
%
%
\begin{Definition}\label{def:normal1}
Given a function $F: M \to \R$ on a manifold, we say that $p_0$ is an \textbf{embryonic critical point}, if in any chart $\varphi:U \subset M \to \Omega \subset \R^n$ with $\varphi(p_0)=0$, we have 
\begin{itemize}
\item $0 \in \Crit(F \circ \varphi^{-1})$,
\item $\dim_{\R} \ker \operatorname{Hess} (F\circ \varphi^{-1})(0)=1$, say $\R v = \ker \operatorname{Hess} (F\circ \varphi^{-1})$,
\item $\partial_v^3 (F \circ \varphi^{-1})(0) \neq 0$,
\end{itemize} 
where $\operatorname{Hess} (F\circ \varphi^{-1})(0)=[\partial_{x_i} \partial_{x_j} (F\circ \varphi^{-1})(0)]_{ij}$ is the Hessian matrix.
A generic smooth family $\R \times M \to \R: (\lambda, p) \mapsto F(\lambda,p):=F_\lambda(p)$ containing $F$ as $F_0$, will have the condition
\begin{equation}\label{eq:normal1}
d(\partial_\lambda F|_{\lambda=0} \circ \varphi_0^{-1})(0)v \neq 0,  
\end{equation}
fulfilled. 
We call $p_0$ a \textbf{birth-death critical point} for the family $(F_\lambda)_{\lambda\in \R}$ at $\lambda=0$. 
\end{Definition}
Some authors call birth-death critical points an $A_2$ singularity \cite{Igu1} to highlight the similarity to the complex case.

We need the following deep result of H. Whitney.
%
%
\begin{Theorem}[\textbf{Whitney}]\label{thm:normal1}\cite{Cie1, Mar1, Whi1}
Let $F_\lambda:M\to \R$ be a family with a birth-death critical point $p_0$ at $\lambda=0$ as in Definition \ref{def:normal1}. Then there is $\rho, \epsilon_0>0$, a family of charts $\varphi_\lambda: U\subset M \to \R^n$ with $\varphi_0(p_0)=0$ and a family of affine invertible maps $\chi_{\lambda}:\R \to \R$ such that $B_{\sqrt2 \rho}(0) \subset \varphi_\lambda(U)$ and
\[
\chi_{\lambda} \circ F_{\lambda} \circ \varphi_{\lambda}^{-1}(x,z)=\frac{1}{2} \sum_{i=1}^{n-1} d_i x_i^2 + \frac{1}{3} z^3 \pm \lambda z 
\]
for all $(x,z) \in (\R^{n-1}\times \R)\cap \varphi_\lambda(U)$ and all $\lambda \in \R$ with $\abs{\lambda} \leq \epsilon_0^2$. 
Here $d_i \in \{\pm 1\}$ and 
\[
\pm = - \sgn \left( \partial_v^3 (F_0 \circ \varphi_0^{-1})(0) \cdot d(\partial_\lambda F|_{\lambda=0} \circ \varphi_0^{-1})(0)v \right).
\]
\end{Theorem}
%
%
\begin{Remark}[Origin Of The Terminology]
Assume that $d(\partial_\lambda F|_{\lambda=0} \circ \varphi_0^{-1})(0)v$ and $\partial_v^3 (F_0 \circ \varphi^{-1})(0)$ have opposite signs. By looking at Theorem \ref{thm:normal1}, one sees that the critical points $\varphi_\lambda^{-1} (x,z)$ fulfil the equations $2d_i x_i =0$ for $i=1,\ldots, n-1$ and $z^2-\lambda=0$. Therefore, for $\lambda > 0$, for $\abs{\lambda}$ small, there are always two Morse critical points $p_\pm(\lambda) := \varphi_\lambda^{-1} (0,\pm \sqrt \lambda)$ of $F_\lambda$ near $p_0$ which merge for $\lambda=0$ and then disappear for $\lambda< 0$. This explains the terminology birth side ($\lambda >  0$), embryonic ($\lambda=0$) and death side ($\lambda < 0$). 
\end{Remark}
%
%
Our main object of study in this paper will be a gradient pair $(F_\lambda, G_\lambda)_{\lambda\in \R}$ where $(F_\lambda)_{\lambda\in \R}$ is a family with a birth-death critical point. These have the following joint normal form. 
\begin{Theorem}\label{thm:normal2}
Let $(G_\lambda)_{\lambda\in \R}$ be a smooth family of Riemannian manifold on $M$ and let $(F_\lambda)_{\lambda \in \R}$ be a smooth family of functions with a birth-death critical point $p_0$ at $\lambda=0$ as in Definition \ref{def:normal1}. Then there is $\mathfrak c, \rho, \epsilon_0>0$, a family of charts $\varphi_\lambda: U\subset M \to \R^n$ with $\varphi_0(p_0)=0$ and a family of affine invertible maps $\chi_{\lambda}:\R \to \R$ such that%
\footnote{Here $G^{-1}$ stands for the inverse of the matrix associated with the metric.}
\begin{equation*}
\begin{aligned}
\chi_{\lambda} \circ F_{\lambda} \circ \varphi_{\lambda}^{-1}(x,z) 	&=\frac 1 {\mathfrak c}\bigg(\frac{1}{2} x^\top A x + \frac{1}{3} z^3 - \lambda z\bigg),\\
((\varphi_{\lambda})_\ast G_\lambda)_{(x,z)}^{-1} 					&= \mathfrak c\bigg(\begin{pmatrix} \mathbbm 1 & b \\   b^\top & 1 \end{pmatrix} + h_{\lambda,(x,z)}\bigg), 
\end{aligned}
\quad \text{ and }\quad B_{\sqrt2 \rho}(0) \subset \varphi_\lambda(U),
\end{equation*}
for all $(x,z)\in (\R^{n-1}\times \R)\cap \varphi_\lambda(U)$ and all $\lambda \in \R$ with $\abs{\lambda} \leq \epsilon_0^2$. Here $A \in \R^{(n-1)\times (n-1)}$ is a symmetric invertible matrix, $\mathbbm 1 \in \R^{(n-1)\times (n-1)}$ is the identity matrix, $b \in \R^{n-1}$ is a column vector with $\norm b <1$ and $h_\lambda:  \varphi_\lambda(U) \to \R^{n\times n}$ with $h_{0,(0,0)} = 0$ and $h_\lambda$ symmetric.  Furthermore, 
\[
\pm = - \sgn \left( \partial_v^3 (F_0 \circ \varphi_0^{-1})(0) \cdot d(\partial_\lambda F|_{\lambda=0} \circ \varphi_0^{-1})(0)v \right).
\]
\end{Theorem}
%
%
\begin{proof}
By Theorem \ref{thm:normal1}, we have that there is $\chi_\lambda^{(1)}$ and $\varphi_\lambda^{(1)}$ such that
\[
\chi_{\lambda}^{(1)} \circ F_{\lambda} \circ (\varphi_{\lambda}^{(1)})^{-1}(x,z)=\frac{1}{2} x^\top D x + \frac{1}{3} z^3 - \lambda z,
\]
with $D = \operatorname{diag}(d_1, \ldots, d_{n-1})$. We also get a positive definite symmetric matrix
\[
((\varphi_{0}^{(1)})_\ast G_0)_{(0,0)}^{-1} =: \begin{pmatrix} P & q \\   q^\top & \mathfrak c \end{pmatrix},
\]
which implies that $P \in \R^{(n-1)\times(n-1)}$ is a symmetric positive definite matrix, $q \in \R^{n-1}$ a column vector and $\mathfrak c>0$ such that $\sqrt \mathfrak c\norm{P^{-1/2}q}<1$.

This allows us to define $\varphi^{(2)}: \R^n \to \R^n$ given by $\varphi^{(2)}(x,z)= (\sqrt \mathfrak c P^{-1/2} x, z)$ and $\chi^{(2)}~:~\R \to \R$ given by $\chi^{(2)}(w) = \frac 1 {\mathfrak c} w$. We get
\begin{align*}
(\chi^{(2)} \circ \chi_{\lambda}^{(1)}) \circ F_{\lambda} \circ (\varphi^{(2)} \circ \varphi_{\lambda}^{(1)})^{-1}(x,z) 	&=\frac 1 {\mathfrak c} \left( \frac{1}{2} \sum_{i=1}^{n-1} x^\top \left( \frac{1}{\mathfrak c} P^{1/2} D P^{1/2} \right) x + \frac{1}{3} z^3 - \lambda z \right),\\[0.5em]
((\varphi^{(2)} \circ \varphi_{0}^{(1)})_\ast G_0)_{(0,0)}^{-1} 											&= \mathfrak c  \begin{pmatrix} \mathbbm 1 & \sqrt \mathfrak c P^{-1/2} q \\   \left(\sqrt \mathfrak c P^{-1/2} q\right)^\top & 1 \end{pmatrix}.
\end{align*}
Therefore, setting 
\[
A := \frac{1}{\mathfrak c} P^{1/2} D P^{1/2}, \; b := \sqrt \mathfrak c P^{-1/2} q \text{ and }h_{\lambda,(x,z)} :=  \frac{1}{\mathfrak c} ((\varphi^{(2)} \circ \varphi_{\lambda}^{(1)})_\ast G_\lambda)_{(x,z)}^{-1} - \begin{pmatrix} \mathbbm 1 & b \\   b^\top & 1 \end{pmatrix}, 
\]
finishes the proof. That the matrix $\begin{pmatrix} \mathbbm 1 & b \\   b^\top & 1 \end{pmatrix}$ is still positive definite is equivalent to $\norm b^2 = b^\top b < 1$.\qedhere

\end{proof}
\bibliographystyle{abbrv}
\bibliography{Library}
\end{document}